\documentclass[Threeside,11pt]{article}
\usepackage{etoolbox}
\makeatletter
\patchcmd{\thebibliography}
{\list}
{\small          
	\setstretch{1.2}
	\list}
{}{}
\makeatother
\usepackage{setspace}
\usepackage{tikz}
\usepackage{float}
\definecolor{myPurple}{RGB}{109,90,207}
\definecolor{myGreen}{RGB}{4,130,67}
\definecolor{myGold}{RGB}{253,177,71}
\definecolor{myBurgundy}{RGB}{63,1,44}
\definecolor{myTeal}{RGB}{8,119,127}
\definecolor{myCream}{RGB}{255,216,177}
\definecolor{myOrange}{RGB}{225,119,1}
\usetikzlibrary{arrows.meta,bending,calc}
\tikzset{
	ball3Donly/.style={
		circle, minimum size=1.4cm,
		shading=ball,
		ball color=myPurple,
		draw=myPurple
	}
}
\pgfmathsetmacro{\wNear}{2.0pt}
\pgfmathsetmacro{\wMid}{1.3pt}
\pgfmathsetmacro{\wFar}{0.7pt}
\pgfmathsetmacro{\bend}{6mm}
\usepackage[latin1]{inputenc}
\usepackage[english]{babel}
\usepackage{amsmath}
\usepackage{amsthm}
\IfFileExists{algorithm.sty}{\usepackage{algorithm}}{}
\usepackage[nosort,nocompress]{cite}
\usepackage{amsfonts}
\usepackage{amssymb}
\usepackage{mathrsfs}
\usepackage{graphicx}
\usepackage{colortbl,dcolumn}
\usepackage{marvosym}
\IfFileExists{ifsym.sty}{\usepackage{ifsym}}{}
\usepackage{paralist}
\usepackage{xcolor}
\usepackage{subcaption}
\usepackage{adjustbox}
\usepackage{enumitem}
\usepackage{booktabs,tabularx}
\usepackage[colorlinks,citecolor=blue,
linkcolor=blue,hypertexnames=false]{hyperref}
\allowdisplaybreaks
\IfFileExists{mathx10.tfm}{
\DeclareFontFamily{U}{mathx}{\hyphenchar\font45}
\DeclareFontShape{U}{mathx}{m}{n}{
	<5> <6> <7> <8> <9> <10>
	<10.95> <12> <14.4> <17.28> <20.74> <24.88>
	mathx10
}{}
\DeclareSymbolFont{mathx}{U}{mathx}{m}{n}
\DeclareMathAccent{\widecheck}{0}{mathx}{"71}
}{\let\widecheck\check}
\usepackage{graphicx}
\usepackage{caption}
\usepackage{amsthm}

\newtheorem{lemma}{\bf Lemma}[section]

\newtheorem{proposition}{\bf Proposition}[section]
\newtheorem{corollary}{\bf Corollary}[section]
\newtheorem{definition}{\bf Definition}[section]
\newtheorem{remark}{\bf Remark}[section]
\newtheorem{example}{\bf Example}
\numberwithin{equation}{section}

\allowdisplaybreaks
\newtheorem{theorem}{Theorem}[section]

\theoremstyle{remark}

\newcommand{\Z}{\mathbb Z}\newcommand{\R}{\mathbb R}
\newcommand{\C}{\mathbb C}\newcommand{\T}{\mathbb T}

\allowdisplaybreaks
\theoremstyle{remark}

\newcommand{\av}[1]{\langle #1\rangle}
\hypersetup{
  pdftitle={Analytic KAM Tori for Contact Hamiltonian Systems},
  pdfauthor={Zhicheng Tong, Yong Li, Xuefeng Zhao},
  pdfsubject={Legendrian tori, contact KAM persistence, and local dynamics}
}

\begin{document}
\title{{\sl KAM theory for contact Hamiltonian systems}}
\author{Zhicheng Tong\thanks{College of Mathematics, Jilin University,
		Changchun 130012, P. R. China. Email: \url{tongzc25@jlu.edu.cn}}
	\and
	Yong Li\thanks{College of Mathematics, Jilin University,
		Changchun 130012, P. R. China; School of Mathematics and Statistics,
		and Center for Mathematics and Interdisciplinary Sciences,
		Northeast Normal University, Changchun 130024, P. R. China.
		Corresponding author. Email: \url{liyong@jlu.edu.cn}}
	\and
	Xuefeng Zhao\thanks{College of Mathematics, Jilin University,
		Changchun 130012, P. R. China. Email: \url{zhaoxuef@jlu.edu.cn}}
}
\date{}
\maketitle

\begin{abstract}
In both finite and infinite dimensions, we establish the first frequency-preserving KAM theorems for full-dimensional Legendrian tori in contact Hamiltonian families with general analytic coupling between the momenta and the contact variable. After adjustment of external parameters, the flow on each torus is analytically conjugate in the original time to the prescribed translation, with estimates uniform as the damping parameter crosses zero. This coupling can generate nonzero mean damping at first order even when the damping parameter vanishes. In finite dimensions, nonzero mean damping yields a time-preserving analytic contact linearization near the torus. In both settings the invariant graphs give global analytic solutions of stationary Hamilton--Jacobi equations whose characteristic flows are conjugate to the prescribed translation. A novel consequence in finite dimensions is that, under global Tonelli assumptions and strict increase in the contact variable, this analytic solution is both the unique forward and the unique backward weak KAM solution, and its associated Legendrian graph is an invariant torus carrying quasi-periodic dynamics.
 \vspace{2mm}
\end{abstract}

\tableofcontents

\section{Introduction}\label{sec:introduction}

In this paper, we study the persistence of full-dimensional Legendrian tori with
prescribed frequency in finite- and infinite-dimensional analytic
contact Hamiltonian systems. We first describe the
finite-dimensional model on
$T^*\T^n\times\R$ with contact form $\alpha=dz-p\cdot dq$,
where $\T^n=\R^n/(2\pi\Z)^n$. Our convention is
$\alpha(X_H)=-H$ and $\iota_{X_H}d\alpha=dH-H_z\alpha$.
The equations of motion read
\[
    \dot q=H_p,\quad
    \dot p=-H_q-pH_z,\quad
    \dot z=p\cdot H_p-H.
\]
In particular, $\mathcal L_{X_H}\alpha=-H_z\alpha$, and
\begin{equation}\label{eq:introenergyvolume}
    X_H(H)=-HH_z,\quad
    \mathcal L_{X_H}\bigl(\alpha\wedge(d\alpha)^n\bigr)
    =-(n+1)H_z\,\alpha\wedge(d\alpha)^n.
\end{equation}
Thus the contact framework includes volume contraction and expansion,
as well as strict contact dynamics when $H_z=0$; see
\cite{GeigesContact,Bravetti}.

Our unperturbed Hamiltonian is
$H^{(0)}_\lambda=\omega\cdot p+\lambda z+R_0(q,p,z)$, where
$R_0(q,0,0)=0$ and $D_{(p,z)}R_0(q,0,0)=0$.
We abbreviate these conditions as $R_0=O_2$.
The zero section $\mathcal T_0=\{p=z=0\}$ is a Legendrian
invariant torus with flow $\dot q=\omega$. The background $R_0$
may depend on the angles, so this assumption does not require
integrability on a whole neighborhood of the torus. We study its
perturbation within the family
\begin{equation}\label{eq:family}
    H_{\xi,e,\lambda,\varepsilon}(q,p,z)
    =\xi\cdot p+\lambda z+R_0(q,p,z)
       +\varepsilon P(q,p,z)-e.
\end{equation}
Here $R_0,P$ are fixed analytic functions, $\lambda$ ranges over
$[-\Lambda,\Lambda]$, and $\varepsilon$ is small and real.
The parameters $\xi\in\R^n$ and $e\in\R$ will be selected
near $\omega$ and zero. Note that unlike a constant shift in a symplectic
Hamiltonian, $e$ \textit{changes} the vector field:
$X_{H-e}=X_H+e\partial_z$.

The invariant torus is a Legendrian graph
$j^1\mathfrak u=(q,D\mathfrak u(q),\mathfrak u(q))$.
It is \emph{full-dimensional} in the KAM sense because the angular
projection $(q,p,z)\mapsto q$ restricts to a diffeomorphism from
the torus onto the entire angular base $\T^n$.
Its invariance is equivalent to
$H(q,D\mathfrak u(q),\mathfrak u(q))=0$.
Preserving the specified frequency also requires an analytic torus
diffeomorphism $Q$ for which
\[
    H_p\bigl(Q(\theta),D\mathfrak u(Q(\theta)),
                    \mathfrak u(Q(\theta))\bigr)=DQ(\theta)\omega.
\]
Thus an invariant graph and a graph carrying a prescribed linear
flow are different requirements.

For a Diophantine frequency $\omega$,
the finite-dimensional Theorem~\ref{thm:main}
selects $\xi,e$ and constructs a near-identity analytic contact
embedding $\Phi$, with positive real multiplier $\rho=\rho(q)$,
such that
\begin{equation}\label{eq:target}
    \Phi^*\alpha=\rho\alpha,\quad
    \rho^{-1}H_{\xi,e,\lambda,\varepsilon}\circ\Phi
       =\omega\cdot p+\beta z+O_2.
\end{equation}
Its restriction $K(\theta)=\Phi(\theta,0,0)$ satisfies
$X_H\circ K=DK\,\omega$ and $K^*\alpha=0$, in the original
time. The parameter corrections and the displacement of $\Phi$
are $O(|\varepsilon|)$ on a common smaller complex domain.
Both this domain and the smallness threshold are uniform in real
$\lambda\in[-\Lambda,\Lambda]$, including $\lambda=0$.
The selected parameters and embeddings are continuous in $\lambda$.
No invertibility assumption is imposed on the normal Hessian of
$R_0$; the constants may depend on its fixed analytic bound.

For countably many angles, Theorem~\ref{thm:infinite} gives
persistence of a full-dimensional Legendrian torus in
the spatial framework of \cite{Poschel1990}, under the uniform
spatial composition hypothesis in Remark~\ref{rem:infinite-composition}.
The graph is defined over the entire torus $\T^{\mathcal I}$ and
carries the prescribed frequency vector $(\omega_i)_{i\in\mathcal I}$.
The actions belong to a weighted $\ell^1$ space, the frequency
correction is bounded, and $z$ remains scalar. The spatial weights
control both the perturbation and the small divisors. Because all
spatial blocks share $z$, their contact brackets also join disjoint
supports; we incorporate this interaction in the spatial structure.
The resulting invariant graph and parameter corrections satisfy
uniform $O(|\varepsilon|)$ bounds through zero real damping.

The coefficient $\beta$ is an output of the construction.
Differentiating the contact pullback along the normalized torus gives
\begin{equation}\label{eq:intromultiplier}
    \beta=H_z\circ K+D_\omega\log\rho,
    \quad D_\omega=\omega\cdot D_q.
\end{equation}
Hence $\beta=\langle H_z\circ K\rangle$, where the brackets denote
angular averaging, defined by the zero Fourier mode in the
infinite-dimensional setting. The multiplier removes the
oscillatory damping but does not change its average on this torus.
In particular, input $\lambda=0$ need not give actual mean damping
$\beta=0$. Both persistence theorems select the parameters in
their respective adjustable families for one prescribed frequency.

\subsection{Relation to earlier work}

For geometric background and related developments in contact
Hamiltonian dynamics, see, for instance,
\cite{Bravetti2017,deLeonSardon2017,deLeonLainz2019,
deLeonLainz2019Precontact,VermeerenBravettiSeri2019,deLeonLainz2020,GasetEtAl2020,
SimoesEtAl2020,deLeonEtAl2023,
deLeonLainzMunoz2023}
and the references therein.

Classical KAM theory originates in the work of Kolmogorov, Arnold,
and Moser \cite{Kolmogorov1954,Moser1962,Arnold1963}. For
\[
    h_\varepsilon(q,I)=h_0(I)+\varepsilon f(q,I),
    \quad (q,I)\in\T^n\times\R^n,
\]
the Kolmogorov condition $\det D^2h_0(I_0)\ne0$ permits one to
compensate a frequency change by varying the action. Parameter
formulations are discussed in \cite{Poschel,Bounemoura}.
A global Legendrian graph does not have an independent constant
momentum displacement: periodicity and $D\mathfrak u\equiv a$
imply $a=0$. In \eqref{eq:family}, the external parameter $\xi$
supplies the freedom needed for the angular average equations.

The proof follows the modifying-terms strategy in
\cite{Moser} and its parameterized development in
\cite{Wagener}. The latter is a structural persistence theorem and
allows certain vanishing normal Floquet exponents. The present work
gives a direct contact implementation for the specific family
\eqref{eq:family}. In finite dimensions, we compute the parameter response
after every accumulated contact transformation and solve the $n+1$
scalar and angular average equations using $(e,\xi)$, and retain the last
average as the coefficient $\beta$. No general normal matrix is
prescribed: contact geometry makes its normal linear part
$-\beta I_{n+1}$. The required parameter matrix starts at the
identity and remains invertible by an estimate with a fixed analytic
margin. This estimate gives a common perturbation threshold for all real
damping parameters in the prescribed interval, including zero.

The theorem allows any fixed analytic background $R_0=O_2$. A constant
contact dilation makes $R_0$ small on a rescaled normal domain, as
described in Remark~\ref{rem:backgroundscaling}. Here we keep the
original coordinates and retain the full Hessian in the linear
correction equations. The working radius and perturbation threshold
are chosen in terms of the analytic bound on $R_0$.

An a posteriori KAM theory for conformally symplectic systems was
established in \cite{Calleja2013}. Local normal forms giving
rotation--contraction coordinates and phase variables were obtained
in \cite{CallejaLocal2013}, and the limit of small dissipation was
studied in \cite{Calleja2017}. Theorem~\ref{thm:main} treats
coupling between momentum and the contact variable in the full contact
space. The weighted pullback in \eqref{eq:target} preserves the
original time, and the local linearization preserves the contact
structure. Related persistence and
parameterization results for normally hyperbolic tori are treated in
\cite{Canadell,Broer}.

Contact integrability is developed in
\cite{KhesinTabachnikov,Boyer2011,ColomboContact}.
In particular, \cite[Theorem~5(iii)]{Jovanovic2025} identifies
invariant Legendrian tori with linearized dynamics under integrability
and compactness assumptions at a zero momentum level. Here we start
with one such local dynamical model and study its perturbation,
without an integrable foliation in a neighborhood. KAM methods also
occur in closing problems \cite{AsaokaIrie,Xue}, but the torus
considered here is not a Reeb torus: a Reeb field is transverse to
the contact distribution, whereas an invariant Legendrian torus for
$X_H$ lies in $H^{-1}(0)$. Section~\ref{sec:symplectizationdetail}
explains the corresponding noncompact lift in the symplectization. Uniform $q$-contact structures have also been studied in
the context of contact foliations \cite{Finamore2024}.
Within this geometric framework, Hamiltonian dynamics was developed
and Liouville--Arnold-type theorems on complete integrability were
established in \cite{ZhaoLi}.

The contact weak KAM and variational literature includes
\cite{SuWangYan2016,WangWangYan2017,CannarsaEtAl2019,
WWY2019,CannarsaEtAl2020,WangWangYan2022}.
For related weak KAM results in almost periodic Hamiltonian
systems, see \cite{NiuWangLi2025}.
The Legendrian graphs give analytic classical solutions of the
stationary Hamilton--Jacobi equation in finite and infinite dimensions, with all
characteristics conjugate in the original time to the prescribed
translation. For countably many angles the Fr\'echet derivative is
represented by a momentum in $P_w$, so the Hamilton--Jacobi equation
holds pointwise on the entire torus. In finite dimensions, global
Tonelli and strict monotonicity assumptions identify the analytic
solution with every forward and backward weak KAM solution and its
graph with the contact Aubry and Mather sets.

\subsection{Proof outline and organization}

The iteration removes only the normal first-jet defect. At step $j$
we write $G_j=\omega\cdot p+\beta_jz+f_j+R_j$, with
$f_j=a_j+b_j\cdot p+c_jz$ and $R_j=O_2$.
A correction is the time-one contact map $\Psi_j$ generated by
$F_j=u_j+v_j\cdot p+w_jz$, and the accumulated map is updated
by $\Phi_{j+1}=\Phi_j\circ\Psi_j$.
The original parameters are adjusted before this pullback, so every
$G_j$ remains a transformed member of \eqref{eq:family}.
The quadratic blocks $R_{j,pp}|_0$, $R_{j,pz}|_0$, and
$R_{j,zz}|_0$ all enter the linear equations derived in
Section~\ref{sec:linear}.

In the finite-dimensional proof, the scalar and angular
averages are eliminated before Fourier inversion. The nonzero divisors are $\beta_j+ik\cdot\omega$
and $ik\cdot\omega$. For real $\beta_j$,
$|\beta_j+ik\cdot\omega|\geqslant|k\cdot\omega|$, so their
estimates are uniform through zero damping. An exact integral
remainder yields $E_{j+1}\leqslant C\sigma_j^{-m_*}E_j^2$, where
$E_j$ measures the first-jet defect. The contact-map increments have
$C^2$ norm at most $C\sigma_j^{-\ell_*}E_j$.
Geometric strip losses and quadratic error reduction make both the
losses and the corrections summable. The normal radius stays fixed
because each accumulated affine contact map is estimated against
the original analytic domain.

The rest of this paper is organized as follows. Section~\ref{sec:main} states our main results, namely the finite- and
infinite-dimensional persistence theorems. Sections~\ref{sec:geometry}--\ref{sec:zero}
develop the finite-dimensional argument and local dynamics.
Section~\ref{sec:geometry} establishes the contact identities
and Fourier estimates. Sections~\ref{sec:linear} and
\ref{sec:iteration} give the correction equations, the KAM step,
and convergence, including dependence on real damping.
Section~\ref{sec:extensions} treats alternative parameter choices
and computes the leading corrections, including an example with a
mixed $pz$ background. For nonzero actual mean damping,
Section~\ref{sec:linearization} proves a time-preserving analytic
contact linearization and asymptotic phase in the attracting case.
Section~\ref{sec:zero} describes the variational shear at zero
mean damping and gives exact quadratic models.
Section~\ref{sec:infinite-proof} proves the
infinite-dimensional theorem, including the weighted space estimates
and the choice of losses for general approximation functions.
Section~\ref{sec:hamiltonjacobi} derives the finite- and
infinite-dimensional Hamilton--Jacobi equations in
Corollaries~\ref{cor:hamiltonjacobi} and~\ref{cor:infinite-hj}
and identifies the weak KAM solutions and Aubry--Mather sets in
finite dimensions.
The appendix records
supplementary coefficient, composition, and tail estimates.

\section{Analytic frameworks and persistence theorems}\label{sec:main}

\subsection{Finite-dimensional notation and fixed data}

Fix $n\geqslant1$ and $s_0,r_0,\Lambda>0$. For $s,r>0$, set
\[
        \T_s^n =\bigl\{q\in\C^n/(2\pi\Z)^n: |\operatorname{Im}q|_\infty<s\bigr\},\quad D_{s,r} =\T_s^n\times \bigl\{(p,z)\in\C^n\times\C: |p|_\infty<r,\ |z|<r\bigr\}.
\]
We write $\|f\|_s=\sup_{\T_s^n}|f|$ for an angular function and
$\|G\|_{s,r}=\sup_{D_{s,r}}|G|$ for a function of all variables.
For a nonnegative integer $j$, the angular derivative norm is
$\|f\|_{C^j(s)}=\sum_{\nu\in\mathbb N_0^n,\,|\nu|\leqslant j}\|\partial_q^\nu f\|_s$.
Full-space $C^j$ norms are defined in the same way, with derivatives
in all variables $(q,p,z)$. We use the maximum norm for vectors and
the associated operator norms for matrices and multilinear maps.

For a periodic function, define
\[
    \av{f}=(2\pi)^{-n}\int_{\T^n}f(q)\,dq,
    \quad \Pi f=f-\av{f},
    \quad D_\omega=\omega\cdot D_q.
\]
The prescribed frequency $\omega\in\R^n$ satisfies the Diophantine
condition
\begin{equation}\label{eq:dio}
    |k\cdot\omega|\geqslant\gamma |k|_1^{-\tau},
    \quad k\in\Z^n\setminus\{0\},
    \quad \gamma>0,\quad \tau>n-1,
\end{equation}
where $|k|_1=\sum_{i=1}^n|k_i|$. Fourier series are taken with respect
to the $2\pi$-periodic angles:
\[
    f(q)=\sum_{k\in\Z^n}f_ke^{ik\cdot q},
    \quad
    f_k=(2\pi)^{-n}\int_{\T^n}f(q)e^{-ik\cdot q}\,dq.
\]
Thus $f_0=\av{f}$, and the Fourier multiplier of $D_\omega$ is
$ik\cdot\omega$.

For a multi-index $\nu\in\mathbb N_0^n$, write
$|\nu|=\sum_{i=1}^n\nu_i$ and $\nu!=\prod_{i=1}^n\nu_i!$.
Cauchy's estimate gives, for $0<\sigma<s$,
\[
    \|\partial_q^\nu f\|_{s-\sigma}
    \leqslant\nu!\,\sigma^{-|\nu|}\|f\|_s.
\]
Similarly, for $\mu\in\mathbb N_0^{n+1}$ and $0<r'<r$,
\[
    \|\partial_{(p,z)}^\mu G\|_{s,r'}
    \leqslant\mu!\,(r-r')^{-|\mu|}\|G\|_{s,r}.
\]
These estimates apply componentwise to vector- and matrix-valued
functions. The scalar supremum norm satisfies
$\|fg\|_s\leqslant\|f\|_s\|g\|_s$; product estimates for derivatives
involve only fixed constants depending on the dimension and the
derivative order.

Gradients are column vectors, and $DQ$ has entries
$\partial Q_i/\partial q_j$; we write $DQ^{-T}=(DQ)^{-T}$.
For a torus map close to the identity, choose a lift
$Q(q)=q+h(q)$ with periodic $h$ and set
$\|Q-\mathrm{id}\|_{C^j(s)}=\|h\|_{C^j(s)}$.
The same convention is used for the angular component of
$\Phi-\mathrm{id}$. The pullback of a vector field is
$\Phi^*X=(D\Phi)^{-1}(X\circ\Phi)$, with both $D\Phi$ and the
pullback evaluated at the input point.

In Sections~\ref{sec:geometry}--\ref{sec:zero}
and the appendix, $C$ denotes a positive constant that may
change from one estimate to the next. Unless specified otherwise,
it depends only on
$n,\omega,\gamma,\tau,s_0,r_0,\Lambda$ and the analytic bounds
$\|R_0\|_{s_0,r_0}$ and $\|P\|_{s_0,r_0}$.
It is independent of $\lambda\in[-\Lambda,\Lambda]$, of the
sufficiently small real parameter $\varepsilon$, and of the
iteration index. Dependence on a stepwise strip loss is displayed
explicitly. Constants given a subscript, such as $C_*$ or $B_0$,
are fixed when introduced.

\subsection{Normal order and parameters}\label{sec:notationguide}

The variables $y=(p,z)$ are normal to the zero torus
$\{p=z=0\}$. For a holomorphic function $R$, the notation $R=O_2$
means that $R(q,0)=0$ and $D_yR(q,0)=0$ for every $q$.
Equivalently, its Taylor expansion in $y$ has no constant or linear
term. Taylor's formula gives
\begin{equation}\label{eq:O2integral}
    R(q,y)
    =\int_0^1(1-t)D_y^2R(q,ty)[y,y]\,dt.
\end{equation}
In particular,
$|R(q,y)|\leqslant\tfrac12\sup_{0\leqslant t\leqslant1}
\|D_y^2R(q,ty)\|\,|y|_\infty^2$.
The condition $R=O_2$ does not require its normal Hessian to be small
or invertible.

For a holomorphic function $G$, define its first jet in the normal
variables by
\[
    \mathcal J_1G(q,p,z)
    =G(q,0,0)+G_p(q,0,0)\cdot p+G_z(q,0,0)z.
\]
Thus $G-\mathcal J_1G=O_2$. This projection keeps the full angular
dependence. A subscript $|_0$ denotes restriction to $p=z=0$.

We consider the family \eqref{eq:family}, with $\omega$, $R_0$,
and $P$ fixed. The quantities $\lambda\in[-\Lambda,\Lambda]$ and
$\varepsilon\in\R$ are inputs. The vector $\xi\in\R^n$ and the
scalar $e\in\R$ are selected to obtain an invariant torus with
frequency $\omega$. The coefficient $\beta$ is determined by the
normal form of the selected system; it is not an additional
parameter in \eqref{eq:family}.
We use $q$ for a coordinate in a contact chart and $\theta$ for the
angle parameter on the invariant torus. The final graph function is
denoted by $\mathfrak u$, to distinguish it from the scalar
coefficient $u$ of a generating Hamiltonian.

\subsection{The finite-dimensional persistence theorem}

\begin{theorem}
\label{thm:main}
Let $\omega$ satisfy \eqref{eq:dio}, and let $R_0$ and $P$ be bounded
and holomorphic on $D_{s_0,r_0}$, with real values at real points.
Assume that $R_0=O_2$; no normal-order condition is imposed on $P$.

There exist $0<s_*<s_0$, $0<r_*<r_0$, $\varepsilon_*>0$, and $C>0$
such that, for every real $\lambda\in[-\Lambda,\Lambda]$ and
$|\varepsilon|<\varepsilon_*$, one can choose
$\xi\in\R^n$, $e,\beta\in\R$, and a holomorphic embedding
$\Phi:D_{s_*,r_*}\to D_{s_0,r_0}$, real on real points, of the form
\begin{equation}\label{eq:affine}
    \Phi(q,p,z)
    =\Bigl(
        Q(q),\,
        DQ(q)^{-T}\bigl[\rho(q)p+DU(q)+zD\rho(q)\bigr],\,
        U(q)+\rho(q)z
      \Bigr).
\end{equation}
The functions $Q-\mathrm{id}$, $\rho$, and $U$ are holomorphic on
$\T_{s_*}^n$ and real on the real torus. The matrix $DQ$ is invertible
and $\rho$ is nonzero on this strip. On $\T^n$, the map $Q$ is a
real-analytic diffeomorphism homotopic to the identity and $\rho>0$.

Writing $H=H_{\xi,e,\lambda,\varepsilon}$ for the selected
Hamiltonian, one has $\Phi^*\alpha=\rho\alpha$ and
\[
    \rho^{-1}H\circ\Phi
    =\omega\cdot p+\beta z+R_*(q,p,z),
    \quad R_*=O_2.
\]
Moreover,
\begin{equation}\label{eq:mainbound}
        |\xi-\omega|+|e|+|\beta-\lambda| +\|\Phi-\mathrm{id}\|_{C^2(D_{s_*,r_*})} +\|\rho-1\|_{C^2(s_*)} \leqslant C|\varepsilon|.
\end{equation}

The map $K(\theta)=\Phi(\theta,0,0)$ parameterizes a full-dimensional
invariant Legendrian torus and satisfies
\[
        X_H\circ K=DK\,\omega, \quad K^*\alpha=0,\quad H\circ K=0, \quad \beta=\av{H_z\circ K}.
\]
Its image is the graph
$j^1\mathfrak u=\{(q,D\mathfrak u(q),\mathfrak u(q)):q\in\T^n\}$,
where $\mathfrak u=U\circ Q^{-1}$ is real analytic.
The conjugacy $\theta\mapsto K(\theta)$ preserves the original time.

The domain and constants are uniform in
$\lambda\in[-\Lambda,\Lambda]$. For each fixed real
$\varepsilon$ with $|\varepsilon|<\varepsilon_*$, the parameters
and embeddings can be chosen continuously in $\lambda$, with the
embeddings continuous in the $C^2$ norm on a common smaller complex
domain.
\end{theorem}

For $H=h(q,p)+\lambda z$ with constant $\lambda$, the projected
equations $\dot q=h_p$, $\dot p=-h_q-\lambda p$ define a flow
$\varphi_t$ on $(q,p)$ satisfying $\varphi_t^*\Omega=e^{-\lambda t}\Omega$,
where $\Omega=\sum_{i=1}^n dq_i\wedge dp_i$. This is the conformally
symplectic setting of \cite{Calleja2013,Calleja2017}.
Theorem~\ref{thm:main} also allows mixed terms in either $R_0$ or $P$.
For example, $H=\xi\cdot p+\lambda z+\varepsilon p_1z-e$ gives
$\dot q_1=\xi_1+\varepsilon z$. For $\varepsilon\ne0$, two points
above the same $(q,p)$ can thus have different angular velocities,
so their projected motion admits no autonomous equation on $(q,p)$.
Example~\ref{ex:mixedcoupling} shows how a mixed background term can
also produce nonzero mean damping at first order.

\begin{remark}\label{rem:parameters}
The parameter $\xi$ is used to solve the $n$ angular average equations,
and $e$ is used to solve the scalar graph average equation.
The latter changes the vector field, since
$X_{H-e}=X_H+e\partial_z$.
The coefficient $\beta$ is the average of $H_z$ in the parameter
angle $\theta$ for which the torus flow is linear, and it need not
equal $\lambda$. No uniqueness is asserted for the selected
parameters or torus parameterizations.
\end{remark}

\begin{remark}\label{rem:backgroundscaling}
The theorem is local in the normal variables. For $h>0$, the dilation
$\mathcal D_h(q,P,Z)=(q,hP,hZ)$ has contact factor $h$ and gives
\[
    h^{-1}H^{(0)}_\lambda\circ\mathcal D_h
    =\omega\cdot P+\lambda Z+h^{-1}R_0(q,hP,hZ).
\]
Fix a smaller angular strip and a normal radius in $(P,Z)$.
Taylor's formula \eqref{eq:O2integral} bounds the last term by
$Ch$ when $h$ is small enough to map this domain inside the
original one. The perturbation becomes
$h^{-1}\varepsilon P(q,hP,hZ)$, so its threshold then depends
on the chosen $h$. Thus the result permits any fixed bounded
quadratic background, but is not uniform over backgrounds of
arbitrarily large norm. The proof keeps the original coordinates
and includes the Hessian terms in \eqref{eq:homological}; the
external parameter response, rather than a twist assumption,
provides the required average invertibility.
\end{remark}

\subsection{Spatial structure and full-dimensional tori}
\label{sec:infinite-setting}

We now extend the persistence theorem to countably many angles,
using the spatial expansions and weighted action norms of
\cite[Sections~1--2]{Poschel1990}. The contact variable remains
scalar.

Let $\mathcal I$ be a countable set and let $\mathcal S$ be
a family of finite subsets of $\mathcal I$, containing
$\varnothing$, covering $\mathcal I$, and closed under finite
unions and intersections. Let $[A]\geqslant 0$ be a weight on
$\mathcal S$ satisfying
\begin{equation}\label{eq:infinite-weights}
 [\varnothing]=0,\quad
 A\subseteq B\Rightarrow [A]\leqslant[B],\quad
 [A\cup B]+[A\cap B]\leqslant[A]+[B].
\end{equation}
We impose the last inequality also when $A\cap B=\varnothing$,
since the common variable $z$ couples Hamiltonian terms with
disjoint spatial supports.
In particular, \eqref{eq:infinite-weights} implies
\[
 [A\cup B]\leqslant[A]+[B],
 \quad A,B\in\mathcal S.
\]
We use this subadditivity to group products and contact brackets
by the union of their spatial supports.

For $i\in\mathcal I$ and
$k\in\mathbb Z^{(\mathcal I)}$, define
\[
 [i]=\inf_{A\in\mathcal S,\,i\in A}[A],
 \quad
 [k]=\inf_{A\in\mathcal S,\,\operatorname{supp}k\subseteq A}[A],
 \quad
 |k|_1=\sum_{i\in\mathcal I}|k_i|.
\]
Here $\mathbb Z^{(\mathcal I)}$ denotes the finitely supported
integer sequences. A concrete example is
$\mathcal I=\mathbb Z^d$, with $\mathcal S$ the family of
all finite subsets and
\[
 [A]=1+\sum_{i\in A}
       \bigl(\log(1+|i|)\bigr)^a,\quad A\ne\varnothing,
\]
where $a>2$; see \cite[pp.~360--361]{Poschel1990}.

Fix $w\geqslant0$. We use the weighted action space
\[
 P_w=
 \left\{p=(p_i)_{i\in\mathcal I}:
 |p|_w:=\sum_{i\in\mathcal I}|p_i|e^{w[i]}<+\infty
 \right\},
 \quad
 V=\ell^\infty(\mathcal I),
 \quad
 E=V\times P_w\times\mathbb C.
\]
The space $E$ is equipped with the maximum product norm.
We regard the real angular torus as
\[
 \mathbb T^{\mathcal I}
 =
 \ell^\infty(\mathcal I;\mathbb R)/
 2\pi\ell^\infty(\mathcal I;\mathbb Z).
\]
We use the Banach quotient topology, which is finer than the
coordinatewise product topology and is noncompact when
$\mathcal I$ is infinite. Its complex analytic neighborhoods
are described by bounded complex angular displacements:
\[
 \mathbb T_s^{\mathcal I}
 =
 \left\{q\in
 \ell^\infty(\mathcal I;\mathbb C)/
 2\pi\ell^\infty(\mathcal I;\mathbb Z):
 \sup_{i\in\mathcal I}|\operatorname{Im}q_i|<s
 \right\}.
\]
We set
\[
 \mathcal D_{s,r}
 =
 \mathbb T_s^{\mathcal I}
 \times\{p\in P_w:|p|_w<r\}
 \times\{z\in\mathbb C:|z|<r\}.
\]
Holomorphic maps on these domains are understood in the
Fr\'echet sense in the bounded angular displacement and normal
variables. Their derivatives are measured in the corresponding
Banach operator norms.

The pairing between $P_w$ and $V$ is
\[
 \langle p,v\rangle
 =\sum_{i\in\mathcal I}p_iv_i,
 \quad
 |\langle p,v\rangle|
 \leqslant |p|_w\|v\|_\infty.
\]
This pairing is continuous and separates points in both
spaces.

We first specify the weak contact geometry used below.

\begin{definition}
\label{def:infinite-weak-contact}
A smooth one-form $\alpha$ on a Banach manifold $M$ is called
weakly contact in the sense used in this paper if it is
nowhere zero and its tangent bundle admits the splitting
\[
 TM=\ker\alpha\oplus\ker d\alpha.
\]
This splitting implies that $d\alpha$ is weakly non-degenerate
on $\ker\alpha$: if $v\in\ker\alpha_x$ and
$d\alpha_x(v,h)=0$ for every $h\in\ker\alpha_x$, then $v=0$.
Indeed, $d\alpha_x(v,\cdot)$ also vanishes on $\ker d\alpha_x$,
so $v\in\ker\alpha_x\cap\ker d\alpha_x=\{0\}$.
\end{definition}

On the model space $E$, consider
\[
 \alpha=dz-\sum_{i\in\mathcal I}p_i\,dq_i,
 \quad
 d\alpha=\sum_{i\in\mathcal I}dq_i\wedge dp_i.
\]
Both expressions are well defined by the continuous pairing
between $P_w$ and $V$. For tangent vectors
$u=(h,a,c)$ and $v=(k,b,d)$ in $E$, one has
\[
 \alpha(u)=c-\langle p,h\rangle,
 \quad
 d\alpha(u,v)=\langle h,b\rangle-\langle k,a\rangle.
\]
Since the pairing separates points, these identities give
\[
 \ker d\alpha=\operatorname{span}\{\partial_z\}.
\]
Moreover, every tangent vector admits the decomposition
\[
 (h,a,c)
 =
 \left(h,a,\langle p,h\rangle\right)
 +
 \left(0,0,c-\langle p,h\rangle\right).
\]
Consequently,
\[
 TE=\ker\alpha\oplus\operatorname{span}\{\partial_z\},
\]
and $\alpha$ is weakly contact in the sense of
Definition~\ref{def:infinite-weak-contact}. Its Reeb vector
field is
\[
 R_\alpha=\partial_z,
 \quad
 \alpha(R_\alpha)=1,
 \quad
 \iota_{R_\alpha}d\alpha=0.
\]

For a Hamiltonian $H$ in the spatial analytic class
specified below, we use the convention
\[
 \alpha(X_H)=-H,\quad
 \iota_{X_H}d\alpha=dH-H_z\alpha.
\]
Whenever the indicated derivatives belong to the required
spaces, the coordinate expression is
\[
 X_H=
 \left(H_p,-H_q-pH_z,\langle p,H_p\rangle-H\right).
\]
Here $H_p\in V$, $H_q\in P_w$, and the last component is
scalar. We work with Hamiltonians whose derivatives have these
representations; the spatial estimates below verify this
property for the Hamiltonians used in the iteration.

For $A\in\mathcal S$ we identify $\Z^A$ with the integer sequences
supported in $A$ by extension by zero. For a fixed spatial expansion,
write
\begin{equation}\label{eq:infinite-norm}
 \begin{aligned}
 &J=\sum_{A\in\mathcal S}J_A(q_A,p_A,z),\quad J_A=\sum_{k\in\Z^A} J_{A,k}(p_A,z)e^{ik\cdot q},\\
 &|J|_{m,s,r} = \sum_{A\in\mathcal S}e^{m[A]} \sum_{k\in\Z^A} e^{s|k|_1} \sup_{|p|_w<r,\,|z|<r} |J_{A,k}(p_A,z)|.
 \end{aligned}
\end{equation}
We fix expansions for the data $R_0,P$. For a function or coefficient
constructed in the iteration, its spatial norm is the infimum of the
expansion norms over all representations with the indicated support and
coefficient spaces. For affine Hamiltonians we take affine
representatives. We use this convention also for the contact and
operator norms below. Summation is continuous in the supremum norm,
so its kernel is closed and the quotient of the completed expansion
space is Banach. Products are grouped by $A\cup B$.
Terms depending only on $z$ are assigned to
$A=\varnothing$. In this block,
\[
 |c|_{m,s,r}=|c|,
 \quad
 |cz|_{m,s,r}=r|c|.
\]
In particular, the scalar parameter $-e$ is included in
this class and satisfies $X_{-e}=e\partial_z$.
The bounded frequency term
$\xi\cdot p$, with $\xi\in V$, is kept separately.

The condition $m>w$ provides a spatial weight margin for
the Hamiltonian derivatives. On smaller angular strips and normal
balls, the estimates for both $H_q\in P_w$ and $H_p\in V$
cost at most a factor $e^{w[A]}$ on the block with support $A$.
Thus a scalar spatial norm at weight $m$ controls the Hamiltonian
and its vector field at weight $m-w$, with the usual angular and
normal Cauchy losses. We use this conversion once, then estimate
the contact brackets in the resulting norm. The estimates are
proved in Section~\ref{sec:infinite-brackets}.

Let $\Delta:[0,\infty)\to[1,\infty)$ be nondecreasing,
with $\Delta(0)=1$, and assume
\begin{equation}\label{eq:infinite-approximation}
 \frac{\log\Delta(t)}{t}\searrow0
 \quad(t\to\infty),
 \quad
 \int_1^\infty\frac{\log\Delta(t)}{t^2}\,dt<+\infty.
\end{equation}
We prescribe a frequency
$\omega\in V\cap\mathbb R^{\mathcal I}$ satisfying
the nonresonance condition of \cite{Poschel1990}
\begin{equation}\label{eq:infinite-nonresonance}
 |k\cdot\omega|
 \geqslant
 \frac{\gamma}{\Delta([k])\Delta(|k|_1)},
 \quad
 0\ne k\in\mathbb Z^{(\mathcal I)},
 \quad \gamma>0.
\end{equation}
The dot product is well defined because $k$ has finite
support.

For an absolutely convergent spatial Fourier expansion,
we define its average by its zero Fourier coefficients:
\[
 \av J
 =
 \sum_{A\in\mathcal S}J_{A,0}.
\]
The sum converges absolutely whenever the spatial norm is finite.
For functions of finitely many angles, it agrees with the average
defined above. We use this zero-mode definition in the
infinite-dimensional setting.

\begin{remark}\label{rem:infinite-composition}
For the infinite-dimensional iteration, we use the following
spatial composition property on a fixed range of working widths.
Put
\[
\widehat m=\frac34(m_0-w),\quad
\widehat s=\frac34s_0,\quad
\widehat r=\frac12r_0,\quad
m_-=\frac{\widehat m}{2},\quad s_-=\frac{\widehat s}{4},
\]
and choose $r>0$ such that $16r<\widehat r$.
For an affine contact map $\Phi$ of the form \eqref{eq:affine}, write
\[
 A_\Phi=U+(Q-\mathrm{id})\cdot p+(\rho-1)z,
 \quad Y_{m,s}(\Phi)=|A_\Phi|^{c,2}_{m,s,8r},
\]
where the norm $|\cdot|^{c,2}$ is defined in
Section~\ref{sec:infinite-inverse-average}.
The coefficients have spatial expansions, with $Q-\mathrm{id}$
supported in the coordinates of its block. Products are grouped
by unions of supports. Their operator coefficients are local
on $V$ and $P_w$, together with scalar multiples of the identity;
the transposes are taken with respect to the $P_w,V$ pairing.
The coefficient extraction and Neumann-series estimates in
Section~\ref{sec:infinite-inverse-average} apply at each working
width. In the full Jacobian, a term consisting of the common
momentum $p$ times a scalar functional is assigned the support
of that functional and measured in the corresponding Banach
operator norm on the fixed normal ball.

We assume that there exist $Y_*,C_{\mathrm{sc}}>0$ such that,
for every
\[
 m_-\leqslant m\leqslant\widehat m,\quad
 s_-\leqslant s\leqslant\widehat s/2,
\]
and every such affine contact map satisfying
$Y_{m,s}(\Phi)\leqslant Y_*$, the map $\Phi$ sends
$\mathcal D_{s,8r}$ into
$\mathcal D_{\widehat s,\widehat r}$, and the expansions
obtained by substitution satisfy
\begin{equation}\label{eq:infinite-composition-hypothesis}
 |\rho^{-1}J\circ\Phi|^c_{m,s,8r}
 \leqslant C_{\mathrm{sc}}|J|_{m_0,s_0,r_0},
 \quad J=R_0,P.
\end{equation}
The constants are uniform over these widths and maps.
The contact norm in \eqref{eq:infinite-composition-hypothesis}
controls the transformed Hamiltonian and its vector field;
the map is controlled in the stronger $c,2$ norm.

For example, this property holds for
\[
 R_0=a z^2,\quad
 P=\sum_{i\in\mathbb Z}a_i\cos q_i,\quad
 a_i=2^{-|i|-2}e^{-(m_0+w+1)[i]},\quad a\in\mathbb R,
\]
with the preceding spatial structure on $\mathbb Z$ and
$\mathcal S$ consisting of all finite subsets. Indeed,
$\rho^{-1}R_0\circ\Phi
=a\rho z^2+2aUz+a\rho^{-1}U^2$.
Writing $h=Q-\mathrm{id}$, the exponential expansion of
$e^{\pm i(q_i+h_i)}$ and the union-weight inequality control
$\rho^{-1}P\circ\Phi$. Its angular derivative is the sum of
$(DQ)^T(DP)(Q)/\rho$ and $-P(Q)D\rho/\rho^2$.
The transpose bounds control these expressions in $P_w$;
the original factor $e^{w[i]}$ is absorbed by
$m_0-(m+w)\geqslant(m_0-w)/4>0$.
The remaining scalar and contact-vector-field components
are bounded by the same product estimates. This proves
\eqref{eq:infinite-composition-hypothesis} uniformly, and
small $Y_*$ gives the stated domain inclusion.
The example depends on infinitely many angular coordinates.
\end{remark}

We can now state the persistence theorem.

\begin{theorem}
\label{thm:infinite}
Assume \eqref{eq:infinite-weights},
\eqref{eq:infinite-approximation}, and
\eqref{eq:infinite-nonresonance}. Let $m_0>w$,
$s_0,r_0,\Lambda>0$, and suppose that $R_0$ and $P$
belong to the real spatial analytic class defined by
\eqref{eq:infinite-norm} on $\mathcal D_{s_0,r_0}$, with
\[
|R_0|_{m_0,s_0,r_0}+|P|_{m_0,s_0,r_0}<+\infty,
\quad R_0=O_2.
\]
The notation $R_0=O_2$ has the same meaning as in
Section~\ref{sec:notationguide}, with Fr\'echet derivatives
in $(p,z)$.

Assume, in addition, the spatial composition property
stated in Remark~\ref{rem:infinite-composition}, with
constants depending only on the fixed data and the displayed
spatial analytic norms.

Consider the family
\[
H_{\xi,e,\lambda,\varepsilon}
=
\sum_{i\in\mathcal I}\xi_i p_i
+\lambda z+R_0(q,p,z)
+\varepsilon P(q,p,z)-e,
\]
where $\xi\in V\cap\mathbb R^{\mathcal I}$ and $e\in\mathbb R$.

There exist $0<s_*<s_0$, $0<r_*<r_0$, and
$\varepsilon_*,C>0$, depending only on the fixed data,
the spatial composition constants, and the displayed
spatial analytic norms, such that for every
$\lambda\in[-\Lambda,\Lambda]$ and
$|\varepsilon|<\varepsilon_*$ there are
\[
\xi\in V\cap\mathbb R^{\mathcal I},
\quad e,\beta\in\mathbb R,
\]
and a holomorphic embedding
\[
\Phi:\mathcal D_{s_*,r_*}
\to\mathcal D_{s_0,r_0}
\]
of the affine contact form \eqref{eq:affine}. Here
$Q-\mathrm{id}$ is a bounded holomorphic $V$-valued
displacement, $DU$ and $D\rho$ are $P_w$-valued
holomorphic maps, and $DQ^{-T}$ is a bounded operator
on $P_w$ compatible with the dual pairing between $P_w$
and $V$. The spatially weighted operator coefficients
of $\Phi$ and $(D\Phi)^{-1}$ belong to the algebra specified
in Remark~\ref{rem:infinite-composition}. The map is real
on real points, $Q$ is an analytic diffeomorphism of the
real torus, and $\rho>0$ there.

Writing $H=H_{\xi,e,\lambda,\varepsilon}$ for the selected
Hamiltonian, one has
\[
\Phi^*\alpha=\rho\alpha,\quad
\rho^{-1}H\circ\Phi=\omega\cdot p+\beta z+R_*,
\quad R_*=O_2.
\]
The normal-form identity is understood on
$\mathcal D_{s_*,r_*}$, with the transformed remainder
belonging to the corresponding smaller spatial analytic
class.

Moreover,
\begin{equation}\label{eq:infinite-bound}
\|\xi-\omega\|_\infty+|e|+|\beta-\lambda| +\|\Phi-\mathrm{id}\|_{C^2(\mathcal D_{s_*,r_*};E)} +\|\rho-1\|_{C^2(\mathbb T_{s_*}^{\mathcal I})} \leqslant C|\varepsilon|.
\end{equation}
The $C^2$ norms are taken with respect to the induced
Banach operator norms, and the angular displacement is
measured using its bounded periodic lift.

The embedding $K(\theta)=\Phi(\theta,0,0)$ satisfies
\[
X_H\circ K=DK\,\omega,\quad
K^*\alpha=0,\quad H\circ K=0,
\quad
\beta=\av{H_z\circ K}.
\]
Its image is the full-dimensional invariant Legendrian torus
\[
j^1\mathfrak u
=
\{(q,D\mathfrak u(q),\mathfrak u(q)):
q\in\mathbb T^{\mathcal I}\},
\quad
\mathfrak u=U\circ Q^{-1},
\]
where $D\mathfrak u$ takes values in $P_w$.

The domain and constants are uniform in
$\lambda\in[-\Lambda,\Lambda]$. For each fixed real
$\varepsilon$ with $|\varepsilon|<\varepsilon_*$, the
parameters and embeddings can be chosen continuously
in $\lambda$, with the embeddings continuous in the
$C^2$ norm on a common smaller complex domain.
\end{theorem}

\begin{remark}\label{rem:infinite-polynomial}
A different KAM theorem for infinite-dimensional contact
Hamiltonian systems can also be established using the
setting introduced in \cite{TongLiBreathers}.
\end{remark}

\section{Contact geometry and uniform cohomological estimates}
\label{sec:geometry}

\subsection{Contact Hamiltonian fields and invariant graphs}
\label{sec:contactfoundations}

Let $\alpha$ be a contact form on a $(2n+1)$-dimensional manifold
$M$, so that $\alpha\wedge(d\alpha)^n$ is nowhere zero.
Its Reeb field $R_\alpha$ is determined by
$\alpha(R_\alpha)=1$ and $\iota_{R_\alpha}d\alpha=0$.
We use the convention
\begin{equation}\label{eq:intrinsiccontact}
    \alpha(X_H)=-H,\quad
    \iota_{X_H}d\alpha=dH-(R_\alpha H)\alpha.
\end{equation}
These equations determine $X_H$ uniquely: its Reeb component is
$-H R_\alpha$, and its component in $\ker\alpha$ is determined by
the nondegeneracy of $d\alpha$ on that distribution.

For $\alpha=dz-p\cdot dq$, one has $R_\alpha=\partial_z$ and
$d\alpha=\sum_{i=1}^n dq_i\wedge dp_i$. Writing
$X_H=A\partial_q+B\partial_p+C\partial_z$, the second equation
in \eqref{eq:intrinsiccontact} gives
$A=H_p$ and $B=-H_q-pH_z$; the first then gives
$C=p\cdot H_p-H$. Thus
\begin{equation}\label{eq:contact}
        X_H=H_p\partial_q-(H_q+pH_z)\partial_p +(p\cdot H_p-H)\partial_z,\quad \mathcal L_{X_H}\alpha=-H_z\alpha.
\end{equation}
The second identity follows from Cartan's formula and
\eqref{eq:intrinsiccontact}.

\begin{lemma}
\label{lem:energyReeb}
The vector field \eqref{eq:intrinsiccontact} satisfies
\eqref{eq:introenergyvolume}, with $H_z$ replaced by $R_\alpha H$
on a general contact manifold. Its flow preserves $H^{-1}(0)$
wherever it is defined. On a region where $H\ne0$, the field $X_H$
is the Reeb field of $-\alpha/H$. A Reeb field cannot be tangent
to a Legendrian submanifold for the same contact distribution.
\end{lemma}

\begin{proof}
Put $b=R_\alpha H$. Evaluating the second equation in
\eqref{eq:intrinsiccontact} on $X_H$ gives
$X_H(H)=-Hb$. Hence, along any trajectory $x(t)$,
\[
    H(x(t))=H(x(0))
        \exp\left(-\int_0^t b(x(a))\,da\right).
\]
In particular, the zero level is invariant, without any
regular-value assumption.

Cartan's formula gives $\mathcal L_{X_H}\alpha=-b\alpha$ and
$\mathcal L_{X_H}d\alpha=-db\wedge\alpha-bd\alpha$.
Consequently,
\[
    \mathcal L_{X_H}\bigl(\alpha\wedge(d\alpha)^n\bigr) =-b\alpha\wedge(d\alpha)^n +n\alpha\wedge(-db\wedge\alpha-bd\alpha) \wedge(d\alpha)^{n-1} =-(n+1)b\,\alpha\wedge(d\alpha)^n.
\]

Where $H\ne0$, the form $\widetilde\alpha=-\alpha/H$ is contact
and satisfies $\widetilde\alpha(X_H)=1$. Moreover,
\[
    \mathcal L_{X_H}\widetilde\alpha
    =\frac{X_H(H)}{H^2}\alpha
       -\frac1H\mathcal L_{X_H}\alpha=0.
\]
Cartan's formula therefore gives
$\iota_{X_H}d\widetilde\alpha=0$, proving the Reeb assertion.
Finally, a field tangent to a Legendrian submanifold is annihilated
by its contact form, whereas a Reeb field has value one under that
form.
\end{proof}

The Reeb interpretation uses a contact form that is singular at
$H=0$. The invariant Legendrian tori considered here lie precisely
on this level, as the following proposition shows.

\begin{proposition}
\label{prop:graphHJ}
Let $\mathfrak u\in C^2(\T^n)$ and let $H$ be $C^2$ on a
neighborhood of its graph. The embedding
$K_{\mathfrak u}(q)=(q,D\mathfrak u(q),\mathfrak u(q))$
is Legendrian. The field $X_H$ is tangent to its image if and only if
\begin{equation}\label{eq:graphHJcompatibility}
    H(q,D\mathfrak u(q),\mathfrak u(q))=0,
    \quad q\in\T^n.
\end{equation}
When this holds, the induced vector field on the base is
$X(q)=H_p(q,D\mathfrak u(q),\mathfrak u(q))$.
A $C^1$ torus diffeomorphism $Q$ conjugates this field to
$\dot\theta=\omega$ in the same time precisely when
\begin{equation}\label{eq:baseconjugacydetail}
    X(Q(\theta))=DQ(\theta)\omega.
\end{equation}
\end{proposition}

\begin{proof}
The identity
$K_{\mathfrak u}^*\alpha=d\mathfrak u-D\mathfrak u\cdot dq=0$
proves the Legendrian property. If $X_H$ is tangent to the graph,
then $\alpha(X_H)=0$ there, so \eqref{eq:intrinsiccontact} gives
\eqref{eq:graphHJcompatibility}.

Conversely, differentiation of \eqref{eq:graphHJcompatibility}
gives $H_q+D^2\mathfrak u\,H_p+H_zD\mathfrak u=0$ along the
graph. Substitution into \eqref{eq:contact} yields
\[
    X_H\circ K_{\mathfrak u}
    =\bigl(H_p,D^2\mathfrak u\,H_p,
                    D\mathfrak u\cdot H_p\bigr)
    =DK_{\mathfrak u}\,X.
\]
This proves tangency and identifies the induced field.
For $K=K_{\mathfrak u}\circ Q$, the equation
$X_H\circ K=DK\,\omega$ is then equivalent to
\eqref{eq:baseconjugacydetail}.
\end{proof}

\begin{remark}\label{rem:actionobstruction}
A graph $(q,P(q),Z(q))$ is Legendrian if and only if $P=DZ$.
Its momentum one-form is therefore exact. In particular, if
$D\mathfrak u=a$ is constant, periodicity gives
$0=\mathfrak u(q+2\pi e_i)-\mathfrak u(q)=2\pi a_i$
for each $i$, and hence $a=0$.
Similarly, the map $(q,p,z)\mapsto(q,p+a,z+f(q))$ is strict
contact if and only if $df=a\cdot dq$, which has a periodic
solution only when $a=0$. Thus a constant momentum translation
is not an independent global Legendrian graph parameter.
\end{remark}

\subsection{The symplectization lift}
\label{sec:symplectizationdetail}

On $M\times\R_+$, with coordinate $\upsilon>0$, set
$\Omega_{\rm sym}=d(\upsilon\alpha)
=d\upsilon\wedge\alpha+\upsilon\,d\alpha$.
It is symplectic, since
\[
    \Omega_{\rm sym}^{n+1}
    =(n+1)\upsilon^n
       d\upsilon\wedge\alpha\wedge(d\alpha)^n\ne0.
\]
With the convention
$\iota_{X_{\widetilde H}}\Omega_{\rm sym}=d\widetilde H$,
the Hamiltonian $\widetilde H=\upsilon H$ has vector field
\[
    \widetilde X=X_H+\upsilon(R_\alpha H)\partial_\upsilon.
\]
Indeed,
\[
    \iota_{\widetilde X}\Omega_{\rm sym} =\upsilon(R_\alpha H)\alpha-\alpha(X_H)d\upsilon +\upsilon\,\iota_{X_H}d\alpha =H\,d\upsilon+\upsilon\,dH =d(\upsilon H).
\]
The projection of the lifted flow is the contact flow with the
same time parameter.

Over the zero torus of
$H=\omega\cdot p+\lambda z+O_2$, one has
$\upsilon(t)=\upsilon(0)e^{\lambda t}$.
For $\lambda\ne0$, its lift at fixed $\upsilon$ is not invariant.
When $\lambda=0$, a fixed-$\upsilon$ lift is invariant but has
dimension $n$, rather than the Lagrangian dimension $n+1$.
The cylinder $\mathcal T_0\times\R_+$ is Lagrangian, since
$\upsilon\alpha$ restricts to zero on it, but it is noncompact.
Thus this lift does not directly supply a compact Lagrangian torus
to which the usual full-dimensional Hamiltonian KAM theorem applies.

\subsection{Contact pullback and composition}

Consider a map of the form \eqref{eq:affine}, on a domain where
$DQ$ is invertible and $\rho$ is nonzero. Direct differentiation gives
\[
    d(U+\rho z)
      -\bigl[DQ^{-T}(\rho p+DU+zD\rho)\bigr]\cdot dQ
    =\rho(dz-p\cdot dq).
\]
Hence $\Phi^*\alpha=\rho\alpha$. Conversely, prescribing
$q'=Q(q)$, $z'=U(q)+\rho(q)z$, and this contact identity forces
$DQ^Tp'=\rho p+DU+zD\rho$, and therefore determines the momentum
component in \eqref{eq:affine}.

Let $G=\rho^{-1}H\circ\Phi$ and set
$Y=(D\Phi)^{-1}(X_H\circ\Phi)$.
Then $\alpha(Y)=-G$, while naturality of Lie derivatives gives
\[
    \mathcal L_Y(\rho\alpha)
    =\Phi^*(\mathcal L_{X_H}\alpha)
    =-(H_z\circ\Phi)\rho\alpha.
\]
Expanding the left-hand side shows that
$\mathcal L_Y\alpha=\mu\alpha$, where
$\mu=-H_z\circ\Phi-\rho^{-1}Y(\rho)$.
Cartan's formula now gives
$\iota_Yd\alpha=dG+\mu\alpha$.
Evaluation on $\partial_z$ implies $\mu=-G_z$, so uniqueness in
\eqref{eq:intrinsiccontact} yields
\[
    \Phi^*X_H=X_{\rho^{-1}H\circ\Phi}.
\]
Consequently, $\Phi\circ\varphi_t^G=\varphi_t^H\circ\Phi$
wherever both sides are defined. The factor $\rho^{-1}$ is part of
the Hamiltonian transformation law; no change of time is made.

If $\Phi_1$ has triple $(Q_1,\rho_1,U_1)$, the composite
$\Phi_+=\Phi\circ\Phi_1$ has triple
\begin{equation}\label{eq:composition}
        Q_+=Q\circ Q_1,\quad \rho_+=(\rho\circ Q_1)\rho_1,\quad U_+=(\rho\circ Q_1)U_1+U\circ Q_1.
\end{equation}
The base and vertical components give these formulas, and
$\Phi_+^*\alpha=(\rho\circ Q_1)\rho_1\alpha$ determines the
momentum component. The affine contact class is therefore closed
under the compositions used in the iteration.

\subsection{A uniform Fourier inverse}

\begin{lemma}\label{lem:fourier}
Assume \eqref{eq:dio} and set $\kappa=\tau+n$.
For a bounded holomorphic function $f$ on $\T_s^n$ and
$\beta\in\R$, define
\[
    T_\beta f(q)
    =\sum_{k\in\Z^n\setminus\{0\}}
       \frac{f_k}{\beta+ik\cdot\omega}e^{ik\cdot q}.
\]
The series converges locally uniformly on $\T_s^n$ and defines
the unique zero-average holomorphic solution of
$(D_\omega+\beta)u=\Pi f$.
For $0<\sigma<\min\{s,1\}$ and $\beta,\beta'\in\R$,
\begin{align}
    \|T_\beta f\|_{s-\sigma}
    &\leqslant C\gamma^{-1}\sigma^{-\kappa}\|f\|_s,
        \label{eq:Tbound}\\
    \|(T_\beta-T_{\beta'})f\|_{s-\sigma}
    &\leqslant C\gamma^{-2}\sigma^{-(2\tau+n)}
               |\beta-\beta'|\,\|f\|_s.
        \label{eq:Tcontinuous}
\end{align}
Here $C$ depends only on $n$ and $\tau$.
Moreover, for fixed $a>0$ and an integer $j\geqslant0$, there is a
constant $C_{a,j}$, depending only on $n,\tau,\gamma,a,j$,
such that, whenever $s\geqslant a$,
\begin{equation}\label{eq:fixedmargin}
    \|T_\beta f\|_{C^j(\T^n)}\leqslant C_{a,j}\|f\|_s.
\end{equation}
All these bounds are uniform in real $\beta$ and $\beta'$.
If $f$ is real on the real torus, then so is $T_\beta f$.
\end{lemma}

\begin{proof}
Shift the Fourier integral to
$\operatorname{Im}q_i=-s'\operatorname{sgn}k_i$, with $s'<s$,
and let $s'\uparrow s$. This gives
$|f_k|\leqslant\|f\|_s e^{-s|k|_1}$.
For real $\beta$, the Diophantine condition implies
\[
    |\beta+ik\cdot\omega|
    =\bigl(\beta^2+(k\cdot\omega)^2\bigr)^{1/2}
    \geqslant\gamma|k|_1^{-\tau},\quad k\ne0.
\]
It follows that
\[
    \|T_\beta f\|_{s-\sigma}
    \leqslant\gamma^{-1}\|f\|_s
       \sum_{k\in\Z^n\setminus\{0\}}|k|_1^\tau e^{-\sigma|k|_1}.
\]
Write $N_m=\#\{k\in\Z^n:|k|_1=m\}$.
Counting the nonnegative integer solutions of
$|k_1|+\cdots+|k_n|=m$ and then allowing all choices of signs gives
$N_m\leqslant2^n\binom{m+n-1}{n-1}\leqslant C_nm^{n-1}$.
Consequently,
\[
    \sum_{k\in\Z^n\setminus\{0\}}|k|_1^\tau e^{-\sigma|k|_1}
    =\sum_{m=1}^{\infty}N_m m^\tau e^{-\sigma m}
    \leqslant C_n\sum_{m=1}^{\infty}m^{\tau+n-1}e^{-\sigma m}.
\]
For $a\geqslant0$, $0<\sigma\leqslant1$, and $x\in[m,m+1]$,
$m^a e^{-\sigma m}\leqslant e^\sigma x^a e^{-\sigma x}$.
Integrating on these intervals and summing yields
\begin{equation}\label{eq:sumshell}
    \sum_{m=1}^\infty m^a e^{-\sigma m}
    \leqslant e\int_0^\infty x^a e^{-\sigma x}\,dx
    =e\,\Gamma(a+1)\sigma^{-(a+1)}.
\end{equation}
Taking $a=\tau+n-1$ proves \eqref{eq:Tbound}.
The same argument with extra powers of $|k|_1$ proves local
uniform convergence of the differentiated series.
Termwise application of $D_\omega+\beta$ therefore gives $\Pi f$.
The zero Fourier coefficient is zero, and the nonzero denominators
also prove uniqueness under this normalization.

For the parameter difference, use
\[
    \frac1{\beta+ik\cdot\omega}
      -\frac1{\beta'+ik\cdot\omega}
    =\frac{\beta'-\beta}
           {(\beta+ik\cdot\omega)(\beta'+ik\cdot\omega)}.
\]
Thus
\[
    \|(T_\beta-T_{\beta'})f\|_{s-\sigma}
    \leqslant\gamma^{-2}|\beta-\beta'|\,\|f\|_s
          \sum_{k\in\Z^n\setminus\{0\}}|k|_1^{2\tau}e^{-\sigma|k|_1}.
\]
Applying the shell count and \eqref{eq:sumshell} with
$a=2\tau+n-1$ proves \eqref{eq:Tcontinuous}.

For $s\geqslant a>0$, differentiation on the real torus gives
\[
    \|T_\beta f\|_{C^j(\T^n)}
    \leqslant C(n,j)\gamma^{-1}\|f\|_s
       \sum_{k\in\Z^n\setminus\{0\}}(1+|k|_1^j)|k|_1^\tau e^{-a|k|_1}.
\]
The sum is finite and independent of $s$ and $\beta$, proving
\eqref{eq:fixedmargin}.
Finally, if $f$ is real, then $f_{-k}=\overline{f_k}$ and the
denominators at $-k$ and $k$ are conjugate. The same symmetry
holds for the coefficients of $T_\beta f$.
\end{proof}

For two sources $f,f'$ and real coefficients $\beta,\beta'$, the
same estimates give
\begin{equation}\label{eq:Tjointcontinuous}
        \|T_\beta f-T_{\beta'}f'\|_{s-\sigma} \leqslant C\gamma^{-1}\sigma^{-\kappa}\|f-f'\|_s +C\gamma^{-2}\sigma^{-(2\tau+n)} |\beta-\beta'|\,\|f'\|_s.
\end{equation}
Indeed, split the difference as
$T_\beta(f-f')+(T_\beta-T_{\beta'})f'$.
Cauchy's estimate gives the corresponding derivative convergence
on a further smaller strip. This joint estimate will be used for
the finite-stage parameter dependence.

Estimate \eqref{eq:fixedmargin} uses a fixed lower bound for the
strip width, rather than the loss $\sigma$ at a particular step.
It will be used to estimate the finite-dimensional average equations.
The oscillatory corrections are instead controlled on smaller
complex strips by \eqref{eq:Tbound}.

\subsection{Normalization of the damping coefficient}
\label{sec:multiplierdetail}

Suppose that the zero torus already carries the prescribed angular
flow and that
\[
    G(q,p,z)=\omega\cdot p+c(q)z+R(q,p,z),\quad R=O_2.
\]
For a real-analytic periodic function $g$, consider
\begin{equation}\label{eq:dampinggmap}
    \Psi(q,I,W)
    =\bigl(q,e^{-g(q)}[I-WDg(q)],e^{-g(q)}W\bigr).
\end{equation}
This map fixes the zero section pointwise and satisfies
\[
    d(e^{-g}W)-e^{-g}(I-WDg)\cdot dq
    =e^{-g}(dW-I\cdot dq).
\]
Its contact multiplier is therefore $e^{-g}$, and the transformed
Hamiltonian is
\begin{equation}\label{eq:dampingpullback}
    e^gG\circ\Psi
    =\omega\cdot I+(c-D_\omega g)W+e^gR\circ\Psi.
\end{equation}
The last term is $O_2$, since $\Psi$ is linear in $(I,W)$ and
preserves its origin. The coefficient of $W$ is constant precisely
when $D_\omega g=c-\beta$. Averaging gives $\beta=\av{c}$,
and Lemma~\ref{lem:fourier} provides the zero-average solution
$g=T_0c$ on any smaller angular strip.
On a fixed narrower strip where $g$ and $Dg$ are bounded, reducing
the normal radius makes \eqref{eq:dampinggmap} a well-defined
analytic contact map into the original domain. Its inverse is
$(q,p,z)\mapsto(q,e^{g(q)}[p+zDg(q)],e^{g(q)}z)$.
No smallness of $c$ is needed for this normalization, although
the domain and bounds may depend on it.

Strict contact transformations do not provide this freedom within
the affine class \eqref{eq:affine}.
Indeed, suppose such a transformation preserves the zero section,
is homotopic to the identity on the base, and preserves the
prescribed flow on that section. Then $U=0$, $\rho=1$, and
$DQ\,\omega=\omega$.
Writing $Q(q)=q+h(q)$ gives $D_\omega h=0$; by nonresonance,
$h$ is constant. Such a transformation only translates the argument
of $c$ and cannot make a nonconstant $c$ constant.
Allowing $\rho=e^{-g}$ adds the term
$D_\omega\log\rho=-D_\omega g$ in \eqref{eq:dampingpullback}.
At the infinitesimal level, a generator
$F=u(q)+v(q)\cdot p+w(q)z$ satisfies
$\mathcal L_{X_F}\alpha=-w\alpha$; the coefficient $w$ supplies
this change of multiplier in the KAM step.

For later use, the corresponding identity can be computed for any
affine contact map. Let $G=\rho^{-1}H\circ\Phi$ and
$K(q)=\Phi(q,0,0)$. Differentiating in $p$ and $z$ gives
\begin{equation}\label{eq:normaldampingtransform}
    \begin{aligned}
 &G_p|_0 =(DQ)^{-1}(H_p\circ K),\\
 &G_z|_0 =\rho^{-1}(H_p\circ K)\cdot DQ^{-T}D\rho +H_z\circ K =H_z\circ K+(G_p|_0)\cdot D\log\rho.
 \end{aligned}
\end{equation}
Here $D\log\rho$ denotes $\rho^{-1}D\rho$.
If the current angular coefficient is $G_p|_0=\omega+b$, the
identity reads
\[
    G_z|_0=H_z\circ K+D_\omega\log\rho
                         +b\cdot D\log\rho.
\]
Thus the term containing $b$ is present before the angular defect
has been removed. Once
$\mathcal J_1G=\omega\cdot p+\beta z$, this reduces to
\eqref{eq:intromultiplier}.
Since $\rho>0$ on the real torus, its real logarithm is periodic,
and averaging gives $\beta=\av{H_z\circ K}$.
For the selected torus with this frequency parameterization, the
multiplier changes the oscillatory coefficient but not its average.

\subsection{Parameter response in accumulated coordinates}

Fix an accumulated transformation $\Phi$ of the form
\eqref{eq:affine}, and let
$G=\rho^{-1}H_{\xi,e,\lambda,\varepsilon}\circ\Phi$.
Changing $(\xi,e)$ to $(\xi+\eta,e+\zeta)$ gives the exact variation
\[
    \Delta G
    =\rho^{-1}
       \bigl(H_{\xi+\eta,e+\zeta,\lambda,\varepsilon}
                      -H_{\xi,e,\lambda,\varepsilon}\bigr)
           \circ\Phi.
\]
If $p'$ denotes the momentum component of $\Phi$, this is
$\rho^{-1}(\eta\cdot p'-\zeta)$.
Substituting \eqref{eq:affine} yields
\[
    \rho^{-1}(\eta\cdot p'-\zeta) =\bigl((DQ)^{-1}\eta\bigr)\cdot p +\eta\cdot\rho^{-1}DQ^{-T}DU +z\,\eta\cdot DQ^{-T}D\log\rho-\rho^{-1}\zeta.
\]
Hence
\begin{equation}\label{eq:paramresponse}
    \Delta G=\ell\cdot\eta-\chi\zeta
                +(B\eta)\cdot p+(m\cdot\eta)z,
\end{equation}
where
\begin{equation}\label{eq:paramcoeff}
        B=(DQ)^{-1},\quad \chi=\rho^{-1},\quad \ell=\rho^{-1}DQ^{-T}DU, \quad m=DQ^{-T}D\log\rho.
\end{equation}
Here $B$ is an $n\times n$ matrix, $\ell$ and $m$ are column
vectors, and $\chi$ is scalar. All coefficients are evaluated at
the current angle $q$.

Measure the size of the accumulated transformation by
\begin{equation}\label{eq:Y}
    Y_s=\|Q-\mathrm{id}\|_{C^2(s)}
        +\|\rho-1\|_{C^2(s)}+\|U\|_{C^2(s)}.
\end{equation}
For sufficiently small $Y_s$,
\[
    \|B-I\|_s+\|\chi-1\|_s+\|\ell\|_s+\|m\|_s
    \leqslant C Y_s.
\]
Indeed, if $\|DQ-I\|_s\leqslant1/2$, the Neumann series gives
$\|(DQ)^{-1}\|_s\leqslant2$ and
$\|(DQ)^{-1}-I\|_s\leqslant2\|DQ-I\|_s$.
Similarly, $\|\rho-1\|_s\leqslant1/2$ gives
$\|\rho^{-1}\|_s\leqslant2$ and
$\|\rho^{-1}-1\|_s\leqslant2\|\rho-1\|_s$.
The remaining bounds follow from \eqref{eq:paramcoeff},
$D\log\rho=\rho^{-1}D\rho$, and the fixed dimensional norm
bounds for transposes.
At the identity triple $(Q,\rho,U)=(\mathrm{id},1,0)$,
one has $B=I$, $\chi=1$, and $\ell=m=0$.
Finally, \eqref{eq:paramresponse} is affine in $(p,z)$, so
$D^2_{(p,z)}(G+\Delta G)=D^2_{(p,z)}G$.
Thus the parameter update leaves the normal Hessian unchanged;
this will be used when deriving the linear correction equations.

\section{Linear correction and the parameter matrix}
\label{sec:linear}

\subsection{The first-jet equations}

In the current coordinates, write
\begin{equation}\label{eq:current}
    G=a(q)+(\omega+b(q))\cdot p+(\beta+c(q))z+R(q,p,z),
    \quad R=O_2.
\end{equation}
The first-jet defect is measured by
$E_s=\|a\|_s+\|b\|_s+\|c\|_s$.
Put $N=\omega\cdot p+\beta z$ and
$f=a+b\cdot p+cz$, so that $G=N+R+f$.
The normal Hessian blocks are
$S=R_{pp}|_0$, $t=R_{pz}|_0$, and $r_2=R_{zz}|_0$.
Here $t$ is an $n$-vector, not a flow parameter, and $r_2$ is a
scalar coefficient, not the working normal radius.

Choose a generating Hamiltonian
\begin{equation}\label{eq:generator}
    F(q,p,z)=u(q)+v(q)\cdot p+w(q)z,
    \quad \av{u}=\av{v}=\av{w}=0.
\end{equation}
The scalar functions $u,w$ and the vector function $v$ are to be
determined together with the parameter increments $\eta,\zeta$
and the normal-form increment $\delta\beta$.
By \eqref{eq:contact}, its contact flow $\Psi_\tau$ satisfies
\begin{equation}\label{eq:Fflow}
    \dot q=v(q),\quad
    \dot p=-Du-(Dv)^Tp-zDw-wp,\quad
    \dot z=-u-wz.
\end{equation}
The angular equation is independent of the normal variables, and
the vertical equation is affine in $z$. Together with the contact
identity, this places the local flow in the class
\eqref{eq:affine}. Write
$\Psi_\tau^*\alpha=\widehat\rho_\tau\alpha$.
For a holomorphic test Hamiltonian $J$, define
\begin{equation}\label{eq:Lie}
        \mathscr P_\tau J =\widehat\rho_\tau^{-1}J\circ\Psi_\tau, \quad \mathscr L_FJ=X_FJ+wJ,\quad \frac{d}{d\tau}\mathscr P_\tau J =\mathscr P_\tau\mathscr L_FJ =\mathscr L_F\mathscr P_\tau J.
\end{equation}
Indeed, $\partial_\tau\widehat\rho_\tau
=-(w\circ\Psi_\tau)\widehat\rho_\tau$, which gives the first
operator identity by differentiation. The second follows from the
local group property of the weighted pullbacks. These identities
hold on any common domain on which the flow is defined; the
uniform time-one domains are established in
Lemma~\ref{lem:flow}.

The action on the normal form is
\[
    \mathscr L_FN
    =-(D_\omega+\beta)u-(D_\omega v)\cdot p-(D_\omega w)z.
\]
The terms $-w\omega\cdot p$ and $-\beta wz$ in $X_FN$ cancel
the corresponding terms in $wN$.
For the nonlinear remainder, expansion of its derivatives gives
$R_p=Sp+tz+O_2$, $R_z=t\cdot p+r_2z+O_2$, and $R_q=O_2$.
Only the normal constant components $-Du$ and $-u$ of $X_F$
can therefore contribute to the first jet of $\mathscr L_FR$.
Using symmetry of $S$, one obtains
\begin{equation}\label{eq:LR}
    \mathcal J_1\mathscr L_FR
    =-(SDu+tu)\cdot p-(t\cdot Du+r_2u)z.
\end{equation}
These terms are linear in the generating coefficients. Since the
Hessian of $R$ is not assumed small, they must be included in the
correction equations.

The parameter update is the exact affine response
\eqref{eq:paramresponse}, with coefficients
\eqref{eq:paramcoeff}. We require
$\mathcal J_1[G+\Delta G+\mathscr L_F(N+R)]
=N+\delta\beta\,z$.
Equating its constant and normal-linear coefficients gives
\begin{equation}\label{eq:homological}
    \begin{aligned}
        (D_\omega+\beta)u
            &=a+\ell\cdot\eta-\chi\zeta,\\
        D_\omega v
            &=b+B\eta-SDu-tu,\\
        D_\omega w
            &=c+m\cdot\eta-t\cdot Du-r_2u-\delta\beta.
    \end{aligned}
\end{equation}
The updated system parameters are $(\xi+\eta,e+\zeta)$, and
the new normal-form coefficient is $\beta+\delta\beta$.
The term $\mathscr L_F(f+\Delta G)$ is not included in
\eqref{eq:homological}; it contains a defect coefficient and a
correction, and is estimated as part of the quadratic remainder in
Section~\ref{sec:iteration}.

\subsection{Elimination of the averages}

The first equation in \eqref{eq:homological}, together with
$\av u=0$, requires
$\av a+\av\ell\cdot\eta-\av\chi\,\zeta=0$.
For a sufficiently small $Y_s$, one has $\av\chi\geqslant1/2$.
Set
\begin{equation}\label{eq:meanelim}
        \zeta=\frac{\av a+\av\ell\cdot\eta}{\av\chi},\quad a^\circ=a-\chi\frac{\av a}{\av\chi}, \quad \ell^\circ=\ell-\chi\frac{\av\ell}{\av\chi}.
\end{equation}
Then $\av{a^\circ}=\av{\ell^\circ}=0$, and
Lemma~\ref{lem:fourier} gives
\begin{equation}\label{eq:usplit}
    u=u_0+V\eta,\quad
    u_0=T_\beta a^\circ,\quad
    V_j=T_\beta\ell_j^\circ,\quad 1\leqslant j\leqslant n.
\end{equation}
Here $V$ is the row $(V_1,\ldots,V_n)$, so
$V\eta=\sum_{j=1}^nV_j\eta_j$.
Taking the average of the second equation in
\eqref{eq:homological} yields the finite-dimensional system
\begin{equation}\label{eq:matrix}
        \mathsf M\eta=-\av{b-SDu_0-tu_0},\quad \mathsf M e_j=\av{Be_j-SDV_j-tV_j}, \quad 1\leqslant j\leqslant n.
\end{equation}
Once this system has been solved, \eqref{eq:meanelim} and
\eqref{eq:usplit} determine $\zeta$ and $u$. Set
\begin{equation}\label{eq:dbeta}
    \delta\beta=\av{c+m\cdot\eta-t\cdot Du-r_2u}.
\end{equation}
The remaining zero-average solutions are
\begin{equation}\label{eq:vwexplicit}
        v=T_0\bigl(b+B\eta-SDu-tu\bigr),\quad w=T_0\bigl(c+m\cdot\eta-t\cdot Du-r_2u-\delta\beta\bigr).
\end{equation}
The two right-hand sides have zero average by
\eqref{eq:matrix} and \eqref{eq:dbeta}, respectively.
Thus the first $n+1$ average equations determine the genuine system
increments $(\eta,\zeta)$; the final average equation determines
$\delta\beta$. No division by $\beta$ is used.

At the identity accumulated map, $B=I$ and $\ell=0$, so
$V=0$ and $\mathsf M=I$, independently of $S,t,r_2$.
The following estimates preserve this invertibility for small
accumulated transformations.

\subsection{Linear estimates}

\begin{lemma}
\label{lem:linear}
Assume $s\geqslant s_0/4$ and $\beta\in\R$.
Suppose that $a,b,c,S,t,r_2$ are bounded and holomorphic on
$\T_s^n$, real on the real torus, and that
$\|S\|_s+\|t\|_s+\|r_2\|_s\leqslant H_2$ for a fixed constant
$H_2$.
Let the parameter coefficients be given by
\eqref{eq:paramcoeff} for a triple real on the real torus.
There exist $\delta>0$ and $C>0$, depending on $H_2$ and the
fixed analytic and frequency data, such that $Y_s\leqslant\delta$
implies unique solvability of \eqref{eq:homological} with
$\av u=\av v=\av w=0$.
The constants $\eta,\zeta,\delta\beta$ are real, and the
functions are holomorphic on smaller strips and real on the real
torus. For $0<\sigma<\min\{1,s/12\}$, they satisfy
\begin{align}
    \|\mathsf M-I\|&\leqslant CY_s,\quad
       \|\mathsf M^{-1}\|\leqslant2,
       \label{eq:Minv}\\
    |\eta|+|\zeta|+|\delta\beta|&\leqslant CE_s,
       \label{eq:parambound}\\
    \|u\|_{s-\sigma}&\leqslant C\sigma^{-\kappa}E_s,
       \label{eq:ubound}\\
    \|v\|_{s-3\sigma}+\|w\|_{s-3\sigma}
       &\leqslant C\sigma^{-(2\kappa+1)}E_s.
       \label{eq:vwbound}
\end{align}
In particular,
\begin{equation}\label{eq:AF}
        A_F:=\|u\|_{C^2(s-5\sigma)} +\|v\|_{C^2(s-5\sigma)} +\|w\|_{C^2(s-5\sigma)} \leqslant C\sigma^{-(2\kappa+3)}E_s.
\end{equation}
The constants are independent of real $\beta$ and of $\sigma$.
\end{lemma}

\begin{proof}
The parameter bounds following \eqref{eq:Y} and
\eqref{eq:meanelim} give, for a small fixed $\delta$,
$\av\chi\geqslant1/2$, $\|a^\circ\|_s\leqslant CE_s$, and
$\|\ell^\circ\|_s\leqslant CY_s$.
Apply the fixed-margin estimate \eqref{eq:fixedmargin} with
distance $s_0/4$ and derivatives through order one. It gives
\begin{equation}\label{eq:realuv}
    \|u_0\|_{C^1(\T^n)}\leqslant CE_s,\quad
    \|V\|_{C^1(\T^n)}\leqslant CY_s.
\end{equation}
Only these real-torus estimates are needed in the averages.
For each column of the matrix in \eqref{eq:matrix},
\[
    |(\mathsf M-I)e_j|
    \leqslant\|B-I\|_{\T^n}
       +C H_2\bigl(\|DV_j\|_{\T^n}+\|V_j\|_{\T^n}\bigr)
    \leqslant CY_s.
\]
After accounting for the fixed matrix norm, choose $\delta$ so
that $\|\mathsf M-I\|\leqslant1/2$.
The Neumann series proves \eqref{eq:Minv}.
The right-hand side of \eqref{eq:matrix} is bounded by $CE_s$,
so $|\eta|\leqslant CE_s$.
Equations \eqref{eq:meanelim}, \eqref{eq:usplit}, and
\eqref{eq:realuv} then give
$|\zeta|+\|u\|_{C^1(\T^n)}\leqslant CE_s$.
Taking the average in \eqref{eq:dbeta} proves
$|\delta\beta|\leqslant CE_s$. This proves \eqref{eq:parambound}
without a loss depending on $\sigma$.

On complex strips, the source
$a^\circ+\ell^\circ\cdot\eta$ has norm at most $CE_s$.
The Fourier estimate \eqref{eq:Tbound} gives
\eqref{eq:ubound}; one additional Cauchy estimate gives
$\|Du\|_{s-2\sigma}\leqslant C\sigma^{-(\kappa+1)}E_s$.
The sources in \eqref{eq:vwexplicit} are consequently bounded by
$C\sigma^{-(\kappa+1)}E_s$ on $s-2\sigma$.
Applying $T_0$ with one further loss $\sigma$ proves
\eqref{eq:vwbound}. Two angular derivatives across the remaining
margin from $s-3\sigma$ to $s-5\sigma$ give \eqref{eq:AF}.
The corresponding estimate for $u$ has a smaller power of
$\sigma^{-1}$ and is absorbed in the same bound because
$\kappa>0$ and $0<\sigma<1$.

Finally, the average equations force \eqref{eq:meanelim},
\eqref{eq:matrix}, and \eqref{eq:dbeta} for every normalized
solution. The invertibility of $\mathsf M$ and uniqueness of the
zero-average Fourier inverses then imply uniqueness of all the
corrections. Real averaging, real matrix inversion, and
Lemma~\ref{lem:fourier} preserve reality throughout the construction.
\end{proof}

\section{A quantitative KAM step and convergence}
\label{sec:iteration}

\subsection{Fixed domains and the generating flow}

We first choose a normal radius $r>0$ and a geometric bound
$\bar\delta>0$ such that every affine map \eqref{eq:affine} with
$s\leqslant s_0/2$ and $Y_s\leqslant2\bar\delta$ satisfies
\begin{equation}\label{eq:backgrounddomain}
    \Phi(D_{s,8r})\subset D_{3s_0/4,r_0/2},\quad
    \|(DQ)^{-1}\|_s\leqslant2,\quad |\rho|\geqslant\tfrac12.
\end{equation}
For the norms fixed in Section~\ref{sec:main}, sufficient choices are
\[
    r\leqslant\frac{r_0}{128n},\quad
    \bar\delta\leqslant
    \min\left\{\frac18,\frac{s_0}{16},\frac{r_0}{64n}\right\}.
\]
Indeed, on $D_{s,8r}$,
\begin{align*}
 &|\operatorname{Im}Q|_\infty \leqslant s+2\bar\delta<3s_0/4,\quad |U+\rho z| \leqslant2\bar\delta+(1+2\bar\delta)8r<r_0/2,\\
 &|DQ^{-T}(\rho p+DU+zD\rho)|_\infty \leqslant\frac{n}{1-2\bar\delta} \bigl[8r(1+4\bar\delta)+2\bar\delta\bigr] <r_0/2.
 \end{align*}
The factor $n$ in the last line bounds the maximum-norm operator
norm of a transpose. The Neumann series for $DQ$ and the bound
$|\rho-1|\leqslant2\bar\delta$ give the other assertions in
\eqref{eq:backgrounddomain}.

For $|\xi-\omega|+|e|\leqslant2$, $|\varepsilon|\leqslant1$, and
$|\lambda|\leqslant\Lambda$, fix $B_0\geqslant1$ such that
\begin{equation}\label{eq:uniformB}
    \|\rho^{-1}H_{\xi,e,\lambda,\varepsilon}\circ\Phi\|_{s,8r}
    \leqslant B_0.
\end{equation}
This bound follows from \eqref{eq:backgrounddomain} and the original
analytic bounds. Cauchy's estimate across the fixed normal margin
from $8r$ to $4r$ gives uniform bounds for the normal derivatives
through order three. In particular, it gives a fixed bound for the
Hessian blocks $S,t,r_2$ in Section~\ref{sec:linear}.
We now choose $0<\delta\leqslant\bar\delta$ small enough for
Lemma~\ref{lem:linear} with these Hessian bounds. This last reduction
does not change $r$, $B_0$, or the derivative bounds already fixed.
All constants below are independent of the current strip loss and
of $\lambda\in[-\Lambda,\Lambda]$.

\begin{lemma}
\label{lem:flow}
Let $0<\sigma<\min\{1,s/12\}$, and let $F$ be the affine
Hamiltonian in \eqref{eq:generator}, with $A_F$ as in
\eqref{eq:AF}. There is a fixed $c_{\rm fl}>0$ such that, if
$A_F\leqslant c_{\rm fl}\min\{\sigma,r,1\}$, its contact flow
$\Psi_t$ is defined for $0\leqslant t\leqslant1$ on $D_{s-6\sigma,2r}$.
Write its affine triple as
$(\widehat Q_t,\widehat\rho_t,\widehat U_t)$. Then
\begin{equation}\label{eq:flowdomain}
    \Psi_t(D_{s-8\sigma,2r})\subset D_{s-7\sigma,4r},\quad
    \|\widehat\rho_t^{-1}\|_{s-8\sigma}\leqslant2,
\end{equation}
and
\begin{equation}\label{eq:flowC2}
    \|\widehat Q_t-\mathrm{id}\|_{C^2(s-6\sigma)}
    +\|\widehat\rho_t-1\|_{C^2(s-6\sigma)}
    +\|\widehat U_t\|_{C^2(s-6\sigma)}
    \leqslant CA_F.
\end{equation}
For an accumulated map $\Phi$ with $Y_s\leqslant\delta$, let
$\Phi_+=\Phi\circ\Psi_1$ and denote its triple by
$(Q_+,\rho_+,U_+)$. One has
\begin{equation}\label{eq:compositionbound}
        \|Q_+-Q\|_{C^2(s-12\sigma)} +\|\rho_+-\rho\|_{C^2(s-12\sigma)} +\|U_+-U\|_{C^2(s-12\sigma)} \leqslant C\sigma^{-1}A_F.
\end{equation}
These estimates do not require the generating coefficients to have
zero average.
\end{lemma}

\begin{proof}
The equations \eqref{eq:Fflow} give
\begin{equation}\label{eq:tripleODE}
    \dot{\widehat Q}_t=v(\widehat Q_t),\quad \widehat Q_0=\mathrm{id}, \quad \dot{\widehat\rho}_t=-w(\widehat Q_t)\widehat\rho_t,\quad \widehat\rho_0=1, \quad \dot{\widehat U}_t=-u(\widehat Q_t) -w(\widehat Q_t)\widehat U_t,\quad \widehat U_0=0.
\end{equation}
Since $\|v\|_{s-5\sigma}\leqslant A_F$, a trajectory starting in
$\T^n_{s-6\sigma}$ has displacement at most $A_F$ before it
leaves $\T^n_{s-5\sigma}$. Taking $A_F<\sigma/2$ excludes such
an exit for $0\leqslant t\leqslant1$. The solutions depend holomorphically
on the initial angle.

For $\mathsf P_t=D\widehat Q_t$ and
$\mathsf Z_t=D^2\widehat Q_t$, the derivative equations are
\[
    \dot{\mathsf P}_t=Dv(\widehat Q_t)\mathsf P_t,
    \quad
    \dot{\mathsf Z}_t=Dv(\widehat Q_t)\mathsf Z_t
       +D^2v(\widehat Q_t)[\mathsf P_t,\mathsf P_t],
\]
with $\mathsf P_0=I$ and $\mathsf Z_0=0$.
Gronwall's inequality gives
$\|\mathsf P_t-I\|+\|\mathsf Z_t\|\leqslant CA_F$ for small
$A_F$. The other two components have the explicit expressions
\begin{equation}\label{eq:rhoUexplicit}
        \widehat\rho_t =\exp\left(-\int_0^t w(\widehat Q_a)\,da\right),\quad \widehat U_t =-\widehat\rho_t \int_0^t\widehat\rho_a^{-1}u(\widehat Q_a)\,da.
\end{equation}
Differentiating these formulas twice in the initial angle and using
the bounds for $\mathsf P_t$ and $\mathsf Z_t$ proves
\eqref{eq:flowC2}. The affine momentum formula then gives, for
$|p|_\infty,|z|<2r$,
\[
    |p_t-p|_\infty+|z_t-z|\leqslant C(1+r)A_F.
\]
A sufficiently small choice of $c_{\rm fl}$ makes this bound
smaller than the available normal margin. Together with the angular
displacement bound and \eqref{eq:rhoUexplicit}, it proves
\eqref{eq:flowdomain}. The affine equations also give existence of
the full flow on the stated domain. Reality of the coefficients
implies reality of the flow, and $\widehat\rho_t>0$ on real points.

It remains to estimate composition. Put
$h=\widehat Q_1-\mathrm{id}$. For a periodic coefficient $f$ of
the old triple,
\[
    f(q+h(q))-f(q)
    =\int_0^1 Df(q+ah(q))h(q)\,da.
\]
The chain rule through order two implies
$\|f\circ(\mathrm{id}+h)-f\|_{C^2}
 \leqslant C\|f\|_{C^3}\|h\|_{C^2}$ on domains containing these
segments. On the target strip $s-12\sigma$, all the segments lie
in $s-11\sigma$. Cauchy's estimate therefore bounds the required
$C^3$ norm by $C\sigma^{-1}\|f\|_{C^2(s)}$.
For $Q=\mathrm{id}+k$, apply this estimate to the periodic part
$k$, treating the identity part separately. Formula
\eqref{eq:composition} gives
\begin{align*}
 &Q_+-Q=h+k\circ\widehat Q_1-k,\quad \rho_+-\rho =(\rho\circ\widehat Q_1-\rho) +(\rho\circ\widehat Q_1)(\widehat\rho_1-1),\\
 &U_+-U=(U\circ\widehat Q_1-U) +(\rho\circ\widehat Q_1)\widehat U_1.
 \end{align*}
The product rule, $Y_s\leqslant\delta$, and \eqref{eq:flowC2} now
yield \eqref{eq:compositionbound}.
\end{proof}

\subsection{The one-step lemma}

\begin{lemma}\label{lem:step}
Suppose that
$G=\rho^{-1}H_{\xi,e,\lambda,\varepsilon}\circ\Phi$
is defined on $D_{s,8r}$, has the decomposition \eqref{eq:current},
and satisfies \eqref{eq:uniformB}. Assume that all data are real on
real points and that
\[
    s_0/4\leqslant s-12\sigma<s\leqslant s_0/2,\quad
    0<\sigma<1,\quad Y_s\leqslant\delta,
\]
with $|\beta|\leqslant\Lambda+1$, $|\xi-\omega|+|e|\leqslant1$,
$|\lambda|\leqslant\Lambda$, and $|\varepsilon|\leqslant1$.
Set $\ell_*=2\kappa+6$ and $m_*=4\kappa+16$.
There are fixed constants $c_*>0$ and $C_*\geqslant1$ such that, if
\[
    E_s\leqslant c_*\sigma^{\ell_*},
\]
the correction equations of Section~\ref{sec:linear} determine a
parameter increment $(\eta,\zeta)$, a coefficient increment
$\delta\beta$, and a generating flow with time-one map $\Psi_1$.
Let $s_+=s-12\sigma$, $\beta_+=\beta+\delta\beta$, and
$\Phi_+=\Phi\circ\Psi_1$. The updated Hamiltonian is defined
on $D_{s_+,8r}$ by
$G_+=\rho_+^{-1}H_{\xi+\eta,e+\zeta,\lambda,\varepsilon}
\circ\Phi_+$ and satisfies \eqref{eq:uniformB} there.
Its first-jet defect relative to $\omega\cdot p+\beta_+z$ obeys
\begin{align}
    E_+&\leqslant C_*\sigma^{-m_*}E_s^2,
        \label{eq:quadratic}\\
    |\eta|+|\zeta|+|\delta\beta|&\leqslant C_*E_s,
        \label{eq:stepchanges}\\
    Y_{s_+}(\Phi_+)&\leqslant Y_s(\Phi)
                       +C_*\sigma^{-\ell_*}E_s.
        \label{eq:Ystep}
\end{align}
The accumulated triple also satisfies the increment estimate
\begin{equation}\label{eq:tripleincrement}
        \|Q_+-Q\|_{C^2(s_+)} +\|\rho_+-\rho\|_{C^2(s_+)} +\|U_+-U\|_{C^2(s_+)} \leqslant C_*\sigma^{-\ell_*}E_s.
\end{equation}
In addition, $Y_{s_+}(\Phi_+)\leqslant2\delta$ and
$|\xi+\eta-\omega|+|e+\zeta|\leqslant2$.
The next step is available whenever the updated quantities satisfy
the smaller bounds in the hypotheses.
\end{lemma}

\begin{proof}
Let $C_{\rm lin}$ and $C_{\rm par}$ be fixed constants in
\eqref{eq:AF} and \eqref{eq:parambound} of
Lemma~\ref{lem:linear}, and let $C_{\rm comp}$ be the constant
in \eqref{eq:compositionbound}. Since
$\ell_*-(2\kappa+3)=3$, the smallness assumption gives
$A_F\leqslant C_{\rm lin}c_*\sigma^3$.
Choose $c_*$ small enough that
\[
    C_{\rm lin}c_*\leqslant c_{\rm fl}\min\{r,1\},\quad
    C_{\rm comp}C_{\rm lin}c_*\leqslant\delta,\quad
    C_{\rm par}c_*\leqslant1.
\]
Then Lemma~\ref{lem:flow} applies, the updated system parameters
remain in the larger range used in \eqref{eq:uniformB}, and
\[
    Y_{s_+}(\Phi_+)
    \leqslant\delta+C_{\rm comp}C_{\rm lin}c_*\sigma^2
    \leqslant2\delta.
\]

Write $N=\omega\cdot p+\beta z$ and
$f=a+b\cdot p+cz$, so $G=N+R+f$ with $R=O_2$.
After the parameter change, put
$G^\sharp=G+\Delta G=N+R+\widetilde g$, where
\[
    \widetilde g
    =a+\ell\cdot\eta-\chi\zeta
       +(b+B\eta)\cdot p+(c+m\cdot\eta)z.
\]
The sum of the norms of these three coefficients is at most $CE_s$.
Since the parameter change is made in the original family,
$\|G^\sharp\|_{s,8r}\leqslant B_0$ follows from the larger parameter
range and \eqref{eq:uniformB}.

For the generating flow, write
$\mathscr P_tJ=\widehat\rho_t^{-1}J\circ\Psi_t$ and
$\mathscr L_FJ=X_FJ+wJ$, as in \eqref{eq:Lie}.
The correction equations give
\begin{equation}\label{eq:exactcancel}
    \mathcal J_1\bigl[G^\sharp+\mathscr L_F(N+R)\bigr]
    =\omega\cdot p+(\beta+\delta\beta)z.
\end{equation}
Integrating $\partial_t\mathscr P_tJ=\mathscr P_t\mathscr L_FJ$
twice gives the exact identity
\[
    G_+=\mathscr P_1G^\sharp
    =G^\sharp+\mathscr L_FG^\sharp
      +\int_0^1(1-t)\mathscr P_t\mathscr L_F^2G^\sharp\,dt
\]
on the smaller domain where the flow is estimated.
Subtracting \eqref{eq:exactcancel}, the new first-jet defect is
\begin{equation}\label{eq:twoerrors}
    f_+
    =\mathcal J_1\mathscr L_F\widetilde g
      +\mathcal J_1\int_0^1(1-t)
                    \mathscr P_t\mathscr L_F^2G^\sharp\,dt.
\end{equation}

To estimate the first term, write
$\widetilde g=\widetilde a+\widetilde b\cdot p+\widetilde c z$.
Its image under $\mathscr L_F$ is affine, with coefficients
\[
        d_0=v\cdot D\widetilde a-Du\cdot\widetilde b -u\widetilde c+w\widetilde a,\quad d_p=D\widetilde b\,v-Dv\,\widetilde b,\quad d_z=v\cdot D\widetilde c-\widetilde b\cdot Dw.
\]
Cauchy's estimate for the defect coefficients gives
$\|d_0\|_{s-6\sigma}+\|d_p\|_{s-6\sigma}
 +\|d_z\|_{s-6\sigma}\leqslant C\sigma^{-1}A_FE_s$.

For the integral term, the expression
\[
    \mathscr L_FJ
    =v\cdot J_q
      -[Du+(Dv)^Tp+zDw+wp]\cdot J_p
      -(u+wz)J_z+wJ
\]
gives, on the domains used below,
\[
    \|\mathscr L_FJ\|_{a-\sigma,b'}
    \leqslant CA_F\bigl(1+\sigma^{-1}+(b-b')^{-1}\bigr)
             \|J\|_{a,b},
\]
where $0<\sigma<a\leqslant s-5\sigma$ and $0<b'<b\leqslant8r$.
The constant incorporates the fixed bound on the normal coordinates.
Applying this estimate successively from $D_{s-5\sigma,8r}$ to
$D_{s-6\sigma,6r}$ and then to $D_{s-7\sigma,4r}$ gives
\[
    \|\mathscr L_FG^\sharp\|_{s-6\sigma,6r} \leqslant C\sigma^{-1}A_FB_0,\quad \|\mathscr L_F^2G^\sharp\|_{s-7\sigma,4r} \leqslant C\sigma^{-2}A_F^2B_0.
\]
The second application differentiates the entire function
$\mathscr L_FG^\sharp$, so it also accounts for derivatives of
the generating coefficients.
By \eqref{eq:flowdomain}, the integral in \eqref{eq:twoerrors}
is bounded by $C\sigma^{-2}A_F^2B_0$ on $D_{s-8\sigma,2r}$.
Taking its first jet costs only a fixed normal-radius factor:
\[
    \|J|_0\|_a+\|J_p|_0\|_a+\|J_z|_0\|_a
    \leqslant\left(1+\frac{n+1}{2r}\right)\|J\|_{a,2r}.
\]
After restricting to $s_+$, we obtain
\begin{equation}\label{eq:sharpstep}
    E_+\leqslant C\bigl(\sigma^{-1}A_FE_s+\sigma^{-2}A_F^2\bigr)
        \leqslant C\sigma^{-(4\kappa+8)}E_s^2.
\end{equation}
Since $0<\sigma<1$, this implies \eqref{eq:quadratic} with the
chosen $m_*$. Lemma~\ref{lem:linear} gives
\eqref{eq:stepchanges}. Combining \eqref{eq:compositionbound}
with \eqref{eq:AF} bounds the triple increment by
$C\sigma^{-(2\kappa+4)}E_s$.
Since $2\kappa+4<\ell_*$ and $\sigma<1$, this proves
\eqref{eq:tripleincrement}; the triangle inequality gives
\eqref{eq:Ystep}. Increase the fixed $C_*$ to cover all these
bounds.

It remains to recover the domain needed for the next step.
The new triple is holomorphic on $\T^n_{s_+}$ and has size at most
$2\delta\leqslant2\bar\delta$. Its affine formula therefore defines
$\Phi_+$ on $D_{s_+,8r}$ with image in the original domain.
On that larger domain set
\[
    G_+=\rho_+^{-1}
          H_{\xi+\eta,e+\zeta,\lambda,\varepsilon}\circ\Phi_+.
\]
By \eqref{eq:composition}, this definition agrees with
$\mathscr P_1G^\sharp$ on their common domain. In particular, it
has the first-jet estimate already proved and the bound
\eqref{eq:uniformB}. Define
$a_+=G_+|_0$, $b_+=G_{+,p}|_0-\omega$,
$c_+=G_{+,z}|_0-\beta_+$, and
$R_+=G_+-\mathcal J_1G_+$.
Then $R_+=O_2$, and its normal Hessian has the same fixed bounds
as that of $G_+$. The larger normal domain is thus recovered from
the accumulated map and the original Hamiltonian, rather than by
requiring $\Psi_1$ to preserve the radius $8r$.
\end{proof}

\subsection{Iteration and convergence}

Use $s^{(j)}$ for the widths of the successive angular strips, so
that $s_0$ continues to denote the original width. Set
\begin{equation}\label{eq:strips}
    \sigma_0=\min\{1/16,s_0/256\},\quad
    \sigma_j=\sigma_0 2^{-j},\quad
    s^{(j)}=s_0/2-12\sum_{i=0}^{j-1}\sigma_i.
\end{equation}
Then $s^{(j)}\searrow s_\infty=s_0/2-24\sigma_0>s_0/4$.
Initialize with
\[
    \Phi_0=\mathrm{id},\quad
    \xi_0=\omega,\quad e_0=0,\quad\beta_0=\lambda.
\]
Thus $G_0=H_{\omega,0,\lambda,\varepsilon}$ and
\[
    a_0=\varepsilon P|_0,\quad
    b_0=\varepsilon P_p|_0,\quad
    c_0=\varepsilon P_z|_0.
\]
The initial normal remainder is
$\mathcal R_0=R_0+\varepsilon(P-\mathcal J_1P)$, rather than just
$R_0$. Cauchy's estimate in the original normal domain gives
$E_0\leqslant C_0|\varepsilon|$ with a fixed $C_0\geqslant1$.

At each step, let $\Psi_j$ be the time-one map supplied by
Lemma~\ref{lem:step}, and let its multiplier be
$\widehat\rho_j$. The updates are
\[
    \begin{aligned}
        \xi_{j+1}&=\xi_j+\eta_j,&
        e_{j+1}&=e_j+\zeta_j,&
        \beta_{j+1}&=\beta_j+\delta\beta_j,\\
        \Phi_{j+1}&=\Phi_j\circ\Psi_j,&
        G_{j+1}&=\widehat\rho_j^{-1}
                    (G_j+\Delta G_j)\circ\Psi_j.
    \end{aligned}
\]
In particular,
$G_j=\rho_j^{-1}H_{\xi_j,e_j,\lambda,\varepsilon}\circ\Phi_j$
throughout. Write $E_j$ for its first-jet defect on $s^{(j)}$.
If $E_0=0$, all corrections can be taken to be zero. Suppose
henceforth that $E_0>0$.

Put $\mathcal A=C_*\sigma_0^{-m_*}$ and
$\mathcal B=2^{m_*}$. On every finite segment of the construction,
\eqref{eq:quadratic} gives
$E_{j+1}\leqslant\mathcal A\mathcal B^jE_j^2$.
For $d_j=\mathcal A\mathcal B^{j+1}E_j$, this becomes
\begin{equation}\label{eq:recursion}
    d_{j+1}\leqslant d_j^2,\quad
    E_j\leqslant\frac{d_0^{2^j}}{\mathcal A\mathcal B^{j+1}}
         =E_0\mathcal B^{-j}d_0^{2^j-1}.
\end{equation}
If $d_0\leqslant1/2$, the numerical majorant on the right has a
summable weighted series. More precisely, for each fixed $L\geqslant0$
define
\[
    W_L=\sigma_0^{-L}
          \sum_{j=0}^\infty
          2^{(L-m_*)j}2^{-(2^j-1)}<+\infty.
\]
For every constructed finite segment,
\begin{equation}\label{eq:weighted}
    \sum_{j=0}^{N}\sigma_j^{-L}E_j\leqslant W_LE_0.
\end{equation}
All the $W_L$ are fixed before iteration.

Choose
\begin{equation}\label{eq:epsilonbudget}
        E_{\rm crit} =\min\left\{ \frac1{2\mathcal A\mathcal B},\, \frac{c_*}{W_{\ell_*}},\, \frac{\delta}{2C_*W_{\ell_*}},\, \frac1{2C_*W_0} \right\},\quad \varepsilon_*\leqslant\min\{1,E_{\rm crit}/C_0\}.
\end{equation}
We verify the induction before taking an infinite composition.
Assume that the quantities at stage $j$ have been constructed.
The recurrence from the preceding steps and \eqref{eq:weighted}
give
\[
    \begin{aligned}
 &Y_{s^{(j)}}(\Phi_j) \leqslant C_*\sum_{i=0}^{j-1}\sigma_i^{-\ell_*}E_i \leqslant\delta/2,\\
 &|\xi_j-\omega|+|e_j|+|\beta_j-\lambda| \leqslant C_*\sum_{i=0}^{j-1}E_i\leqslant1/2,\quad \sigma_j^{-\ell_*}E_j\leqslant W_{\ell_*}E_0\leqslant c_*.
 \end{aligned}
\]
The domain condition follows from \eqref{eq:strips}. Hence all
hypotheses of Lemma~\ref{lem:step} hold, and it constructs stage
$j+1$. This proves existence of the entire sequence. Letting
$N\to\infty$ in \eqref{eq:weighted}, and using
\eqref{eq:stepchanges} and \eqref{eq:tripleincrement}, gives
summability of the parameter and triple increments. In particular,
\[
        \sum_{j=0}^{\infty}\bigl(\|Q_{j+1}-Q_j\|_{C^2(s_\infty)} +\|\rho_{j+1}-\rho_j\|_{C^2(s_\infty)} +\|U_{j+1}-U_j\|_{C^2(s_\infty)}\bigr) \leqslant C_*W_{\ell_*}E_0.
\]
This is a Cauchy estimate for the accumulated triples, not merely
a bound on their total size. Its right-hand side is
$O(|\varepsilon|)$ uniformly in $\lambda$.

\subsection{The limiting torus and parameter dependence}

The parameter sequences converge to $(\xi,e,\beta)$, and the
triples converge in $C^2$ on $\T^n_{s_\infty}$ to $(Q,\rho,U)$.
The limits are holomorphic and real on the real torus. On every
fixed smaller strip, Cauchy's estimate gives convergence of all
fixed-order angular derivatives. In particular,
\[
    |\xi-\omega|+|e|+|\beta-\lambda|
    +\|Q-\mathrm{id}\|_{C^2(s_\infty)}
    +\|\rho-1\|_{C^2(s_\infty)}
    +\|U\|_{C^2(s_\infty)}
    \leqslant C|\varepsilon|.
\]
Take $s_*=s_0/8$ and $r_*=r$, decreasing $\varepsilon_*$ when
necessary in what follows.

Write $Q=\mathrm{id}+h$. We may assume
$\|Dh\|_{s_\infty}<1/2$ and $\|\rho-1\|_{s_\infty}<1/2$.
On the universal covering strip, integration along line segments
gives
$|Q(q)-Q(q')|_\infty\geqslant\tfrac12|q-q'|_\infty$.
Periodicity then proves injectivity modulo $2\pi\Z^n$.
The real restriction is a local diffeomorphism of degree one and
hence a diffeomorphism of $\T^n$.
To obtain an inverse on a fixed complex strip, choose $d>0$ with
$s_*+d<s_\infty$ and reduce $\varepsilon_*$ so that
$\|h\|_{s_\infty}<d/2$.
For $q'\in\T^n_{s_*}$, the map $q\mapsto q'-h(q)$ is a
contraction on the closed ball $|q-q'|_\infty\leqslant d$ in the
covering space. Its fixed point depends holomorphically on $q'$
and defines $Q^{-1}$. The inverse is compatible with the lattice
translations.

The affine formula defines $\Phi$ on $D_{s_*,r_*}$, with image
in $D_{s_0,r_0}$. For a point in its image, the full inverse is
\[
        q=Q^{-1}(q'),\quad z=\frac{z'-U(q)}{\rho(q)},\quad p=\rho(q)^{-1} \bigl[DQ(q)^Tp'-DU(q)-zD\rho(q)\bigr].
\]
Thus $\Phi$ is an embedding and $\rho>0$ on the real torus.
The fixed gap $s_\infty-s_*$ bounds the third angular derivatives
of the limiting triple by $C|\varepsilon|$. These are the derivatives
needed to estimate the full-space $C^2$ norm of the momentum
component of \eqref{eq:affine}. The product and inverse rules
therefore give \eqref{eq:mainbound}.

On a fixed smaller domain, the maps $\Phi_j$ have images in a
common compact subset of the original analytic domain, modulo the
real angular lattice. The parameter and map convergence therefore
imply
\[
    G_j\to
    G_\infty=\rho^{-1}H_{\xi,e,\lambda,\varepsilon}\circ\Phi.
\]
Cauchy's estimate passes this convergence to the first normal jet.
Since $E_j\to0$ and $\beta_j\to\beta$, one obtains
$\mathcal J_1G_\infty=\omega\cdot p+\beta z$.
This proves \eqref{eq:target}. The zero torus for $G_\infty$
carries the vector field $(\omega,0,0)$, and the contact pullback
identity gives $X_H\circ K=DK\,\omega$ for
\[
    K(\theta)
    =\bigl(Q(\theta),DQ(\theta)^{-T}DU(\theta),U(\theta)\bigr).
\]
The restriction of $\alpha$ to the zero torus is zero, so
$K^*\alpha=0$. Evaluation of the normal form there gives
$H\circ K=0$. For $\mathfrak u=U\circ Q^{-1}$, the chain rule
gives $D\mathfrak u(Q(\theta))=DQ(\theta)^{-T}DU(\theta)$;
hence the image of $K$ is the full Legendrian graph
$j^1\mathfrak u$.

Using $H_p\circ K=DQ\,\omega$ in
\eqref{eq:normaldampingtransform} gives
\begin{equation}\label{eq:meanidentity}
    \beta=H_z\circ K
            +(H_p\circ K)\cdot DQ^{-T}D\log\rho
          =H_z\circ K+D_\omega\log\rho.
\end{equation}
The real logarithm of $\rho$ is periodic. Its directional derivative
has zero average, proving $\beta=\av{H_z\circ K}$, with the
average taken in the parameter angle $\theta$.

\begin{lemma}
\label{lem:finitecontinuity}
Fix $\varepsilon$ satisfying \eqref{eq:epsilonbudget} and use the
zero-average normalization at every step. For each finite $j$, the
parameters $(\xi_j,e_j,\beta_j)$ are continuous in real
$\lambda\in[-\Lambda,\Lambda]$. The accumulated triple is
continuous in $C^2$ on every strip strictly narrower than
$\T^n_{s^{(j)}}$, and $G_j$ is continuous in the supremum norm
on each $D_{a,b}$ with $a<s^{(j)}$ and $b<8r$.
\end{lemma}

\begin{proof}
The initial quantities have this property. Assume it at step $j$.
By Cauchy's estimate across a positive margin, the coefficients of
the first jet, the Hessian blocks, and the parameter response are
continuous on each smaller angular strip. The joint estimate
\eqref{eq:Tjointcontinuous} handles both the variation of the source
and that of $\beta_j$. It shows that $u_0$ and $V$ depend
continuously on $\lambda$, with their first angular derivatives
continuous on the real torus. Averaging is a bounded linear
operation, and matrix inversion is continuous on the uniformly
invertible matrices in \eqref{eq:Minv}. Equations
\eqref{eq:meanelim}--\eqref{eq:vwexplicit} therefore give
continuous increments and generating coefficients, with their
required derivatives, on smaller strips.

For a fixed smaller target at stage $j+1$, the generating
trajectories stay inside a common coefficient domain with positive
margin, by Lemma~\ref{lem:flow}. Subtracting the integral equations
\eqref{eq:tripleODE} for two damping values and using the uniform
bound on the first derivatives of the coefficients gives continuous
dependence of the flows by Gronwall's inequality. Cauchy's estimate
on a further, fixed margin gives continuity of their first two
angular derivatives. The composition formulas then prove continuity
of the next triple. The definition of $G_{j+1}$ as a pullback of
the original family proves its asserted continuity on smaller full
space domains. Every operation at a finite stage has a strictly
positive analytic margin; no continuity in the norm on an open
boundary strip is being assumed. This completes the induction.
\end{proof}

To finish the parameter assertion of Theorem~\ref{thm:main}, fix
an intermediate strip $s_*<a<s_\infty$. The increment estimate
\eqref{eq:tripleincrement} and \eqref{eq:weighted} bound the tails
of the triples in $C^2(s_\infty)$ uniformly in $\lambda$.
The numerical majorant in \eqref{eq:recursion} makes these bounds
tend to zero as the starting index tends to infinity, also uniformly
in $\lambda$. Cauchy's estimate transfers them to $C^3(a)$.
The same argument applies to the parameter tails with weight zero.
Given a tolerance, choose one finite stage so that both uniform tails
are smaller than a third of it. Lemma~\ref{lem:finitecontinuity}
controls the difference of that finite stage at two nearby damping
values. The triangle inequality gives continuity of the limit.
The affine formula, together with uniform invertibility of $DQ$
and nonvanishing of $\rho$, transfers this to the full-space
$C^2(D_{s_*,r_*})$ norm of $\Phi$.
This completes the proof of Theorem~\ref{thm:main}.
The argument also applies at $\lambda=0$ and requires no common
complex neighborhood in the damping parameter. The corresponding
output $\beta_{0,\varepsilon}$ is $O(|\varepsilon|)$, not
necessarily zero.

\section{Mean damping and parameter compensation}
\label{sec:extensions}

\subsection{Prescribing the mean damping}

\begin{corollary}\label{cor:prescribed}
Under the hypotheses of Theorem~\ref{thm:main}, consider
$H=\xi\cdot p+(\lambda+d)z+R_0+\varepsilon P-e$ with an
additional real parameter $d$.
For $\lambda\in[-\Lambda/2,\Lambda/2]$ and sufficiently small
real $\varepsilon$, there are $\xi,e,d$ and an analytic contact
embedding $\Phi$, with multiplier $\rho$, such that
\[
    |\xi-\omega|+|e|+|d|\leqslant C|\varepsilon|,\quad
    \rho^{-1}H\circ\Phi=\omega\cdot p+\lambda z+O_2.
\]
The domain and the constants are uniform on this interval, and
$\Phi$ has the same near-identity bounds as in
Theorem~\ref{thm:main}.
\end{corollary}

\begin{proof}
Denote the input damping in Theorem~\ref{thm:main} by $t$.
For fixed $\varepsilon$, the construction gives a continuous
function $\beta(t)=t+b(t)$ with $|b(t)|\leqslant C_d|\varepsilon|$,
where $C_d\geqslant1$ is fixed. For $\varepsilon\ne0$, set
$t_\pm=\lambda\pm2C_d|\varepsilon|$.
Reducing the perturbation threshold ensures
$t_\pm\in[-\Lambda,\Lambda]$ and gives
$\beta(t_-)<\lambda<\beta(t_+)$.
The intermediate value theorem supplies a $t_*$ between these
points with $\beta(t_*)=\lambda$. Take $d=t_*-\lambda$ and
use the parameters and embedding corresponding to $t_*$.
For $\varepsilon=0$, take $d=e=0$, $\xi=\omega$, and
$\Phi=\mathrm{id}$.
\end{proof}

This argument gives existence and uniform bounds, but not uniqueness
or continuity of a selected $d$.
The additional parameter cannot in general be replaced by $\xi$
or $e$: for $H=\xi\cdot p+(\lambda+\varepsilon)z-e$,
\eqref{eq:meanidentity} forces $\beta=\lambda+\varepsilon$,
whereas adding $d=-\varepsilon$ gives $\beta=\lambda$.

\subsection{Fixed angle-dependent damping}

\begin{proposition}\label{prop:variable}
Let $c_0$ be a fixed real-analytic function on $\T^n$ with
$\av{c_0}=0$, and consider
\[
    H=\xi\cdot p+(\lambda+c_0(q))z
        +R_0(q,p,z)+\varepsilon P(q,p,z)-e.
\]
Under the remaining hypotheses of Theorem~\ref{thm:main}, there
are selected parameters and an invariant Legendrian torus with the
same prescribed frequency, uniformly for
$\lambda\in[-\Lambda,\Lambda]$. The parameters satisfy
$|\xi-\omega|+|e|+|\beta-\lambda|\leqslant C|\varepsilon|$
and can be chosen, together with the torus parameterization,
continuously in real $\lambda$.

On a fixed smaller domain, the normalizing contact map is an
$O(|\varepsilon|)$ perturbation in $C^2$ of a fixed map $\Psi_g$
defined below. Its restriction $K_\varepsilon$ to the zero torus
satisfies $K_\varepsilon(\theta)-(\theta,0,0)=O(|\varepsilon|)$
in $C^2$ on a fixed smaller angular strip.
The domains and constants may depend on $c_0$ and its analytic
extension.
\end{proposition}

\begin{proof}
Choose a complex strip on which $c_0$ is holomorphic and bounded,
and work on a smaller strip also contained in the original angular
domain. Lemma~\ref{lem:fourier} gives a real-analytic solution
$g=T_0c_0$ of $D_\omega g=c_0$, with $\av g=0$.
Set
\[
    \Psi_g(q,I,W)
    =\bigl(q,e^{-g(q)}(I-WDg(q)),e^{-g(q)}W\bigr).
\]
As in \eqref{eq:dampingpullback},
$\Psi_g^*\alpha=e^{-g}\alpha$.
Its inverse is
$(q,p,z)\mapsto(q,e^g(p+zDg),e^gz)$.
The function $g$ need not be small. After fixing a further angular
margin and a sufficiently small normal radius, both maps and all
coefficients needed below have fixed analytic bounds.

Put $h=e^g$, $\widetilde R_0=e^gR_0\circ\Psi_g$, and
$\widetilde P=e^gP\circ\Psi_g$. Direct substitution gives
\begin{equation}\label{eq:normalizedvar}
        \widetilde H:=e^gH\circ\Psi_g =\omega\cdot I+\lambda W +\widetilde R_0+\varepsilon\widetilde P +(\xi-\omega)\cdot I -(\xi-\omega)\cdot Dg\,W-e h.
\end{equation}
Here $\widetilde R_0=O_2$, because $\Psi_g$ is linear in the
normal variables and fixes their origin.
The parameter dependence in \eqref{eq:normalizedvar} differs from
\eqref{eq:family}; we therefore check the average equations for
this transformed family.

At the identity accumulated map, its parameter response has
$B=I$, $\chi=h$, $\ell=0$, and $m=-Dg$.
For a subsequent near-identity triple $(Q,\rho,U)$, the response
still has the affine form \eqref{eq:paramresponse}, now with
\[
    \begin{aligned}
        B&=(DQ)^{-1},&
        \chi&=\rho^{-1}e^{g\circ Q},\\
        \ell&=\rho^{-1}DQ^{-T}DU
                   -\rho^{-1}U\,(Dg)\circ Q,&
        m&=DQ^{-T}D\log\rho-(Dg)\circ Q.
    \end{aligned}
\]
Here $(Dg)\circ Q$ is the gradient of $g$ evaluated at $Q(q)$,
not the derivative of the composite function $g\circ Q$.
These formulas follow by pulling back
$\eta\cdot I-\eta\cdot Dg\,W-\zeta h$.

On the fixed transformed domain,
$\|B-I\|_s+\|\ell\|_s+\|\chi-h\|_s\leqslant CY_s$ and
$\|m\|_s\leqslant C$. Since $h_*:=\min_{\T^n}h>0$, a small
fixed bound on $Y_s$ ensures $\av\chi\geqslant h_*/2$.
Eliminate $\zeta$ using the first average equation and set
$\ell^\circ=\ell-\chi\av\ell/\av\chi$.
Then $\|\ell^\circ\|_s\leqslant CY_s$, so the fixed-margin
estimate for $V_j=T_\beta\ell_j^\circ$ gives
$\|V\|_{C^1(\T^n)}\leqslant CY_s$.
Consequently, the reduced angular matrix satisfies
$\|\mathsf M-I\|\leqslant CY_s$. At the identity map it is
exactly $I$, since $\ell=V=0$.
The scalar and angular average equations are therefore solvable
with increments of size $CE_s$.
In the last average equation, $m\cdot\eta$ is still $O(E_s)$
because $m$ is bounded; its smallness is not required.
The same Fourier estimates then give \eqref{eq:AF} on smaller
strips.

The transformed family has fixed analytic bounds, and its parameter
response remains affine in the normal variables. Thus the Hessian
bounds, generating-flow estimates, and quadratic remainder proof
apply with the initial widths replaced by the fixed transformed
widths. The larger normal domain is recovered from
\eqref{eq:normalizedvar} at each step, exactly as in the proof of
Lemma~\ref{lem:step}. The iteration therefore produces a
near-identity map $\Phi_\varepsilon$ in these coordinates, with
the asserted parameter estimates and continuous dependence on
real $\lambda$.

The total normalizing map is
$\Phi_{\rm tot}=\Psi_g\circ\Phi_\varepsilon$.
Composition on a fixed smaller domain gives
$\|\Phi_{\rm tot}-\Psi_g\|_{C^2}\leqslant C|\varepsilon|$.
Since $\Psi_g(\theta,0,0)=(\theta,0,0)$, its restriction has
the stated torus displacement.
The total multiplier is $e^{-g\circ Q}\rho$, which need not be
close to one but is positive on the real torus. Applying
\eqref{eq:normaldampingtransform} to the final normal form gives
$\beta=\av{H_z\circ K_\varepsilon}$.
\end{proof}

\subsection{Fixing the constant parameter away from zero damping}

\begin{proposition}\label{prop:noe}
Fix $0<\lambda_0\leqslant\Lambda$.
For $\lambda_0\leqslant|\lambda|\leqslant\Lambda$, the conclusions of
Theorem~\ref{thm:main} hold with $e=0$, adjusting only $\xi$.
The domain, constants, and perturbation threshold may depend on
$\lambda_0$, and are uniform on this set of real damping values.
\end{proposition}

\begin{proof}
For $\beta\ne0$, define the full inverse of $D_\omega+\beta$ by
\[
    \mathcal R_\beta f(q)=
    \begin{cases}
        \displaystyle\int_0^\infty
            e^{-\beta t}f(q-\omega t)\,dt,&\beta>0,\\[2mm]
        \displaystyle-\int_0^\infty
            e^{\beta t}f(q+\omega t)\,dt,&\beta<0.
    \end{cases}
\]
Real translations preserve the complex angular strip, so
$\|\mathcal R_\beta f\|_s\leqslant|\beta|^{-1}\|f\|_s$.
Integration by parts, justified first on a smaller strip, gives
$(D_\omega+\beta)\mathcal R_\beta f=f$.
In particular, the constant Fourier coefficient is $\av f/\beta$.

Set $\zeta=0$, drop the normalization $\av u=0$, and define
\[
    u_0=\mathcal R_\beta a,\quad
    V_j=\mathcal R_\beta\ell_j,\quad u=u_0+V\eta.
\]
The angular average equation becomes
\[
    \mathsf M_\beta\eta=-\av{b-SDu_0-tu_0},\quad
    \mathsf M_\beta e_j=\av{Be_j-SDV_j-tV_j}.
\]
At the identity triple, $B=I$ and $\ell=0$, so
$\mathsf M_\beta=I$.
For $|\beta|\geqslant\lambda_0/2$, the full inverse estimate and
Cauchy's estimate across a fixed angular margin give
\[
    \|u_0\|_{C^1(\T^n)}\leqslant C(\lambda_0)E_s,\quad
    \|V\|_{C^1(\T^n)}\leqslant C(\lambda_0)Y_s.
\]
Thus $\|\mathsf M_\beta-I\|\leqslant C(\lambda_0)Y_s$, and
reducing the fixed transformation bound makes its inverse have
norm at most two. Solve for $\eta$, define $\delta\beta$ by
\eqref{eq:dbeta}, and obtain the zero-average $v,w$ from
\eqref{eq:vwexplicit}.

The generating-flow and remainder estimates do not use $\av u=0$.
They therefore yield the same quadratic recurrence, with constants
now depending on $\lambda_0$.
Choose the initial threshold so that the sum of the absolute
increments $|\delta\beta_j|$ is less than $\lambda_0/2$.
Starting with $\beta_0=\lambda$, finite-step induction then gives
$|\beta_j|\geqslant\lambda_0/2$ at every stage, while $e_j=0$
throughout. The preceding convergence proof applies.
The full inverse depends continuously on $\beta$ on each of the
two real parameter intervals away from zero, so the same argument
also proves parameter continuity there. No uniform estimate as
$\lambda_0\downarrow0$ follows, because the inverse of the
constant mode contains $1/\beta$.
\end{proof}

\subsection{First-order expansions and the role of the averages}

We return to the branch of Theorem~\ref{thm:main} constructed with
zero-average generating coefficients.
All derivatives in the coefficient formulas below are evaluated at
$(q,0,0)$. Put $P_0=P|_0$, $S_0=R_{0,pp}|_0$,
$t_0=R_{0,pz}|_0$, and $r_{20}=R_{0,zz}|_0$.
On a fixed smaller angular strip, the graph and parameters satisfy
\[
        \mathfrak u=\varepsilon u_1+O(\varepsilon^2),\quad \xi=\omega+\varepsilon\xi_1+O(\varepsilon^2),\quad e=\varepsilon e_1+O(\varepsilon^2),\quad \beta=\lambda+\varepsilon\beta_1+O(\varepsilon^2),
\]
where
\begin{equation}\label{eq:firstcoeff}
        u_1=-T_\lambda P_0,\quad e_1=\av{P_0},\quad \xi_1=-\av{P_p+S_0Du_1+t_0u_1},\quad \beta_1=\av{P_z+t_0\cdot Du_1+r_{20}u_1}.
\end{equation}
The remainders are uniform for real $\lambda\in[-\Lambda,\Lambda]$.
For the graph, they are bounded in the analytic supremum norm and,
on a further fixed smaller strip, in every fixed-order derivative
norm.

These expansions follow from the first step and the quadratic bound
on all subsequent corrections. At the identity accumulated map,
$B=I$, $\chi=1$, and $\ell=m=0$. The scalar generator and
constant increment are exactly
$u^{\rm gen}=\varepsilon T_\lambda P_0$ and
$\zeta=\varepsilon\av{P_0}$.
The normal Hessian at this step is that of $R_0+\varepsilon P$;
hence its blocks equal $S_0,t_0,r_{20}$ up to $O(\varepsilon)$.
The remaining correction formulas consequently have first-order
expansions with $O(\varepsilon^2)$ remainders on fixed smaller
strips.

The first time-one flow has vertical displacement
$\widehat U_1=-u^{\rm gen}+O(\varepsilon^2)$ and base map
$\widehat Q_1=\mathrm{id}+v^{\rm gen}+O(\varepsilon^2)$.
Moreover, Lemma~\ref{lem:step} with the fixed initial loss
$\sigma_0$ gives $E_1\leqslant C\varepsilon^2$.
Applying the weighted convergence estimate starting at stage one,
whose initial loss $\sigma_1$ is still fixed, shows that the sum
of all later parameter and triple increments is $O(\varepsilon^2)$.
Composition with the inverse base map therefore gives
$\mathfrak u=-u^{\rm gen}+O(\varepsilon^2)$.
This gives \eqref{eq:firstcoeff} and justifies
the expansions without assuming analytic dependence of the limiting
branch on $\varepsilon$.

For completeness, write
$Q(\theta)=\theta+\varepsilon q_1(\theta)+O(\varepsilon^2)$.
The first-order Hamilton--Jacobi and angular conjugacy equations are
\[
    (D_\omega+\lambda)u_1+P_0-e_1=0,\quad
    D_\omega q_1=\xi_1+P_p+S_0Du_1+t_0u_1.
\]
The chosen normalization gives $\av{u_1}=0$.
Averaging the second equation proves the formula for $\xi_1$.
Finally,
\[
    H_z\circ K
    =\lambda+\varepsilon
       (P_z+t_0\cdot Du_1+r_{20}u_1)+O(\varepsilon^2).
\]
The unperturbed value of $H_z$ is the constant $\lambda$, so the
first-order change of the angular parameterization contributes
nothing to this leading term. Averaging and using
\eqref{eq:meanidentity} gives the formula for $\beta_1$.

\begin{example}
Let $n=1$, $\lambda=0$, and
$H=\xi p+\tfrac12p^2+\varepsilon\cos q-e$, with fixed
$\xi>0$. A graph close to the zero section must satisfy
\[
    \mathfrak u'(q)
    =-\xi+\sqrt{\xi^2+2e-2\varepsilon\cos q}.
\]
If $e=0$ and $0<2|\varepsilon|<\xi^2$, strict concavity gives
\[
    \av{\sqrt{\xi^2-2\varepsilon\cos q}}
       <\sqrt{\av{\xi^2-2\varepsilon\cos q}}=\xi,
\]
contradicting $\av{\mathfrak u'}=0$.
For small $(e,\varepsilon)$, the derivative with respect to $e$
of the average equation at $(0,0)$ is $1/\xi$.
The implicit function theorem therefore determines a unique small
$e=e(\varepsilon)$.
The Hamilton--Jacobi equation and periodicity give the exact
identity $e=\tfrac12\av{(\mathfrak u')^2}$.
Since $\mathfrak u'=-\varepsilon\cos q/\xi+O(\varepsilon^2)$,
one obtains $e=\varepsilon^2/(4\xi^2)+O(\varepsilon^3)$.
Thus a zero first-order average does not remove the need for
higher-order parameter compensation.
\end{example}

\begin{example}
Let $V$ be bounded and holomorphic on $\T_s^n$, real on the real
torus, and consider
$H=\omega\cdot p+\lambda z+\varepsilon V(q)-e$.
Set $e=\varepsilon\av V$ and
\[
    \mathfrak u_\lambda=-\varepsilon T_\lambda V
      =-\varepsilon\sum_{k\in\Z^n\setminus\{0\}}
           \frac{V_k}{\lambda+ik\cdot\omega}e^{ik\cdot q}.
\]
Then $(D_\omega+\lambda)\mathfrak u_\lambda
       +\varepsilon V-e=0$.
The strict contact map
$\Psi_\lambda(q,p,z)
 =(q,p+D\mathfrak u_\lambda,z+\mathfrak u_\lambda)$ therefore
satisfies $H\circ\Psi_\lambda=\omega\cdot p+\lambda z$.
For $0<\sigma<\min\{s,1\}$, Lemma~\ref{lem:fourier} gives
\[
    \|\mathfrak u_\lambda-\mathfrak u_0\|_{s-\sigma}
    \leqslant C|\varepsilon\lambda|\gamma^{-2}
             \sigma^{-(2\tau+n)}\|V\|_s.
\]
\end{example}

\begin{example}\label{ex:mixedcoupling}
Let $n\geqslant2$, choose $k\in\Z^n\setminus\{0\}$, and put
$\nu=k\cdot\omega\ne0$ and $\phi=k\cdot q$.
Let $S$ be any fixed real symmetric matrix, $a\in\R^n$, and
$c\in\R$. Consider
\[
        H=\xi\cdot p+\lambda z+\tfrac12p^TSp +(a\cdot p)z\cos\phi+\tfrac c2 z^2 +\varepsilon\cos\phi-e.
\]
The Hessian blocks at the zero torus are $S_0=S$, $t_0=a\cos\phi$,
and $r_{20}=c$.
No invertibility of $S$ or of the full normal Hessian is required.
Solving $(D_\omega+\lambda)u_1=-\cos\phi$ with zero average gives
\[
    u_1=-\frac{\lambda\cos\phi+\nu\sin\phi}{\lambda^2+\nu^2},
    \quad
    Du_1=k\,\frac{\lambda\sin\phi-\nu\cos\phi}
                       {\lambda^2+\nu^2}.
\]
Since $P_p=P_z=0$ and $\av{Du_1}=\av{u_1}=0$, formula
\eqref{eq:firstcoeff} yields
\[
    e_1=0,\quad
    \xi_1=\frac{\lambda a}{2(\lambda^2+\nu^2)},\quad
    \beta_1=-\frac{\nu(a\cdot k)}{2(\lambda^2+\nu^2)}.
\]
These identities follow from $\av{\cos^2\phi}=1/2$ and
$\av{\sin\phi\cos\phi}=0$.
For input $\lambda=0$, the actual mean damping is therefore
$\beta=-\varepsilon(a\cdot k)/(2\nu)+O(\varepsilon^2)$.
If $a\cdot k\ne0$, a mean damping of first order is created by
the displacement through the mixed background, even though the
perturbing function has no $z$-dependence. This example also has
no closed $(q,p)$ dynamics in general, since
$H_p=\xi+Sp+az\cos\phi$ depends on $z$.
\end{example}

\subsection{Persistence of the graph and of its frequency}
\label{sec:comparison}

The persistence of an invariant graph does not by itself preserve a
prescribed frequency. To see this, let $n\geqslant2$, fix $\lambda>0$,
and choose $\omega$ satisfying \eqref{eq:dio}. Consider
\[
    H_\varepsilon(q,p,z)
    =\omega\cdot p+\frac12\sum_{i=1}^n p_i^2
       +\lambda z+\varepsilon p_1.
\]
The contact equations are
\[
    \dot q=\omega+p+\varepsilon e_1,\quad
    \dot p=-\lambda p,\quad
    \dot z=\frac12\sum_{i=1}^n p_i^2-\lambda z.
\]
The zero graph remains invariant for every $\varepsilon$, but its
frequency is $\Omega_\varepsilon=\omega+\varepsilon e_1$.
Moreover, every trajectory contained in a compact invariant set for
all positive and negative times must have $p=0$, since
$p(t)=e^{-\lambda t}p(0)$. Its vertical coordinate then satisfies
$z(t)=e^{-\lambda t}z(0)$, so it must also have $z=0$.
Thus every compact fully invariant set lies in the zero graph.

Nonresonance of $\omega$ implies $\omega_2\ne0$ and
$\omega_1/\omega_2\notin\mathbb Q$.
Choose coprime integers $a_j,b_j$, with $b_j>0$, such that
$a_j/b_j\to\omega_1/\omega_2$, and set
\[
        \varepsilon_j=\omega_2\frac{a_j}{b_j}-\omega_1 \to0,\quad k_j=(b_j,-a_j,0,\ldots,0),\quad k_j\cdot\Omega_{\varepsilon_j}=0.
\]
At these perturbation values, every trajectory on the graph stays
in a proper subtorus $k_j\cdot q=\mathrm{constant}\pmod{2\pi}$.
The graph flow therefore cannot be conjugate to the nonresonant
linear flow with frequency $\omega$.
There is no other compact invariant embedded $n$-torus on which to
recover that flow: any such torus would be contained in the zero
graph and, being a compact embedded submanifold of the same
dimension, would coincide with it.

Replacing $\omega\cdot p$ by $\xi\cdot p$ and taking
$\xi=\omega-\varepsilon e_1$ restores the frequency exactly.
For comparison, the symplectic Hamiltonian
$h_\varepsilon(p)=\omega\cdot p+\tfrac12\sum_{i=1}^n p_i^2
+\varepsilon p_1$ already has an invariant torus of frequency
$\omega$ at $p=-\varepsilon e_1$.
That momentum level is not invariant for the damped contact system
when $\varepsilon\ne0$, because $\dot p=-\lambda p$.
The external frequency parameter thus supplies a freedom that
cannot be obtained by selecting a nearby constant momentum level
in this example.

\section{Nonzero mean damping and analytic contact linearization}
\label{sec:linearization}

We now study the finite-dimensional local dynamics of
the normal form obtained in Theorem~\ref{thm:main}. In this section the angular frequency may be
any real vector. The arithmetic assumption used to construct the
torus is not needed for its local linearization when the mean damping
is nonzero. Related rotation--contraction normal forms for
conformally symplectic systems are established in
\cite{CallejaLocal2013}. Here we construct the conjugacy in the full
contact space and verify its contact property directly.

\begin{theorem}
\label{thm:linearization}
Let $\omega\in\R^n$ and $\beta\in\R\setminus\{0\}$. Suppose
that $\widehat H(\theta,J,W)=\omega\cdot J+\beta W+R(\theta,J,W)$
is real analytic on a neighborhood of the zero torus
$\mathcal T_0=\{J=W=0\}$, with $R=O_2$ in $(J,W)$.
There are a neighborhood $\mathcal U_0$ of $\mathcal T_0$ and a
real-analytic contact diffeomorphism $F$ from $\mathcal U_0$ onto
its image, fixing $\mathcal T_0$ pointwise, which conjugates the
flow of $\widehat H$ to
\begin{equation}\label{eq:linearflow}
    \dot\Theta=\omega,\quad \dot Y=-\beta Y,
    \quad Y=(Y_1,Y_2)\in\R^n\times\R.
\end{equation}
The conjugacy preserves time. With
$\alpha=dW-J\cdot d\theta$, $\alpha_*=dY_2-Y_1\cdot d\Theta$,
and $H_*(\Theta,Y)=\omega\cdot Y_1+\beta Y_2$, one has
$F^*\alpha_*=f\alpha$ and $H_*\circ F=f\widehat H$ for a
positive real-analytic function $f$.

The zero torus is locally attracting for $\beta>0$ and locally
repelling for $\beta<0$. Its $n+1$ normal Lyapunov exponents are
$-\beta$, and its tangent exponents are zero. No nonresonance
condition on $\omega$ is required. The neighborhood and estimates
may depend on $|\beta|^{-1}$.
\end{theorem}

\begin{proof}
Assume first that $\beta>0$ and write $N=(J,W)$. The contact
equations give
\begin{equation}\label{eq:normalODE}
        \dot\theta=\omega+R_J(\theta,N),\quad \dot N=-\beta N+A(\theta,N),\quad A=\binom{-R_\theta-JR_W}{J\cdot R_J-R}.
\end{equation}
Choose $a,b>0$ such that the coefficients extend holomorphically
to a neighborhood of $\overline{D_{a,b}}$. In this proof, absolute
values of vectors denote the maximum norm. Since $R=O_2$, there
are constants $C_1,C_2>0$ such that
\begin{equation}\label{eq:normalbounds}
    |A(\theta,N)|\leqslant C_2|N|^2,\quad
    |R_J(\theta,N)|+|R_W(\theta,N)|\leqslant C_1|N|
\end{equation}
on $D_{a,b}$. Angular differentiation preserves normal order, so
$R_\theta=O_2$. Thus the normal equation has linear part
$-\beta I_{n+1}$, although the angular equation can still contain
terms linear in $N$.

For complex initial data $x=(\theta,N)$, variation of constants
gives, as long as the trajectory remains in $D_{a,b}$,
\[
    e^{\beta t}|N(t)|
    \leqslant |N|+C_2\int_0^t e^{-\beta s}
                       \bigl(e^{\beta s}|N(s)|\bigr)^2\,ds.
\]
Denote the right-hand side by $B_x(t)$. If $N\ne0$, then
$B_x(0)=|N|$ and $B_x'(t)\leqslant C_2e^{-\beta t}B_x(t)^2$.
Integrating the inequality for $1/B_x$ gives, whenever
$C_2|N|\leqslant\beta/2$,
\begin{equation}\label{eq:Ndecay}
    |N(t)|\leqslant
    \frac{|N|e^{-\beta t}}
         {1-(C_2/\beta)|N|(1-e^{-\beta t})}
    \leqslant2|N|e^{-\beta t}.
\end{equation}
The case $N=0$ is the invariant zero torus. On the universal cover,
the angular equation consequently gives
\begin{equation}\label{eq:thetadrift}
    |\theta(t)-\theta-\omega t|
    \leqslant C_1\int_0^t|N(s)|\,ds
    \leqslant 2C_1\beta^{-1}|N|.
\end{equation}

Fix $0<h<a$ and choose $\delta_\beta>0$ so that
$2\delta_\beta<b$, $2C_2\delta_\beta<\beta$, and
$2C_1\delta_\beta/\beta<h$.
For initial data in $D_{a-h,\delta_\beta}$,
\eqref{eq:Ndecay}--\eqref{eq:thetadrift} exclude a first exit
from $D_{a,b}$: the normal variables remain strictly within radius
$b$, and the imaginary angular displacement is less than $h$.
Here the translation $\omega t$ is real and does not change the
imaginary parts. These bounds also keep each trajectory in a compact
subset of the coefficient domain modulo the angular lattice.
The local existence theorem then extends the trajectory for every
$t\geqslant0$. We have therefore obtained one common complex initial
domain on which all positive-time trajectories are defined.

Let $\varphi_t$ denote this flow and set
$L_t(\Theta,Y)=(\Theta+\omega t,e^{-\beta t}Y)$.
Define
\[
        \Theta(x)=\theta+\int_0^\infty R_J(\varphi_sx)\,ds,\quad Y(x)=N+\int_0^\infty e^{\beta s}A(\varphi_sx)\,ds, \quad F=(\Theta,Y).
\]
The decay estimate implies
\[
        \left|\int_t^\infty R_J(\varphi_sx)\,ds\right| \leqslant 2C_1\beta^{-1}|N|e^{-\beta t},\quad \left|\int_t^\infty e^{\beta s}A(\varphi_sx)\,ds\right| \leqslant 4C_2\beta^{-1}|N|^2e^{-\beta t}.
\]
For finite upper limits the integrals depend holomorphically on
$x$, and these tail bounds give locally uniform convergence.
Thus $F$ is holomorphic and real on real points. The corrections
to $\theta$ and $N$ are periodic in the initial angle, so the map
descends from the universal cover to a neighborhood of the torus.
Variation of constants also gives
\[
    L_{-t}\circ\varphi_t(x)
    =\left(\theta+\int_0^t R_J(\varphi_sx)\,ds,
            N+\int_0^t e^{\beta s}A(\varphi_sx)\,ds\right).
\]
Consequently $F=\lim_{t\to\infty}L_{-t}\circ\varphi_t$.
Cauchy's estimate passes this convergence to every fixed derivative
order on smaller complex domains.

Along the zero torus,
\begin{equation}\label{eq:DF}
    F(\theta,0)=(\theta,0),\quad
    DF(\theta,0)=
    \begin{pmatrix}
        I_n&\mathsf C_\beta(\theta)\\
        0&I_{n+1}
    \end{pmatrix},
\end{equation}
where
\[
    \mathsf C_\beta(\theta)
    =\int_0^\infty
        \bigl(R_{JJ}(\theta+\omega s,0),
              R_{JW}(\theta+\omega s,0)\bigr)e^{-\beta s}\,ds.
\]
Indeed, $Y=N+O(|N|^2/\beta)$ and
$\Theta-\theta=O(|N|/\beta)$. In differentiating the angular
integral at $N=0$, the normal variation is $e^{-\beta s}\delta N$,
while $D_\theta R_J(\theta,0)=0$.
The matrix in \eqref{eq:DF} is invertible even when
$\mathsf C_\beta$ is not small.

The inverse function theorem and compactness give local inverses
around the entire real torus. They yield injectivity on a sufficiently
small common neighborhood. To see this, otherwise there would be
distinct $x_j,y_j$ approaching the torus with $F(x_j)=F(y_j)$.
After taking subsequences, compactness and $F|_{\mathcal T_0}
=\mathrm{id}$ force both limits to be the same point of
$\mathcal T_0$. Both sequences eventually lie in one local inverse
neighborhood, a contradiction. This argument also gives a sufficiently
small complex neighborhood on which $F$ is biholomorphic onto its
image.

For fixed $t$ and points for which $x$ and $\varphi_t x$ lie in
the domain of $F$, the flow property gives
\[
        F(\varphi_t x) =\lim_{s\to\infty}L_{-s}\varphi_{s+t}x =L_t\lim_{s\to\infty}L_{-(s+t)}\varphi_{s+t}x =L_tF(x).
\]
In particular, differentiation at $t=0$ proves
$DF\,X_{\widehat H}=X_{H_*}\circ F$.

It remains to check the contact identity. By \eqref{eq:contact},
\[
    \varphi_t^*\alpha
    =\exp\left(-\beta t-\int_0^tR_W(\varphi_sx)\,ds\right)
       \alpha.
\]
The linear map $L_{-t}$ multiplies the standard contact form by
$e^{\beta t}$, and hence
\[
    (L_{-t}\circ\varphi_t)^*\alpha_*
    =\exp\left(-\int_0^tR_W(\varphi_sx)\,ds\right)\alpha.
\]
The integral converges locally uniformly by
\eqref{eq:normalbounds} and \eqref{eq:Ndecay}. Passing to the
limit using the convergence of first derivatives gives
\[
    F^*\alpha_*=f\alpha,\quad
    f(x)=\exp\left(-\int_0^\infty R_W(\varphi_sx)\,ds\right).
\]
The factor $f$ is holomorphic and nonzero on a smaller complex
domain, and is positive on its real part. Evaluating this identity
on $X_{\widehat H}$ and using the vector-field conjugacy gives
$H_*\circ F=f\widehat H$.

Choose a real product neighborhood $\T^n\times B_\delta$ with
closure in the image of $F$. Its inverse image is forward invariant,
because $L_t$ contracts the normal ball. This proves local
attraction. If $\beta<0$, apply the preceding construction to
$-\widehat H$, whose frequency is $-\omega$, damping is $-\beta>0$,
and remainder is $-R$. Reversing time gives the claimed contact
conjugacy for $\widehat H$ and a backward-invariant neighborhood
on which the zero torus is attracting in reverse time.

Finally, the derivative of the linear flow induces
$e^{-\beta t}I_{n+1}$ on the normal quotient by
$T\mathcal T_0$. The induced derivatives of $F$ and $F^{-1}$
are bounded on the compact torus. The normal Lyapunov exponents
are therefore $-\beta$. The derivative along the torus is the
identity, so the tangent exponents are zero.
\end{proof}

The map $F$ in this theorem is a general local contact
diffeomorphism. It need not belong to the affine class
\eqref{eq:affine} used in the persistence proof; both its angular
component and its multiplier may depend on the normal variables.

\subsection{Asymptotic phase and stable leaves}

Fix a smooth product metric on the real contact space and denote
its distance by $d$.

\begin{corollary}\label{cor:phase}
Suppose that the torus of Theorem~\ref{thm:main} has $\beta>0$.
There are a relatively compact, forward-invariant open neighborhood
$\mathcal U$ of its image and a real-analytic map
$\vartheta:\mathcal U\to\T^n$, with
$\vartheta\circ K=\mathrm{id}$, such that
\begin{equation}\label{eq:tracking}
    d\bigl(\varphi_t x,K(\vartheta(x)+\omega t)\bigr)
    \leqslant C_xe^{-\beta t},\quad t\geqslant0.
\end{equation}
For the fixed parameterization $K$, this phase is unique. Its
fibers are local analytic strong-stable leaves of dimension $n+1$.

The probability measure
\begin{equation}\label{eq:mu}
    \mu=K_*\left(\frac{d\theta}{(2\pi)^n}\right)
\end{equation}
is the unique invariant probability measure on the torus. For every
$x\in\mathcal U$ and $a\in C(\overline{\mathcal U})$,
\begin{equation}\label{eq:timeaverage}
    \lim_{T\to\infty}\frac1T\int_0^T a(\varphi_t x)\,dt
    =\int a\,d\mu.
\end{equation}
In particular, $\mu$ has an open basin of time averages and is a
physical measure with respect to any smooth positive volume density
on a neighborhood of $\overline{\mathcal U}$.
\end{corollary}

\begin{proof}
Let $\Phi$ be the normalizing map of Theorem~\ref{thm:main} and
let $F$ be the map of Theorem~\ref{thm:linearization} for the
resulting normal form. The map
$\mathcal T=\Phi\circ F^{-1}$ takes linear coordinates to the
original coordinates and satisfies $\mathcal T(\Theta,0)=K(\Theta)$.
Choose $\delta>0$ so that $\mathcal T$ is analytic on a neighborhood
of $\T^n\times\overline{B_\delta}$, and put
$\mathcal U=\mathcal T(\T^n\times B_\delta)$.
Define $\vartheta=\operatorname{pr}_\Theta\circ\mathcal T^{-1}$
and denote the normal component of $\mathcal T^{-1}(x)$ by $Y(x)$.
Then
\[
        \vartheta(\varphi_t x) =\vartheta(x)+\omega t\pmod{2\pi},\quad \varphi_t\mathcal T(\Theta,Y) =\mathcal T(\Theta+\omega t,e^{-\beta t}Y).
\]
These identities show forward invariance of $\mathcal U$ and
existence of its trajectories for all $t\geqslant0$.
The derivatives of $\mathcal T$ are bounded on the closed product,
so
\[
    d\bigl(\varphi_t x,K(\vartheta(x)+\omega t)\bigr)
    \leqslant C|Y(x)|e^{-\beta t}.
\]
This proves \eqref{eq:tracking}.
The phase map is a submersion, and each fiber is the image of
$\{\Theta\}\times B_\delta$. Its dimension is $n+1$, and for $t\geqslant0$ the flow carries it
into the fiber with phase $\Theta+\omega t$.
Points in the same fiber converge exponentially to the same torus
trajectory.

If two phases track the same trajectory, the distance between their
corresponding torus trajectories tends to zero. The inverse of $K$
is uniformly continuous on its compact image, so the flat distance
between these two phases after translation by $\omega t$ tends to
zero. That distance is independent of $t$, proving uniqueness.
This part of the argument does not require nonresonance.

For the measure assertion, the frequency in Theorem~\ref{thm:main}
satisfies $k\cdot\omega\ne0$ for every nonzero integer vector $k$.
For each such $k$,
\begin{equation}\label{eq:Fourieraverage}
    \frac1T\int_0^T e^{ik\cdot(\theta+\omega t)}\,dt
    =e^{ik\cdot\theta}
       \frac{e^{iTk\cdot\omega}-1}{iT k\cdot\omega}
    \to0
\end{equation}
uniformly in $\theta$. Approximation by trigonometric polynomials
gives uniform convergence of time averages for every continuous
function on $\T^n$ to its Haar integral. Integrating these averages
against any invariant probability proves that Haar probability is
the only invariant probability of the linear flow. Pushing forward
by $K$ proves the assertion about \eqref{eq:mu}.

For $a\in C(\overline{\mathcal U})$, uniform continuity and
\eqref{eq:tracking} give
$a(\varphi_t x)-a(K(\vartheta(x)+\omega t))\to0$.
The difference is bounded, and its Ces\`aro mean tends to zero.
Together with \eqref{eq:Fourieraverage}, this proves
\eqref{eq:timeaverage}. Thus the basin contains $\mathcal U$,
which has positive measure for every smooth positive volume density.
Only nonresonance, rather than the quantitative Diophantine bound,
is used for these time-average conclusions.
\end{proof}

\begin{remark}\label{rem:basedensity}
Let $K$ be a torus parameterization from Theorem~\ref{thm:main},
and let $\mu=K_*(d\theta/(2\pi)^n)$.
This probability is invariant under the graph flow, independently
of the sign of the mean damping. Write
$K(\theta)=(Q(\theta),D\mathfrak u(Q(\theta)),
\mathfrak u(Q(\theta)))$ and let $\pi(q,p,z)=q$.
Changing variables under $Q$ gives
\[
    d(\pi_*\mu)(q)=\varpi(q)\,dq,\quad
    \varpi(q)=\frac{1}
        {(2\pi)^n|\det DQ(Q^{-1}(q))|}.
\]
The density is positive and real analytic. The absolute value may
be omitted for the near-identity map $Q$ of the main theorem.

For the induced base vector field
$X(q)=H_p(q,D\mathfrak u(q),\mathfrak u(q))$ and any smooth
periodic test function $\phi$,
\[
    \int_{\T^n}D\phi(q)\cdot X(q)\,\varpi(q)\,dq
    =(2\pi)^{-n}\int_{\T^n}D_\omega(\phi\circ Q)(\theta)\,d\theta
    =0.
\]
Integration by parts gives $\operatorname{div}(\varpi X)=0$.
This is an invariant density for the base flow, not preservation of
the ambient contact volume; compare \eqref{eq:introenergyvolume}.
\end{remark}

\subsection{Dependence of the linearization neighborhood on damping}

\begin{example}\label{ex:nonuniform}
Let $H=\omega\cdot p+\beta z+az^2$, where $a\ne0$ is fixed and
$\beta\ne0$. Near the zero torus, define
\[
    F_\beta(q,p,z)
    =\left(q,\frac{p}{(1+az/\beta)^2},
                 \frac{z}{1+az/\beta}\right).
\]
Writing $b_\beta(z)=1+az/\beta$, one has
$d(z/b_\beta)=b_\beta^{-2}dz$ and hence
$F_\beta^*(dY_2-Y_1\cdot dq)=b_\beta^{-2}\alpha$.
Moreover,
\[
    (\omega\cdot Y_1+\beta Y_2)\circ F_\beta
    =\frac{\omega\cdot p+\beta z+az^2}{b_\beta^2}.
\]
Thus $F_\beta$ is a time-preserving contact linearization wherever
it is defined near the zero torus.

The pole at $z=-\beta/a$ reflects an obstruction to a common
linearization neighborhood. The original flow has two disjoint
compact invariant tori,
$\mathcal T_0=\{p=z=0\}$ and
$\mathcal T_1=\{p=0,\ z=-\beta/a\}$, both carrying
$\dot q=\omega$.
Every compact set invariant for all times under the linear flow
\eqref{eq:linearflow} lies in $Y=0$, since a nonzero normal
coordinate is unbounded in one time direction.
A conjugacy fixing $\mathcal T_0$ pointwise already maps it onto
the entire target zero torus. If its domain contained
$\mathcal T_1$, the image of that torus would also have to lie in
$Y=0$, contradicting injectivity. As $\beta\to0$,
$\mathcal T_1$ approaches $\mathcal T_0$. Therefore a
linearization neighborhood of fixed positive normal size cannot in
general be chosen uniformly as the mean damping vanishes.
\end{example}

\section{Zero mean damping: variational shear and quadratic models}
\label{sec:zero}

In this section we consider the finite-dimensional
normal form with actual mean damping $\beta=0$.
This is not implied by setting the input parameter $\lambda=0$;
the distinction is given by \eqref{eq:meanidentity}.
The exponential decay used in Section~\ref{sec:linearization} is
now absent. We first describe the linearized flow and then consider
two quadratic models.

\subsection{The full variational equation}

Let $\widehat H=\omega\cdot J+R(\theta,J,W)$, where $R=O_2$
and $\omega$ satisfies \eqref{eq:dio}.
Writing $N=(J,W)$, the variational equation along
$(\theta_0+\omega t,0)$ is
\begin{equation}\label{eq:variation}
        \delta\dot N=0,\quad \delta\dot\theta=\mathcal B(\theta_0+\omega t)\delta N,\quad \mathcal B(\theta) =\bigl(R_{JJ}(\theta,0),R_{JW}(\theta,0)\bigr).
\end{equation}
Indeed, the normal remainder $A$ in \eqref{eq:normalODE} is
$O_2$, whereas $D_NR_J|_0=\mathcal B$ and
$D_\theta R_J|_0=0$.
The matrix $\mathcal B$ has size $n\times(n+1)$.
Apply Lemma~\ref{lem:fourier} componentwise to obtain the
zero-average analytic matrix $\mathcal V=T_0\mathcal B$,
which satisfies $D_\omega\mathcal V=\mathcal B-\av{\mathcal B}$.
Integration gives
\begin{equation}\label{eq:shear}
    \delta\theta(t)=\delta\theta(0)
       +\bigl[t\av{\mathcal B}
              +\mathcal V(\theta_0+\omega t)
              -\mathcal V(\theta_0)\bigr]\delta N(0).
\end{equation}
Consequently,
\[
    D\varphi_t(\theta_0,0)
    =\begin{pmatrix}
        I_n&t\av{\mathcal B}
              +\mathcal V(\theta_0+\omega t)-\mathcal V(\theta_0)\\
        0&I_{n+1}
      \end{pmatrix}.
\]
The inverse matrix changes the sign of the upper-right block.
Both norms are bounded by $C(1+|t|)$, uniformly in $\theta_0$.
In particular, for every tangent vector $v$,
\[
    \frac{|v|}{C(1+|t|)}
    \leqslant|D\varphi_t(\theta_0,0)v|
    \leqslant C(1+|t|)|v|.
\]
All Lyapunov exponents on the torus are therefore zero.
The linear growth bound itself also follows directly by integrating
the bounded matrix $\mathcal B$; the Diophantine assumption gives
the sharper decomposition into constant shear and a bounded
oscillatory term in \eqref{eq:shear}.

If $\av{\mathcal B}\ne0$, choose $\nu$ with
$\av{\mathcal B}\nu\ne0$. The variation with initial value
$(\delta\theta(0),\delta N(0))=(0,\nu)$ grows linearly in at
least one angular direction. This excludes a time-preserving
$C^1$ conjugacy on a neighborhood of the whole zero torus to
$(\dot\Theta,\dot Y)=(\omega,0)$.
Indeed, the derivative of that model flow is the identity. Such a
conjugacy $\mathcal C$ would give
\[
    D\varphi_t(x)
    =\bigl(D\mathcal C(\varphi_t x)\bigr)^{-1}D\mathcal C(x),
    \quad x\in\mathcal T_0,
\]
which is uniformly bounded in $t$ by compactness of the torus.
If $\av{\mathcal B}=0$, the derivative flow is bounded, but this
alone does not remove the nonlinear terms.

\subsection{Polynomial normal motion}

Consider $H=\omega\cdot J+aW^2$ with a nonzero real constant $a$.
The equations are $\dot\theta=\omega$, $\dot W=-aW^2$, and
$\dot J=-2aWJ$. Their solution is
\begin{equation}\label{eq:polynomial}
    W(t)=\frac{W_0}{1+aW_0t},\quad
    J(t)=\frac{J_0}{(1+aW_0t)^2},\quad
    \theta(t)=\theta_0+\omega t,
\end{equation}
on the connected time interval containing $0$ on which the
denominator does not vanish.
If $aW_0>0$, the solution exists for all positive times, with
$W(t)=O(t^{-1})$ and $J(t)=O(t^{-2})$.
If $aW_0<0$, the normal coordinate $W(t)$ becomes unbounded at
$t_*=-1/(aW_0)>0$. In a local neighborhood, the trajectory leaves
the neighborhood before this time.
On the boundary $W_0=0$, one instead has $W(t)=0$ and
$J(t)=J_0$ for all time.

For this model $\mathcal B=0$, so even the shear in
\eqref{eq:variation} vanishes. Nevertheless, arbitrarily small
initial data on the side $aW_0<0$ leave a fixed small neighborhood
of the zero torus. The torus is therefore not Lyapunov stable and
cannot be locally conjugate, with the zero torus preserved, to the
neutral model $(\dot\Theta,\dot Y)=(\omega,0)$.
Thus neither zero Lyapunov exponents nor vanishing mean shear
implies nonlinear stability or linearization.

\subsection{Exact normalization of a variable quadratic coefficient}

\begin{proposition}\label{prop:quadraticmodel}
Let $H(q,p,z)=\omega\cdot p+a(q)z^2$, where $a$ is real analytic
on $\T^n$ and $\omega$ satisfies \eqref{eq:dio}.
Put $\bar a=\av a$ and let $\psi=T_0a$, so that
$D_\omega\psi=a-\bar a$ and $\av\psi=0$.
On a sufficiently small neighborhood of the zero torus, the map
\begin{equation}\label{eq:quadratictransform}
    \mathcal G(q,p,z)
    =\left(q,
           \frac{p+z^2D\psi(q)}{(1-\psi(q)z)^2},
           \frac{z}{1-\psi(q)z}\right)
\end{equation}
is an analytic contact diffeomorphism onto its image, fixing the
zero torus pointwise. It conjugates $X_H$, in the original time,
to the contact Hamiltonian flow of
$H_{\rm av}(q,J,W)=\omega\cdot J+\bar aW^2$.
\end{proposition}

\begin{proof}
Lemma~\ref{lem:fourier} gives $\psi$ on any fixed smaller angular
strip within a bounded holomorphic extension domain of $a$.
Writing $(q,J,W)=\mathcal G(q,p,z)$, direct differentiation gives
\[
    dW=\frac{dz+z^2D\psi\cdot dq}{(1-\psi z)^2},\quad
    dW-J\cdot dq=\frac{dz-p\cdot dq}{(1-\psi z)^2}.
\]
Moreover,
\[
    H_{\rm av}\circ\mathcal G
    =\frac{\omega\cdot p+(D_\omega\psi+\bar a)z^2}
            {(1-\psi z)^2}
    =\frac{H}{(1-\psi z)^2}.
\]
The inverse leaves $q$ unchanged and is given by
\[
    z=\frac{W}{1+\psi(q)W},\quad
    p=\frac{J-W^2D\psi(q)}{(1+\psi(q)W)^2}.
\]
On a fixed smaller strip, $\psi$ and $D\psi$ are bounded.
Choose the normal domains so that $|\psi z|<1/2$ for
\eqref{eq:quadratictransform} and $|\psi W|<1/2$ for its
inverse. These formulas then define inverse holomorphic maps on
suitable neighborhoods and are real on real points.
The contact factor is $(1-\psi z)^{-2}>0$ on the real domain.
The two pullback identities show that the pullback of the target
Hamiltonian, divided by this contact factor, is exactly $H$;
hence the conjugacy preserves time.
\end{proof}

If $\bar a=0$, this particular family is exactly conjugate to
$H_{\rm av}=\omega\cdot J$ and has neutral normal dynamics.
If $\bar a\ne0$, its dynamics in the normalized coordinates are
given by \eqref{eq:polynomial} with $a=\bar a$.
These conclusions use the special form $a(q)z^2$ and do not
classify a general normal remainder $R=O_2$.

\section{Proof of the infinite-dimensional theorem}
\label{sec:infinite-proof}

We prove Theorem~\ref{thm:infinite} using the correction equations
of Sections~\ref{sec:linear}--\ref{sec:iteration}. We first establish
the spatial estimates for the contact bracket and the Fourier
inverse, then solve the average equations on $V$ and choose the
analytic losses for the iteration.

\subsection{Contact brackets and the action space}
\label{sec:infinite-brackets}

We first estimate contact Hamiltonian vector fields and their
brackets in the spatial norms. Derivatives are taken in the
Fr\'echet sense on the Banach spaces of
Section~\ref{sec:infinite-setting}.

For a spatial block
\[
J_A(q_A,p_A,z)
=
\sum_{k\in\Z^A}
J_{A,k}(p_A,z)e^{ik\cdot q},
\]
the Hamiltonian vector field of an individual Fourier
coefficient is
\[
\mathcal X_{A,k}(p,z)
=
\left(
D_{p_A}J_{A,k},
-ikJ_{A,k}-p\,\partial_zJ_{A,k},
\langle p,D_{p_A}J_{A,k}\rangle-J_{A,k}
\right).
\]
The second component belongs to $P_w$. The common contact
variable produces the term $-p\,\partial_zJ_{A,k}$, whose
support may extend beyond $A$.

We measure a Hamiltonian and its vector field together by
\begin{equation}\label{eq:infinite-contactnorm}
 |J|^c_{m,s,r}
 =
 \sum_{A\in\mathcal S}
 \sum_{k\in\Z^A}
 e^{m[A]+s|k|_1}
 \left(
 \sup_{|p|_w<r,\,|z|<r}|J_{A,k}|
 +
 \sup_{|p|_w<r,\,|z|<r}
 \|\mathcal X_{A,k}\|_E
 \right).
\end{equation}
The displayed sum is the norm of a representative. The blockwise
estimates below pass to the quotient by taking infima. Normal jets
and restrictions commute with summation. For fixed $(p,z)$ the angular
series converges uniformly on the real product torus, where integration
against its characters determines the total Fourier coefficients.
Hence averaging and Fourier inversion depend only on the represented
function. In particular, $R_0=O_2$ gives $\mathcal J_1R_0=0$ in the
quotient.

The following estimate controls the contact norm by the scalar
spatial norm.

\begin{lemma}
\label{lem:infinite-vectorfield}
Let $m\geqslant0$, $0<\sigma<s$, and $0<d<r$. Then
\begin{equation}\label{eq:infinite-initialnorm}
 |J|^c_{m,s-\sigma,r-d}
 \leqslant
 C(1+\sigma^{-1}+d^{-1})
 |J|_{m+w,s,r}.
\end{equation}
The constant is independent of the cardinalities of the
finite supports.
\end{lemma}

\begin{proof}
It is enough to estimate a single Fourier coefficient
$J_{A,k}(p_A,z)e^{ik\cdot q}$. We write
\[
M_{A,k}
=
\sup_{|p|_w<r,\,|z|<r}
|J_{A,k}(p_A,z)|.
\]
Since $J_{A,k}$ depends only on $p_A$, the Banach Cauchy
estimate in the direction of the $i$-th coordinate gives,
on the smaller action ball,
\[
\left|
\frac{\partial J_{A,k}}{\partial p_i}
\right|
\leqslant
C d^{-1}e^{w[i]}M_{A,k},
\quad i\in A.
\]
By the definition of $[i]$ and the monotonicity of the
spatial weight,
\[
[i]\leqslant[A].
\]
Consequently,
\[
\|D_pJ_{A,k}\|_\infty
\leqslant
C d^{-1}e^{w[A]}M_{A,k}.
\]
Similarly,
\[
|\partial_zJ_{A,k}|
\leqslant Cd^{-1}M_{A,k}.
\]

The angular derivative satisfies
\[
    \|D_q(J_{A,k}e^{ik\cdot q})\|_w \leqslant e^{w[A]}|k|_1M_{A,k}e^{s|k|_1}.
\]
The common contact term is controlled by
\[
|p\,\partial_zJ_{A,k}|_w
\leqslant
Cr d^{-1}M_{A,k}.
\]
Finally, the duality estimate between $P_w$ and $V$ yields
\[
|\langle p,D_pJ_{A,k}\rangle|
\leqslant
|p|_w\|D_pJ_{A,k}\|_\infty
\leqslant
Cr d^{-1}e^{w[A]}M_{A,k}.
\]
Combining these estimates gives
\[
    \sup\|\mathcal X_{A,k}\|_E \leqslant C e^{w[A]} (1+d^{-1}+|k|_1)M_{A,k},
\]
where the constant may depend on the fixed normal radius
$r$, but not on $A$ or its cardinality.

We now use
\[
e^{-\sigma|k|_1}|k|_1
\leqslant (e\sigma)^{-1}
\]
and sum over all spatial blocks and Fourier indices.
This proves \eqref{eq:infinite-initialnorm}.
\end{proof}

We next establish the support property of the contact
bracket. For two Hamiltonians $F$ and $J$, define
\[
\mathscr L_FJ=X_F(J)+F_zJ.
\]
In coordinates this is
\begin{equation}\label{eq:infinite-contact-bracket-coordinate}
    \mathscr L_FJ =F_p\cdot J_q-F_q\cdot J_p -F_z\langle p,J_p\rangle +J_z\langle p,F_p\rangle -FJ_z+F_zJ.
\end{equation}

\begin{lemma}
\label{lem:infinite-support}
Let $F_A$ and $J_B$ be Hamiltonian blocks with spatial
supports contained in $A$ and $B$, respectively. Then
\[
\operatorname{supp}(\mathscr L_{F_A}J_B)
\subseteq A\cup B.
\]
Moreover, the spatial weight satisfies
\[
e^{m[A\cup B]}
\leqslant e^{m[A]}e^{m[B]},
\quad m\geqslant0.
\]
\end{lemma}

\begin{proof}
The terms involving $F_p\cdot J_q$ and
$F_q\cdot J_p$ involve only indices in $A\cap B$.
The terms
$\langle p,F_p\rangle$ and $\langle p,J_p\rangle$
involve only the corresponding spatial supports, namely
$A$ and $B$. The remaining terms are products of functions
depending on the variables in $A$ and $B$, together with
the common scalar variable $z$.

Thus every scalar term in
\eqref{eq:infinite-contact-bracket-coordinate} depends only
on the variables indexed by $A\cup B$, and the asserted
support property follows.

The weight estimate follows directly from
\eqref{eq:infinite-weights} and $[A\cap B]\geqslant0$:
\[
[A\cup B]\leqslant[A]+[B].
\]
\end{proof}

For example, if $F=f(q_i)$ and $J=zg(q_j)$ with $i\ne j$, then
\[
\mathscr L_FJ=-f(q_i)g(q_j),
\]
which has spatial support $\{i,j\}$.

We now estimate the full contact bracket, including its
Banach vector-field component.

\begin{lemma}
\label{lem:infinite-bracket}
Let $m>0$, $0<\sigma<s$, and $0<d<r$. Assume that the
Hamiltonians $F$ and $J$ have finite contact action norms
on $\mathcal D_{s,r}$. Then
\begin{equation*}
 |\mathscr L_FJ|^c_{m,s-\sigma,r-d}
 \leqslant
 C(\sigma^{-1}+d^{-1})
 |F|^c_{m,s,r}|J|^c_{m,s,r}.
\end{equation*}
The constant is independent of the cardinalities of the
spatial supports.
\end{lemma}

\begin{proof}
It suffices to estimate a pair of spatial Fourier terms
$f(p_A,z)e^{ik\cdot q}$ and $j(p_B,z)e^{il\cdot q}$.
Let $X_{A,k}$ and $X_{B,l}$ denote their vector-field coefficients
as in the definition of the contact norm, and set
\[
 M_{A,k}=
 \sup_{|p|_w<r,\,|z|<r}|f|
 +\sup_{|p|_w<r,\,|z|<r}\|X_{A,k}\|_E,
 \quad
 M_{B,l}=
 \sup_{|p|_w<r,\,|z|<r}|j|
 +\sup_{|p|_w<r,\,|z|<r}\|X_{B,l}\|_E.
\]
All normal derivatives below are Fr\'echet derivatives on
$P_w\times\mathbb C$, equipped with the maximum norm.

For the scalar component, use
$\mathscr L_FJ=X_F(J)+F_zJ$.
The coefficient of the contribution from this pair is
\[
 i\langle l,X_{A,k}^{q}\rangle j
 +D_{(p,z)}j\,[X_{A,k}^{p},X_{A,k}^{z}]
 +f_zj.
\]
On the normal ball of radius $r-d$, the Banach Cauchy estimates
give $\|D_{(p,z)}j\|\leqslant Cd^{-1}\sup|j|$ and
$|f_z|\leqslant d^{-1}\sup|f|$.
Since $|\langle l,X_{A,k}^{q}\rangle|
\leqslant |l|_1\|X_{A,k}^{q}\|_\infty$, this scalar coefficient
is bounded by $C(|l|_1+d^{-1})M_{A,k}M_{B,l}$.
Only one derivative loss is needed: the undifferentiated
vector field $X_F$ is already controlled by the input contact norm.

For the vector-field component, the contact convention implies
$X_{\mathscr L_FJ}=[X_F,X_J]$.
Indeed, the commutator is a contact vector field and
\[
 \alpha([X_F,X_J])
 =X_F\bigl(\alpha(X_J)\bigr)
  -(\mathcal L_{X_F}\alpha)(X_J)
 =-X_F(J)-F_zJ
 =-\mathscr L_FJ.
\]
Uniqueness of the admissible contact Hamiltonian vector field
therefore gives the stated identity. The coefficient of this
commutator contributed by the same pair is
\[
i\langle l,X_{A,k}^{q}\rangle X_{B,l} +D_{(p,z)}X_{B,l}[X_{A,k}^{p},X_{A,k}^{z}] -i\langle k,X_{B,l}^{q}\rangle X_{A,k} -D_{(p,z)}X_{A,k}[X_{B,l}^{p},X_{B,l}^{z}].
\]
Banach-valued Cauchy estimates thus bound its $E$ norm by
$C(|k|_1+|l|_1+d^{-1})M_{A,k}M_{B,l}$ on the smaller normal ball.
This argument also covers the common-momentum terms in the
vector-field coefficients; no coordinatewise summation of their
operator norms is used.

The scalar bracket has support contained in $A\cup B$ and
Fourier index $k+l$. Put $L=|k|_1+|l|_1$. The union-weight
inequality and $|k+l|_1\leqslant L$ give
\[
\begin{aligned}
 e^{m[A\cup B]+(s-\sigma)|k+l|_1}(L+d^{-1})
 &\leqslant e^{m[A]+s|k|_1}e^{m[B]+s|l|_1}e^{-\sigma L}(L+d^{-1})\\
 &\leqslant C(\sigma^{-1}+d^{-1})e^{m[A]+s|k|_1}e^{m[B]+s|l|_1}.
\end{aligned}
\]
Summing the scalar and vector-field bounds absolutely over all
pairs proves the claimed contact-norm estimate. Every constant
comes from a Banach Cauchy estimate or a continuous pairing and
is independent of the cardinalities of $A$ and $B$.
\end{proof}

The bracket estimate also controls the pullback of an arbitrary
spatial Hamiltonian by a small generating flow.

\begin{lemma}
\label{lem:infinite-pullback}
Let $0<\sigma<s$ and $0<d<r$, and suppose that $F,J$ have
finite contact norms on $\mathcal D_{s,r}$. There is a constant
$C>0$ such that, if
\[
\theta=Ce(\sigma^{-1}+d^{-1})|F|^c_{m,s,r}\leqslant\frac12,
\]
then the flow $\Psi_t$ of $X_F$ is defined on
$\mathcal D_{s-\sigma,r-d}$ for $0\leqslant t\leqslant1$.
Writing $\Psi_t^*\alpha=\widehat\rho_t\alpha$, its weighted
pullback satisfies
\begin{equation*}
\mathscr P_tJ=\widehat\rho_t^{-1}J\circ\Psi_t
=\sum_{N=0}^{\infty}\frac{t^N}{N!}\mathscr L_F^N J,
\quad
|\mathscr P_tJ|^c_{m,s-\sigma,r-d}
\leqslant2|J|^c_{m,s,r}.
\end{equation*}
\end{lemma}

\begin{proof}
For $N\geqslant1$, distribute the angular loss $\sigma$ and
normal loss $d$ equally among $N$ applications of
Lemma~\ref{lem:infinite-bracket}. This gives
\[
|\mathscr L_F^N J|^c_{m,s-\sigma,r-d}
\leqslant
\bigl[CN(\sigma^{-1}+d^{-1})\bigr]^N
\bigl(|F|^c_{m,s,r}\bigr)^N|J|^c_{m,s,r}.
\]
Since $N!\geqslant(N/e)^N$, the Lie series is bounded by
$(1-\theta)^{-1}|J|^c_{m,s,r}$. Increasing $C$ if necessary,
the same smallness condition bounds the displacement of
$X_F$ by a fixed fraction of both available domain margins,
so its Banach-space flow exists on the stated domain.
On intermediate smaller domains the series can be
differentiated termwise. It solves
$\partial_t\mathscr P_tJ=\mathscr L_F\mathscr P_tJ$ with
initial value $J$. The weighted pullback solves the same
transport equation; uniqueness along the characteristics
of $X_F$ identifies the two expressions. All the spatial
series converge absolutely by the displayed majorant.
\end{proof}

For the affine generator
\[
F(q,p,z)=u(q)+v(q)\cdot p+w(q)z,
\]
the momentum equation is
\[
\dot p
=
-Du-(Dv)^Tp-zDw-wp.
\]
Here $(Dv)^T$ is defined relative to the duality pairing
between $P_w$ and $V$. Its boundedness on $P_w$ follows from
the Banach-holomorphic $P_w$-valued momentum component of
$X_F$. Indeed, if $Y^p$ denotes that component, then
$D_pY^p=-(Dv)^T-wI$, where $w$ in this formula is the scalar
coefficient of $z$ in the generator. The Banach Cauchy estimate
therefore gives a bounded operator $(Dv)^T$ on $P_w$ on every
smaller normal ball. The same argument applies to the other
coefficients of the affine flow. For a sufficiently small
flow time, the corresponding linear operators are invertible
by the Neumann series. Solving the affine Banach-space flow
then gives the affine contact form \eqref{eq:affine}, with
$DQ^{-T}$ acting boundedly on $P_w$. The spatial expansions of the affine coefficients and their
finite compositions are controlled by
Lemma~\ref{lem:infinite-accumulated} below. General spatial
Hamiltonians are pulled back by Lemma~\ref{lem:infinite-pullback};
\eqref{eq:infinite-composition-hypothesis} will restore the
original data $R_0,P$ on the full working normal ball.

The contact pullback identity and the composition law
\eqref{eq:composition} follow from the continuous dual
pairing and the absolute convergence of the spatial
expansions. In particular, the limiting affine formula
has $DU,D\rho\in P_w$.

Finally, the bare linear functions
\[
L_\nu=\nu\cdot p,\quad \nu\in V,
\]
are treated separately from the spatial norm. A direct
calculation gives
\begin{equation}\label{eq:infinite-linearpart}
 \mathscr L_FL_\nu=-D_\nu F,
 \quad
 D_\nu=\sum_{i\in\mathcal I}\nu_i\partial_{q_i},
\end{equation}
and the Fourier multiplier estimate yields
\[
 |D_\nu F|^c_{m,s-\sigma,r}
 \leqslant
 (e\sigma)^{-1}\|\nu\|_\infty
 |F|^c_{m,s,r}.
\]
Indeed, on a Fourier mode $k$ the multiplier is
$i k\cdot\nu$, and
\[
 |k\cdot\nu|
 \leqslant \|\nu\|_\infty |k|_1,
 \quad
 |k|_1e^{-\sigma|k|_1}\leqslant(e\sigma)^{-1}.
\]
The same multiplier acts on the vector-field coefficient
because $D_\nu$ commutes with the coordinate derivatives
in the Hamiltonian vector-field formula. Thus the weighted
pullback of $L_\nu$ differs from $L_\nu$ by a convergent
spatial series on smaller domains.

We also use a bounded first-jet projection in the contact norm.
If $J$ is defined on $\mathcal D_{s,2r}$, then
\begin{equation}\label{eq:infinite-jet-bound}
 |\mathcal J_1J|^c_{m,s,r}
 \leqslant C|J|^c_{m,s,2r}.
\end{equation}
To see this, write $X_J=(X_J^q,X_J^p,X_J^z)$. The coordinate
formula gives
\[
 X_{\mathcal J_1J}
 =\bigl(X_J^q|_0,\mathcal J_1X_J^p,
                       \mathcal J_1X_J^z\bigr).
\]
Normal Cauchy estimates bound these coefficient projections on
each spatial Fourier block, and summation proves
\eqref{eq:infinite-jet-bound}. The constant depends on the fixed
radius $r$ and is independent of the spatial supports.

\subsection{Fourier inversion and the average operator}
\label{sec:infinite-inverse-average}

Put
\[
\Gamma(t)=\sup_{x\geqslant0}\Delta(x)e^{-tx},
\quad t>0.
\]
The assumption \eqref{eq:infinite-approximation} implies
that $\Gamma(t)<+\infty$ for every $t>0$.

For a zero-average spatial Fourier series
\[
f(q)=\sum_{A\in\mathcal S}
\sum_{k\in\Z^A\setminus\{0\}}
f_{A,k}e^{ik\cdot q},
\]
define
\[
T_\beta f
=
\sum_{A\in\mathcal S}
\sum_{k\in\Z^A\setminus\{0\}}
\frac{f_{A,k}}{\beta+ik\cdot\omega}e^{ik\cdot q},
\quad \beta\in\mathbb R.
\]
The zero Fourier mode is set equal to zero.

For every nonzero $k\in\mathbb Z^{(\mathcal I)}$,
the nonresonance condition gives
\[
\frac{1}{|\beta+ik\cdot\omega|}
\leqslant
\frac{\Delta([k])\Delta(|k|_1)}{\gamma}.
\]
Since $\operatorname{supp}k\subseteq A$ implies
$[k]\leqslant[A]$, we obtain
\[
e^{-\mu[A]-\sigma|k|_1}
\frac{1}{|\beta+ik\cdot\omega|}
\leqslant
\gamma^{-1}\Gamma(\mu)\Gamma(\sigma).
\]
Consequently,
\begin{equation}\label{eq:infinite-inverse}
 |T_\beta f|^c_{m-\mu,s-\sigma,r}
 \leqslant
 \gamma^{-1}\Gamma(\mu)\Gamma(\sigma)
 |f|^c_{m,s,r},
 \quad 0<\mu<m,\quad 0<\sigma<s.
\end{equation}

The same estimate holds for Banach-valued coefficients in the
corresponding spatial norms. The Fourier multiplier commutes
with all coordinate derivatives in the Hamiltonian vector-field
formula, so the contact norm requires no additional loss here.

We next consider the average equations. Write the first
jet of the current Hamiltonian as
\[
G(q,p,z)
=
\omega\cdot p+\beta z
+a(q)+b(q)\cdot p+c(q)z+R(q,p,z),
\quad R=O_2.
\]
With the fixed working radius $r$ chosen below, the defect is
measured by
\[
 E_{m,s}
 =
 \|\av b\|_\infty+
 |a+(b-\av b)\cdot p+cz|^c_{m,s,r}.
\]
The constant vector $\av b$ is treated separately as an
element of $V=\ell^\infty(\mathcal I)$.

Let the accumulated affine contact map be represented by
$(Q,\rho,U)$, and put
\[
 A_\Phi(q,p,z)
 =U(q)+(Q(q)-q)\cdot p+(\rho(q)-1)z,
 \quad
 Y=|A_\Phi|^{c,2}_{m,s,8r}.
\]
Here $|J|^{c,2}_{m,s,8r}$ is the norm
\eqref{eq:infinite-contactnorm} with the additional factor
$(1+|k|_1)^2$ in each summand. It controls two angular
derivatives in the corresponding Banach operator norms.
Writing $h=Q-\mathrm{id}$, we have
\[
 X_{A_\Phi}
 =\bigl(h,-DU-(Dh)^Tp-zD\rho-(\rho-1)p,
                  -U-(\rho-1)z\bigr).
\]
Coefficient extraction on the fixed normal ball therefore
controls $DU,D\rho\in P_w$, $Dh\in\mathcal L(V)$, and
$(Dh)^T\in\mathcal L(P_w)$ by $CY$.
Each spatial Fourier coefficient $h_{A,k}$ of $h$ takes values
in $V_A$ and depends only on $q_A$. Hence $Dh_{A,k}$ is an
$A$-local operator on $V$, and its transpose is $A$-local on
$P_w$. The scalar coefficient $(\rho-1)I$ is treated as a scalar
spatial multiplier times the identity. Coefficient extraction
from $X_{A_\Phi}$, together with the two angular weights in $Y$,
gives
\[
\|Dh\|_{\mathrm{op},m,s}
+\|(Dh)^T\|_{\mathrm{op},m,s}
+\|DU\|_{P_w,m,s}
+\|D\rho\|_{P_w,m,s}
\leqslant CY.
\]
For sufficiently small $Y$, the Neumann series therefore gives
\[
\|DQ^{-1}-I\|_{\mathrm{op},m,s}
+\|DQ^{-T}-I\|_{\mathrm{op},m,s}
\leqslant CY.
\]
These are the operator bounds used in the affine formula
\eqref{eq:affine}. Differentiating that formula once also gives
\[
\|D\Phi-I\|_{\mathrm{op},m,s,8r}\leqslant CY,
\quad
\|(D\Phi)^{-1}-I\|_{\mathrm{op},m,s,8r}\leqslant CY.
\]
Here we use the same spatial convolution norm for operators on
$E$, with the common momentum terms assigned as in
Remark~\ref{rem:infinite-composition}. The first bound uses
$D(DQ^{-T}),D^2U,D^2\rho$, all controlled by $Y$; the second
follows from the Neumann series.

For clarity, define the spatial operator norm by
\[
\|C\|_{\mathrm{op},m,s}
=\sum_{A\in\mathcal S}
\sum_{k\in\Z^A}
e^{m[A]+s|k|_1}\|C_{A,k}\|,
\]
where the operator norm is taken between the indicated Banach
spaces. Let $\Pi_A^V$ and $\Pi_A^P$ denote the coordinate
projections onto the coordinates in $A$ in $V$ and $P_w$,
respectively. A block operator $C_{A,k}:X\to Y$, with $X,Y$ equal
to $V$ or $P_w$, is called $A$-local if
$C_{A,k}=\Pi_A^Y C_{A,k}\Pi_A^X$. The same convention is used
for mixed operators. Scalar multiplication operators are treated
separately as scalar spatial coefficients times the identity.
Products of local blocks have spatial support contained in the
union of their supports. Thus the union-weight inequality gives
the submultiplicative estimate
\[
\|C_1C_2\|_{\mathrm{op},m,s}
\leqslant \|C_1\|_{\mathrm{op},m,s}
             \|C_2\|_{\mathrm{op},m,s}
\]
for compatible operator-valued spatial expansions, including
mixed operators. The identity is the unit of this algebra, so
Neumann series converge in the spatial operator norm whenever
the non-identity part has norm less than one.

The parameter response is still given by \eqref{eq:paramcoeff},
with $B(q)\in\mathcal L(V)$, $\ell(q),m(q)\in P_w$, and the scalar-valued function $\chi(q)=\rho(q)^{-1}$. The product estimate in the spatial operator algebra gives
\[
\|B-I\|_{\mathrm{op},m,s}
+\|\ell\|_{P_w,m,s}
+\|m\|_{P_w,m,s}
+\|\chi-1\|_{m,s}
\leqslant CY.
\]
We also need the contact norm of the parameter response. For an
affine Hamiltonian $f+g\cdot p+az$, coefficient extraction on
the fixed normal ball shows that its contact norm is equivalent,
with constants depending only on $r$, to the sum of the spatial
norms of $f,Df,g,(Dg)^T,a,Da$, in $\mathbb C,P_w,V,
\mathcal L(P_w),\mathbb C,P_w$, respectively. Put
$v_\eta=B\eta$ and $T_\eta=Dv_\eta$. Differentiating the inverse
matrix identity gives
\[
T_\eta=-B D^2h[\,\cdot,v_\eta],
\quad
D\chi=-\chi^2D\rho,
\quad
\|T_\eta^T\|_{\mathrm{op},m,s}
\leqslant CY\|\eta\|_\infty.
\]
Indeed, the extra Fourier weights in $Y$ control the directional
derivatives of $(Dh)^T$, $DU$ and $D\rho$, with no further
spatial loss. Since
$\ell\cdot\eta=\chi DU\cdot v_\eta$ and
$m\cdot\eta=\chi D\rho\cdot v_\eta$, we have
\[
\begin{aligned}
D(\ell\cdot\eta)
 &= (DU\cdot v_\eta)D\chi
    +\chi D^2U[\,\cdot,v_\eta]+\chi T_\eta^TDU,\\
D(m\cdot\eta)
 &= (D\rho\cdot v_\eta)D\chi
    +\chi D^2\rho[\,\cdot,v_\eta]+\chi T_\eta^TD\rho.
\end{aligned}
\]
All terms belong to the corresponding spatial $P_w$ class.
The product estimate and the preceding coefficient bounds yield
\begin{equation}\label{eq:infinite-param-spatial}
\left|\ell\cdot\eta-(\chi-1)\zeta
 +((B-I)\eta)\cdot p+(m\cdot\eta)z\right|^c_{m,s,8r}
\leqslant CY(\|\eta\|_\infty+|\zeta|).
\end{equation}

We eliminate the first average exactly as in
\eqref{eq:meanelim}:
\[
 a^\circ=a-\chi\frac{\av a}{\av\chi},
 \quad
 \ell^\circ=\ell-\chi\frac{\av\ell}{\av\chi},
 \quad
 u_0=T_\beta a^\circ,
 \quad
 \mathcal V\eta=T_\beta(\ell^\circ\cdot\eta).
\]
Thus both sources have zero average and $u=u_0+\mathcal V\eta$.

We estimate the average equations with fixed spatial and angular
margins, as in Lemma~\ref{lem:linear}. Suppose that the current
widths are bounded below by fixed positive numbers. Applying
\eqref{eq:infinite-inverse} across half of these lower bounds,
and then restricting to the real torus, gives
\begin{equation}\label{eq:infinite-real-estimates}
     \sup_{q\in\mathbb T^{\mathcal I}} \bigl(|u_0(q)|+|Du_0(q)|_w\bigr) \leqslant CE_{m,s},\quad \sup_{q\in\mathbb T^{\mathcal I}} \bigl(|\mathcal V\eta(q)|+|D(\mathcal V\eta)(q)|_w\bigr) \leqslant CY\|\eta\|_\infty.
\end{equation}
The constants are independent of the stepwise losses.
Normal Cauchy estimates in the contact norm give
$S(q)\in\mathcal L(P_w,V)$, $t(q)\in V$, and scalar $r_2(q)$,
with uniform bounds. Consequently, the products
$SDu$, $tu$, and $t\cdot Du$ in the correction equations
are well defined, and their real-torus bounds follow from
\eqref{eq:infinite-real-estimates}.

The average of the second correction equation is
therefore the operator equation
\begin{equation*}
\mathsf M\eta=-\av{b-SDu_0-tu_0}, \quad \mathsf M\eta=\av{B\eta-SD(\mathcal V\eta)-t(\mathcal V\eta)}.
\end{equation*}
Here
\[
\mathsf M:V\to V.
\]
At the identity accumulated map,
\[
B=I,\quad \ell=0,\quad \mathcal V=0,
\]
and hence $\mathsf M=I$.

The parameter coefficient bounds and
\eqref{eq:infinite-real-estimates} give
\[
 \|\mathsf M-I\|_{\mathcal L(V)}
 \leqslant\sup_{q\in\mathbb T^{\mathcal I}}\|B(q)-I\|_{\mathcal L(V)}
       +CY
 \leqslant CY.
\]
If $CY<1/2$, the Neumann series gives
\[
\mathsf M^{-1}
=
\sum_{j=0}^{\infty}(I-\mathsf M)^j,
\quad
\|\mathsf M^{-1}\|_{\mathcal L(V)}\leqslant2.
\]
Thus the average equation has a unique solution
$\eta\in V$, with
\[
\|\eta\|_\infty\leqslant CE_{m,s}.
\]

The first average equation determines
\[
\zeta
=
\frac{\av a+\langle\av\ell,\eta\rangle}{\av\chi},
\]
provided $\av\chi\ne0$. Since $\chi$ is close to $1$,
this denominator is bounded away from zero for sufficiently
small $Y$.

The average of the third correction equation gives
$\delta\beta=\av{c+m\cdot\eta-t\cdot Du-r_2u}$. The
bounds above and the uniform bounds for $m$, $t$, and $r_2$ imply
$|\delta\beta|\leqslant CE_{m,s}$. The remaining zero-average
equations are solved by $T_\beta$ and $T_0$ as in
\eqref{eq:usplit}--\eqref{eq:vwexplicit}. In particular,
the zero Fourier mode of
\[
(\av b+\eta)\cdot p
\]
is not inverted; it is removed by the zero-average
normalization and is controlled through the separate
parameter equation.

Thus the average equations determine the parameter corrections
in $V\times\mathbb R$ and the normal-form increment
$\delta\beta$.

We next prove the estimate for the increment of the
accumulated affine map. This estimate is needed because
the finite-dimensional composition estimate alone does
not control the spatially weighted norm used here.

\begin{lemma}
\label{lem:infinite-accumulated}
Let $\Phi$ be an accumulated affine contact map with
\[
Y=|A_\Phi|^{c,2}_{m,s,8r}\leqslant Y_*,
\]
and let $\Psi_t$ be the flow generated by
\[
F=u+v\cdot p+wz.
\]
Write $\Phi_t=\Phi\circ\Psi_t$, and let $A_t$ be the
affine coefficient associated with $\Phi_t$. There are
constants $c_0,C>0$, independent of the iteration index,
such that, if
\[
|F|^{c,2}_{m,s,8r}\leqslant c_0\sigma,
\quad 0<\sigma<\min\{1,s\},
\]
then
\[
|A_1-A_0|^{c,2}_{m,s-\sigma,8r}
\leqslant
C\sigma^{-1}(1+Y)|F|^{c,2}_{m,s,8r}.
\]
In particular,
\[
Y_+\leqslant
Y+C\sigma^{-1}|F|^{c,2}_{m,s,8r}.
\]
\end{lemma}

\begin{proof}
Write $A_t=U_t+(Q_t-\mathrm{id})\cdot p+(\rho_t-1)z$.
Differentiating the affine composition formula
\eqref{eq:composition} gives the exact equation
\[
\partial_tA_t=\mathcal B_FA_t+b_F,
\quad
\mathcal B_F=v\cdot D_q-(u+wz)\partial_z,
\quad b_F=-u+v\cdot p-wz.
\]
For $A=f+g\cdot p+az$, its image has coefficients
\[
\mathcal B_FA
 =(Df\cdot v-ua)+(Dg\,v)\cdot p+(Da\cdot v-wa)z.
\]
The contact norm controls the spatial norms of
$f,Df,g,(Dg)^T,a,Da$ by coefficient extraction on the fixed
normal ball. In differentiating the displayed coefficients,
we use
\[
D(Df\cdot v)=D^2f[\,\cdot,v]+(Dv)^TDf,
\quad
D(Dg\,v)^T=(D^2g[\,\cdot,v])^T+(Dv)^T(Dg)^T,
\]
and the analogous formula with $f$ replaced by $a$.
These identities involve bounded operators on $P_w$ and
continuous pairings with $V$. The union-weight inequality
controls their spatial products. Each second derivative costs
only one additional angular derivative relative to the
coefficient norms just listed. Thus, for
$0<\delta<\min\{1,s\}$,
\[
|\mathcal B_FA|^{c,2}_{m,s-\delta,8r}
\leqslant C\delta^{-1}|F|^{c,2}_{m,s,8r}
                         |A|^{c,2}_{m,s,8r},
\quad
|b_F|^{c,2}_{m,s,8r}\leqslant C|F|^{c,2}_{m,s,8r}.
\]
Here the factor $(1+|k|_1)^2$ is retained by
$(1+|k+l|_1)^2\leqslant(1+|k|_1)^2(1+|l|_1)^2$;
the remaining derivative is absorbed by $\delta$.
No normal loss is needed for these affine coefficients.

Put $f_*=|F|^{c,2}_{m,s,8r}$. The solution series separates as
\[
A_1-A_0
 =\sum_{n=1}^{\infty}\frac{\mathcal B_F^nA_0}{n!}
  +\sum_{n=0}^{\infty}\frac{\mathcal B_F^nb_F}{(n+1)!}.
\]
Distributing the angular loss $\sigma$ equally among the
$n$ operator factors and using $n!\geqslant(n/e)^n$ bounds
the two sums, respectively, by
\[
Y\sum_{n=1}^{\infty}(Cef_*/\sigma)^n,
\quad
Cf_*\sum_{n=0}^{\infty}(Cef_*/\sigma)^n.
\]
For $f_*\leqslant c_0\sigma$ with sufficiently small $c_0$,
these series converge and give
\[
|A_1-A_0|^{c,2}_{m,s-\sigma,8r}
\leqslant C(Yf_*/\sigma+f_*)
\leqslant C\sigma^{-1}(1+Y)f_*.
\]
The series solves the coefficient evolution equation and
therefore agrees with the affine composition. The final
assertion follows from the triangle inequality.
\end{proof}

\subsection{The step and convergence}

We now combine the spatial estimates with the correction
equations of Sections~\ref{sec:linear} and~\ref{sec:iteration}.
We use the fixed widths and normal radius from
Remark~\ref{rem:infinite-composition}. The current widths
$(m,s)$ and the next widths $(m-12\mu,s-12\sigma)$ are
assumed to lie in its stated range, with
$0<\mu<m/12$ and $0<\sigma<\min\{1,s/12\}$.
All constants below depend only on the fixed data, the
composition constants, and the positive lower bounds for
$m$ and $s$.

For an accumulated map with
$Y=|A_\Phi|^{c,2}_{m,s,8r}\leqslant Y_*$,
\eqref{eq:infinite-composition-hypothesis} gives
\[
|\rho^{-1}J\circ\Phi|^c_{m,s,8r}
\leqslant C_{\mathrm{sc}}|J|_{m_0,s_0,r_0},
\quad J=R_0,P.
\]
Thus the background bound always uses the original scalar
data norms. We keep the bare frequency separate throughout.
In particular, after the parameter change, put
$\nu=\xi+\eta$ and write
\[
G^\sharp=L_\nu+J^\sharp,\quad L_\nu=\nu\cdot p,
\]
where the affine pullback formula gives
\[
J^\sharp={}\ell\cdot\nu+((B-I)\nu)\cdot p +(m\cdot\nu)z+\lambda z+\lambda\chi U -(e+\zeta)\chi +\rho^{-1}(R_0+\varepsilon P)\circ\Phi.
\]
Here $m$ in $m\cdot\nu$ denotes the parameter coefficient
in \eqref{eq:paramcoeff}. The affine coefficient estimates,
\eqref{eq:infinite-param-spatial}, and the composition bound
therefore give
\[
\|\nu\|_\infty+|J^\sharp|^c_{m,s,8r}\leqslant M
\]
with a fixed $M$, as long as the external parameters remain
in fixed neighborhoods of their initial values. Since the parameter increment is affine,
$R=(I-\mathcal J_1)J^\sharp$. Projecting first from radius
$8r$ to $4r$ and extending the affine jet back to $8r$ gives
$|R|^c_{m,s,8r}\leqslant CM$. Normal Cauchy estimates on fixed
smaller balls then give the uniform Hessian bounds.

Write $E=E_{m,s}$. The average equations give
\[
\|\eta\|_\infty+|\zeta|+|\delta\beta|\leqslant CE.
\]
The modified affine defect has the decomposition
\[
\widetilde g=L_{\nu_g}+J_g,\quad
\nu_g=\av b+\eta,
\]
with
\[
J_g=a+\ell\cdot\eta-\chi\zeta
 +\bigl(b-\av b+(B-I)\eta\bigr)\cdot p
 +(c+m\cdot\eta)z.
\]
By \eqref{eq:infinite-param-spatial} and the definition of
$E$,
\[
\|\nu_g\|_\infty+|J_g|^c_{m,s,8r}\leqslant CE.
\]
The change from radius $r$ in the defect norm to radius
$8r$ costs only a fixed factor for these affine functions.

We specify the losses in solving the correction equations.
First, \eqref{eq:infinite-inverse} gives $u$ on
$(m-\mu,s-\sigma)$. The Hessian contribution is
$\mathcal J_1\mathscr L_uR$, as in \eqref{eq:LR}.
One bracket estimate and the first-jet estimate bound it on
$(m-\mu,s-2\sigma)$, using a fixed normal-radius loss.
The average equations remove the constant modes of the
sources for $v$ and $w$. In particular, the constant vector
$\av b+\eta$ cancels with the average of the remaining
source, so the zero-average source is spatially summable.
A second Fourier inverse gives $v,w$ on
$(m-2\mu,s-3\sigma)$. Finally,
\[
(1+|k|_1)^2e^{-\sigma|k|_1}
\leqslant C(1+\sigma^{-2})
\]
provides the two additional angular weights. Extending the
affine coefficients to the fixed normal radius $8r$ and
restricting the spatial width further, we obtain
\[
A_F:=|F|^{c,2}_{m-4\mu,s-4\sigma,8r}
\leqslant K_0(\mu,\sigma)E.
\]
Here $K_0$ is bounded by a fixed power of
$1+\mu^{-1}+\sigma^{-1}$ and
$1+\gamma^{-1}\Gamma(\mu)\Gamma(\sigma/2)$.

We next estimate the exact identity \eqref{eq:twoerrors}.
For a decomposition $G=L_\nu+J$, identity
\eqref{eq:infinite-linearpart} gives
\[
\mathscr L_FG=-D_\nu F+\mathscr L_FJ.
\]
After this first bracket every term belongs to the spatial
contact class. On the fixed spatial width $m-4\mu$, the
first bracket is estimated on
$\mathcal D_{s-5\sigma,7r}$ and the second on
$\mathcal D_{s-6\sigma,6r}$. Hence
\[
|\mathscr L_F\widetilde g|^c_{m-4\mu,s-5\sigma,7r}
\leqslant C(1+\sigma^{-1})A_FE,
\]
and
\[
|\mathscr L_F^2G^\sharp|^c_{m-4\mu,s-6\sigma,6r}
\leqslant CM(1+\sigma^{-1})^2A_F^2.
\]
Apply Lemma~\ref{lem:infinite-pullback} to the latter
spatial Hamiltonian, using the additional margins
$\sigma$ and $r$. If $A_F\leqslant c\sigma$, it bounds
$\mathscr P_t\mathscr L_F^2G^\sharp$ on
$\mathcal D_{s-7\sigma,5r}$. The first-jet projection then
uses a fixed normal-radius ratio to give, on radius $r$,
\[
E_+\leqslant C\bigl[(1+\sigma^{-1})A_FE
 +(1+\sigma^{-1})^2A_F^2\bigr].
\]
Zero-mode extraction is bounded in the same norms, so this
estimate also controls the separately recorded mean of the
new angular defect. The quadratic remainder is thus controlled by the general
generating-flow estimate.

Lemma~\ref{lem:infinite-accumulated}, applied on
$(m-4\mu,s-4\sigma)$ with one further angular loss,
gives
\[
Y_+\leqslant Y+C\sigma^{-1}A_F.
\]
We restrict both outputs to $(m-12\mu,s-12\sigma)$.
The preceding bounds are therefore dominated by
\[
\mathcal K(\mu,\sigma)
=C(1+\mu^{-1}+\sigma^{-1})^L
\left[1+\gamma^{-1}\Gamma(\mu)\Gamma(\sigma/2)\right]^L
\]
for fixed $C\geqslant1$ and an integer $L\geqslant1$.
For
\[
\mathcal K(\mu,\sigma)E\leqslant c\sigma,
\quad Y+\mathcal K(\mu,\sigma)E\leqslant Y_*,
\]
the step satisfies
\begin{equation}\label{eq:infinite-step}
E_+\leqslant\mathcal K(\mu,\sigma)E^2,
\quad
Y_+\leqslant Y+\mathcal K(\mu,\sigma)E.
\end{equation}
The same bound controls the increment of the accumulated
affine coefficients. On the new widths we define
\[
G_+=\rho_+^{-1}
H_{\xi+\eta,e+\zeta,\lambda,\varepsilon}\circ\Phi_+.
\]
We use half the threshold in Remark~\ref{rem:infinite-composition}
as the bound for $Y_+$ in the iteration, so the affine coefficient has
a representative within the threshold of that remark. Substitution
then bounds the spatial part of $G_+$ on the full ball of radius $8r$.
On the smaller ball the affine composition identity identifies $G_+$
with the flow pullback. Its first jet therefore has the representation
estimated by the integral remainder above, while projection of the
substituted Hamiltonian bounds its higher-order part on the working
ball. Thus the next step has
exactly the input bounds used above, uniformly for
$\lambda\in[-\Lambda,\Lambda]$.

We now choose the losses directly from
\eqref{eq:infinite-approximation}. Put
\[
 \varphi(x)=\log\Delta(x),\quad h(x)=\frac{\varphi(x)}{x},
 \quad 1<\vartheta<\kappa_\infty<2.
\]
Choose $T\geqslant1$ so large that $h$ is nonincreasing on $[T,\infty)$,
and set
\[
 \mu_j=h(T\vartheta^j)+T^{-1}\vartheta^{-j},
\quad \sigma_j=2\mu_j,\quad j\geqslant0.
\]
Monotonicity and the integral condition imply
\[
 \sum_{j=0}^{\infty}h(T\vartheta^j)
 \leqslant h(T)+\frac{1}{\log\vartheta}
              \int_T^\infty\frac{\varphi(x)}{x^2}\,dx
 \to0
 \quad(T\to\infty).
\]
For $x\geqslant T\vartheta^j$, one has $\varphi(x)\leqslant\mu_jx$;
for $0\leqslant x\leqslant T\vartheta^j$, one has
$\varphi(x)\leqslant\varphi(T\vartheta^j)$. Hence
\[
 \log\Gamma(\mu_j)=\log\Gamma(\sigma_j/2)
 \leqslant\varphi(T\vartheta^j)\leqslant Th(T)\vartheta^j,
 \quad
 \mu_j^{-1}\leqslant T\vartheta^j,
\quad \sigma_j^{-1}\leqslant T\vartheta^j.
\]
Taking $T$ sufficiently large therefore gives
\begin{equation}\label{eq:infinite-losses}
     12\sum_{j=0}^{\infty}\mu_j<\frac{\widehat m}{2},\quad 12\sum_{j=0}^{\infty}\sigma_j<\frac{\widehat s}{4},\quad \sum_{j=0}^{\infty}\kappa_\infty^{-j} \log\mathcal K(\mu_j,\sigma_j)<+\infty.
\end{equation}
We start at $m^{(0)}=\widehat m$ and $s^{(0)}=\widehat s/2$,
so the iterates retain a fixed margin from the original working
strip and satisfy $m^{(j)}>\widehat m/2$ and
$s^{(j)}>\widehat s/4$.

In particular, the limiting spatial and angular widths
remain positive.

Put
\[
K_j=\mathcal K(\mu_j,\sigma_j),
\quad
E_{j+1}\leqslant K_jE_j^2.
\]
Iteration gives
\begin{equation}\label{eq:infinite-recurrence}
E_j
\leqslant
\exp\left(
2^j\left[
\log E_0+
\sum_{i=0}^{j-1}2^{-i-1}\log K_i
\right]\right).
\end{equation}
Put
\[
 A=\sum_{j=0}^{\infty}2^{-j-1}\log K_j,
 \quad B=\sum_{j=0}^{\infty}\kappa_\infty^{-j}\log K_j,
 \quad \overline E=e^{-A-1}.
\]
Both sums are finite by \eqref{eq:infinite-losses}.
For $0<E_0\leqslant\overline E$, \eqref{eq:infinite-recurrence} gives
\[
 E_j\leqslant
 \left(\frac{E_0}{\overline E}\right)^{2^j}e^{-2^j}
 \leqslant\frac{E_0}{\overline E}e^{-2^j},
 \quad K_j\leqslant e^{B\kappa_\infty^j}.
\]
Consequently,
\begin{equation}\label{eq:infinite-weighted-sum}
 \sum_{j=0}^{\infty}K_jE_j
 \leqslant C_\Delta E_0,
 \quad
 C_\Delta=\overline E^{-1}
       \sum_{j=0}^{\infty}\exp(B\kappa_\infty^j-2^j)<+\infty.
\end{equation}
Since $\sigma_j\geqslant 2T^{-1}\vartheta^{-j}$, the same
estimates give
\[
\frac{K_jE_j}{\sigma_j}
\leqslant C T\vartheta^j
\exp(B\kappa_\infty^j-2^j)E_0.
\]
The sequence on the right is bounded after division by $E_0$,
so decreasing the initial defect ensures the stepwise flow
condition $K_jE_j\leqslant c_0\sigma_j$ at every iteration.
If $E_0=0$, every correction is zero. Otherwise, choose $E_0$
small enough that $C_\Delta E_0\leqslant Y_*/2$ and the
fixed average-operator and parameter bounds hold. The estimate above applies
on every finite segment of the construction and verifies the
next-step hypotheses by induction, exactly as in
Section~\ref{sec:iteration}.

It follows from the second estimate in
\eqref{eq:infinite-step} that the accumulated affine maps
converge in the corresponding weighted analytic norms.
The parameter corrections also converge:
\[
\sum_{j=0}^{\infty}\|\eta_j\|_\infty<+\infty,\quad
\sum_{j=0}^{\infty}|\zeta_j|<+\infty,\quad
\sum_{j=0}^{\infty}|\delta\beta_j|<+\infty.
\]
The limiting spatial and angular widths are positive by
\eqref{eq:infinite-losses}. Initially,
$\Phi_0=\mathrm{id}$, $\xi_0=\omega$, $e_0=0$, and
$\beta_0=\lambda$, so $E_0\leqslant C|\varepsilon|$.
The estimate \eqref{eq:infinite-weighted-sum} therefore gives
$O(|\varepsilon|)$ bounds for the parameter and accumulated-map
increments. Banach-valued Cauchy estimates on a fixed smaller
angular strip, followed by the affine formula \eqref{eq:affine},
give \eqref{eq:infinite-bound}.

The limiting affine map is real on real points and has
the form \eqref{eq:affine}. Its contact factor is positive
on the real domain because it is sufficiently close to $1$.
Writing $Q=\mathrm{id}+h$, the bound $\|Dh\|<1$ makes
$x\mapsto q-h(x)$ a contraction on the real angular lift.
It gives a unique inverse for $Q$, compatible with integer
translations. Thus $Q$ is a real-analytic diffeomorphism of
$\mathbb T^{\mathcal I}$.
The limiting affine coefficient remains $c,2$-small on the retained
widths. Hence the spatial composition hypothesis applies directly
to the limiting map and supplies spatial expansions of
$\rho^{-1}R_0\circ\Phi$ and $\rho^{-1}P\circ\Phi$ in the stated
contact class. The affine parameter terms are controlled by the
same coefficient and operator estimates as at each finite stage.

To pass the first-jet equations to the limit, fix a strictly
smaller spatial width and angular strip and a normal ball whose
closure lies inside the working normal ball. On this domain the
affine maps and their contact factors converge uniformly, and the
contact factors stay bounded away from zero. Their images have a
fixed positive margin inside the original analytic domain.
The Cauchy bounds for the fixed data $R_0$ and $P$, together with
the convergence of $\xi_j$ and $e_j$, therefore imply
\[
 G_j\to
 G_\infty:=\rho^{-1}H_{\xi,e,\lambda,\varepsilon}\circ\Phi
\]
uniformly on this smaller domain. Normal Cauchy estimates pass
the convergence to the first normal derivatives. Since
$E_j\to0$ and $\beta_j\to\beta$, it follows that
$\mathcal J_1G_\infty=\omega\cdot p+\beta z$.
The composition hypothesis gives the required spatial class of
the limit, while this fixed-domain argument identifies its first
jet. Consequently,
\[
\rho^{-1}H_{\xi,e,\lambda,\varepsilon}\circ\Phi
 =\omega\cdot p+\beta z+R_*,
 \quad R_*=O_2.
\]

The zero section is therefore invariant for the transformed
Hamiltonian and carries the linear flow
\[
\dot q=\omega,\quad \dot p=0,\quad \dot z=0.
\]
By the contact conjugacy identity, the embedding
\[
K(\theta)=\Phi(\theta,0,0)
\]
satisfies
\[
X_H\circ K=DK\,\omega,\quad K^*\alpha=0.
\]
Since the transformed Hamiltonian vanishes on the zero
section, $H\circ K=0$. The affine formula gives
\[
 K(\theta)
 =\bigl(Q(\theta),DQ(\theta)^{-T}DU(\theta),U(\theta)\bigr).
\]
With $\mathfrak u=U\circ Q^{-1}$, the chain rule identifies
its image with $j^1\mathfrak u$ and gives
$D\mathfrak u\in P_w$. The angular component of $K$ is the analytic
diffeomorphism $Q$ of $\mathbb T^{\mathcal I}$, so the torus is
full-dimensional.

Finally, differentiating the normal-form identity in the
$z$ direction along the zero section gives
\[
\beta=H_z\circ K+D_\omega\log\rho.
\]
Taking the zero Fourier mode yields
\[
\beta=\av{H_z\circ K}.
\]

All estimates above are uniform for real
$\lambda\in[-\Lambda,\Lambda]$. We give the parameter-continuity
argument on fixed strictly smaller domains. For real
$\beta,\beta'$ and a nonzero Fourier index $k$,
\[
 \frac{1}{\beta+ik\cdot\omega}
 -\frac{1}{\beta'+ik\cdot\omega}
 =\frac{\beta'-\beta}
 { (\beta+ik\cdot\omega)(\beta'+ik\cdot\omega)}.
\]
The nonresonance condition and the definition of $\Gamma$ imply,
for $0<\mu<m$ and $0<\sigma<s$,
\[
|(T_\beta-T_{\beta'})f|^c_{m-\mu,s-\sigma,r} \leqslant \gamma^{-2}|\beta-\beta'| \Gamma(\mu/2)^2\Gamma(\sigma/2)^2 |f|^c_{m,s,r}.
\]
Indeed, for $\operatorname{supp}k\subseteq A$, one has
$[k]\leqslant[A]$, and each of the two divisor factors is
absorbed using half of the available spatial and angular losses.
The multiplier commutes with the derivatives in the contact
vector-field formula, so no further derivative loss is needed.
Combining this estimate with the bound for $T_\beta$ gives
\[
|T_\beta f-T_{\beta'}f'|^c_{m-\mu,s-\sigma,r} \leqslant \gamma^{-1}\Gamma(\mu)\Gamma(\sigma) |f-f'|^c_{m,s,r} + \gamma^{-2}|\beta-\beta'| \Gamma(\mu/2)^2\Gamma(\sigma/2)^2 |f'|^c_{m,s,r}.
\]
Thus the Fourier inverses are jointly continuous in the source
and the real damping parameter, including at $\beta=0$.

At a fixed finite stage, coefficient extraction, averaging and
the spatial products are continuous on smaller domains.
The average operator is uniformly invertible, so its inverse
is continuous as well. The preceding joint estimate and the
integral equations for the generating flow then give, by
induction, continuity of the parameter corrections and the
accumulated affine coefficients at every finite stage.
All these assertions use fixed positive margins. On a smaller normal ball, Lemma~\ref{lem:infinite-pullback}
expresses the updated spatial part as a uniformly convergent series
of continuous iterated brackets, with \eqref{eq:infinite-linearpart}
giving the first bracket with the bare frequency. This gives continuity
in the spatial quotient norm. To recover it on any fixed ball of
radius $b<8r$, expand in homogeneous normal terms. For a term of
degree $d\geqslant1$, the angular component of its vector field has
degree $d-1$ and the other components have degree $d$.
The uniform contact bound on the ball of radius $8r$ thus gives a
geometric tail on radius $b$, while each finite homogeneous sum
extends continuously from the smaller ball by homogeneity.

The numerical majorants for the increments are independent of
$\lambda$ and summable. Their tails therefore tend to zero
uniformly in $\lambda$. The finite-stage continuity just proved
passes to the limiting parameters and affine coefficients by
the uniform-limit argument. Banach-valued Cauchy estimates on
a further fixed angular margin, followed by the affine formula
and the uniform bounds for the inverse operators, yield
continuity of the limiting embedding in the full-space $C^2$
norm on a common smaller complex domain. This completes the proof.

\section{Analytic Hamilton--Jacobi solutions}
\label{sec:hamiltonjacobi}

In both finite and infinite dimensions, the invariant Legendrian torus
is a graph over the entire angular base. Its graph function therefore
solves the stationary Hamilton--Jacobi equation globally and analytically.
Its characteristic flow has the prescribed frequency in the
original time. In finite dimensions we also identify the weak KAM
solutions and the contact Aubry and Mather sets under global convexity
and monotonicity assumptions.

\subsection{Analytic solutions with prescribed characteristic frequency}
For the selected finite-dimensional Hamiltonian in
\eqref{eq:family}, the invariant torus projects diffeomorphically
onto $\T^n$. The identity $H\circ K=0$ gives the stationary equation, and the
angular component of its dynamics prescribes the characteristic flow.

\begin{corollary}\label{cor:hamiltonjacobi}
Under the assumptions of Theorem~\ref{thm:main}, there exist
$s_1>0$, $\varepsilon_1>0$, and $C>0$ such that, for every real
$\lambda\in[-\Lambda,\Lambda]$ and $|\varepsilon|<\varepsilon_1$,
one can choose $\xi\in\R^n$, $e\in\R$, a holomorphic function
$\mathfrak u$ on $\T_{s_1}^n$, real on $\T^n$, and a
real-analytic diffeomorphism $Q:\T^n\to\T^n$ whose displacement
$Q-\mathrm{id}$ extends holomorphically to $\T_{s_1}^n$, satisfying
\begin{equation}\label{eq:hjapplication}
    \lambda\mathfrak u+\xi\cdot D\mathfrak u
    +R_0(q,D\mathfrak u,\mathfrak u)
    +\varepsilon P(q,D\mathfrak u,\mathfrak u)=e.
\end{equation}
For the selected Hamiltonian $H=H_{\xi,e,\lambda,\varepsilon}$,
the characteristic field
$b(q)=H_p(q,D\mathfrak u(q),\mathfrak u(q))$ satisfies
\begin{equation}\label{eq:hjfrequency}
    b(Q(\theta))=DQ(\theta)\omega,
    \quad q(t)=Q(\theta_0+\omega t).
\end{equation}
Moreover,
\[
    |\xi-\omega|+|e|+\|\mathfrak u\|_{C^2(s_1)}
    +\|Q-\mathrm{id}\|_{C^2(s_1)}
    \leqslant C|\varepsilon|.
\]
The constants are uniform in $\lambda$. For each fixed
$\varepsilon$, the selected parameters and functions can be
chosen continuously in $\lambda$, with the functions continuous
in the $C^2$ norm on a common smaller complex strip.
\end{corollary}

\begin{proof}
Let $(Q,\rho,U)$ be the affine coefficients supplied by
Theorem~\ref{thm:main}. Choose $s_1<s_*$ with a fixed positive
margin. After decreasing $\varepsilon_1$, the estimate for
$Q-\mathrm{id}$ gives a holomorphic inverse satisfying
$Q^{-1}(\T_{s_1}^n)\subset\T_{s_*}^n$, with uniform bounds
on a fixed smaller strip. Set $\mathfrak u=U\circ Q^{-1}$.
The chain rule gives
\[
    D\mathfrak u(Q(\theta))=DQ(\theta)^{-T}DU(\theta),
    \quad
    K(\theta)=
    \bigl(Q(\theta),D\mathfrak u(Q(\theta)),
                    \mathfrak u(Q(\theta))\bigr).
\]
Thus $H\circ K=0$ is precisely \eqref{eq:hjapplication}.
Taking the angular component of $X_H\circ K=DK\,\omega$
gives \eqref{eq:hjfrequency}, with no change of time.
The estimates follow from \eqref{eq:mainbound}, the chain rule,
and Cauchy's estimate across a fixed angular margin.
The same inverse and composition estimates transfer the
parameter continuity in Theorem~\ref{thm:main} to
$\mathfrak u$ and $Q$ on a common smaller strip.
\end{proof}

For example, take $R_0=\tfrac12\sum_{i=1}^n p_i^2$ and
$P=V(q)$, with $V$ real analytic. Then
\begin{equation}\label{eq:discountedhj}
    \lambda\mathfrak u+\xi\cdot D\mathfrak u
    +\frac12\sum_{i=1}^n(\partial_{q_i}\mathfrak u)^2+\varepsilon V=e,
    \quad \xi+D\mathfrak u(Q(\theta))=DQ(\theta)\omega.
\end{equation}
For this model, replacing $\mathfrak u$ by
$\mathfrak u-\av{\mathfrak u}$ and $e$ by
$e-\lambda\av{\mathfrak u}$ preserves both equations.
We may therefore impose $\av{\mathfrak u}=0$; averaging the
first equation then gives
$e=\varepsilon\av V+
\tfrac12\av{\sum_{i=1}^n(\partial_{q_i}\mathfrak u)^2}$.
Here the average is taken in the base coordinate $q$.

For $\lambda>0$, \eqref{eq:discountedhj} is a discounted
Hamilton--Jacobi equation. Fix a sufficiently small
$\varepsilon$ and write the selected normalized branch as
$(\xi_\lambda,e_\lambda,\mathfrak u_\lambda,Q_\lambda)$.
Corollary~\ref{cor:hamiltonjacobi} gives continuity in
$\lambda$ on a common smaller strip. In particular, as
$\lambda\downarrow0$, the functions $\mathfrak u_\lambda$
converge in $C^2$ on that strip to $\mathfrak u_0$, and the
parameters converge to $\xi_0$ and $e_0$. Passing to the
limit in \eqref{eq:discountedhj} gives
\[
    \xi_0\cdot D\mathfrak u_0
    +\frac12|D\mathfrak u_0|^2+\varepsilon V=e_0,
    \quad \av{\mathfrak u_0}=0.
\]
The limit follows the jointly selected branch
$(\xi_\lambda,e_\lambda,\mathfrak u_\lambda)$ to zero damping.
For $\lambda>0$, the identity $H_z=\lambda$ gives the strict
monotonicity used in the finite-dimensional weak KAM result below.

\subsection{Hamilton--Jacobi equations on the infinite-dimensional torus}

For countably many angles, the action space $P_w$ represents the
derivative of the solution through its pairing with $V=\ell^\infty(\mathcal I)$.

\begin{corollary}\label{cor:infinite-hj}
Under the assumptions of Theorem~\ref{thm:infinite}, there exist
$s_1,\varepsilon_1,C>0$ such that, for every real
$\lambda\in[-\Lambda,\Lambda]$ and $|\varepsilon|<\varepsilon_1$,
the selected $\xi,e$ admit a holomorphic solution $\mathfrak u$ on
$\T_{s_1}^{\mathcal I}$, real on $\T^{\mathcal I}$, of
\[
 \lambda\mathfrak u+\sum_{i\in\mathcal I}\xi_i\partial_{q_i}\mathfrak u
 +R_0(q,D\mathfrak u,\mathfrak u)+\varepsilon P(q,D\mathfrak u,\mathfrak u)=e.
\]
Its gradient is holomorphic with values in $P_w$ and satisfies
\[
 \sup_{q\in\T_{s_1}^{\mathcal I}}
 \bigl(|\mathfrak u(q)|+|D\mathfrak u(q)|_w\bigr)\leqslant C|\varepsilon|.
\]
The characteristic field $b(q)=H_p(q,D\mathfrak u(q),\mathfrak u(q))$
has the complete flow
$q(t)=Q(Q^{-1}(q_0)+\omega t)$ through every
$q_0\in\T^{\mathcal I}$, where $Q$ is the analytic angular
diffeomorphism in Theorem~\ref{thm:infinite}. The strip and bounds
are uniform in $\lambda$. For fixed $\varepsilon$ the solution and
its $P_w$-valued gradient depend continuously on $\lambda$ on a common
smaller strip.
\end{corollary}

\begin{proof}
Write $Q=\mathrm{id}+h$ and take $s_1=s_*/2$. By
\eqref{eq:infinite-bound}, we may choose $\varepsilon_1$ so that
$\|h\|_{s_*}<s_*/4$ and $\|Dh\|_{s_*}<1/2$. For each lift of
$q\in\T_{s_1}^{\mathcal I}$ the map $\theta\mapsto q-h(\theta)$
contracts the closed ball of radius $s_*/4$ about $q$ into itself.
Its iterates converge uniformly to a holomorphic inverse, compatible
with period translations.
Put $\mathfrak u=U\circ Q^{-1}$.
The chain rule in the $P_w,V$ pairing identifies
$D\mathfrak u(Q(\theta))=DQ(\theta)^{-T}DU(\theta)$ with the
momentum component of $K(\theta)$. Thus \eqref{eq:infinite-bound}
controls both $\mathfrak u$ and its $P_w$-valued gradient, and
\[
 D\mathfrak u(q)[v]=\sum_{i\in\mathcal I}\partial_{q_i}\mathfrak u(q)v_i,
 \quad |D\mathfrak u(q)[v]|\leqslant |D\mathfrak u(q)|_w\|v\|_\infty,
 \quad v\in V.
\]
In particular the series for $v=\xi$ converges absolutely.
Substitution in $H\circ K=0$ gives the equation. The angular component
of $X_H\circ K=DK\,\omega$ gives $b(Q(\theta))=DQ(\theta)\omega$,
which yields the stated flow in the original time for all $t\in\R$.
On a fixed smaller strip the contraction estimate gives continuous
dependence of $Q^{-1}$ on $\lambda$. Composing the vertical and
momentum components of $K$ with this inverse transfers the continuity
in Theorem~\ref{thm:infinite} to $\mathfrak u$ and its gradient.
\end{proof}

For the data $R_0=az^2$ with $a\ne0$ and $P=\sum_{i\in\Z}a_i\cos q_i$
in Remark~\ref{rem:infinite-composition}, the characteristic field
is the constant $\xi$. Taking the zero Fourier mode of
$\xi=DQ(\theta)\omega$ gives $\xi=\omega$. Thus the corollary yields
an analytic solution of the nonlinear equation
\[
    \sum_{i\in\Z}\omega_i\partial_{q_i}\mathfrak u
    +\lambda\mathfrak u+a\mathfrak u^2
    +\varepsilon\sum_{i\in\Z}a_i\cos q_i=e
\]
whose characteristics are the translations $q(t)=q_0+\omega t$.
For $\varepsilon\ne0$, the solution depends on every angle $q_i$:
independence of $q_i$ would leave $\varepsilon a_i\cos q_i$ as the
only term depending on that angle, and $a_i\ne0$.

More generally, if the Hamiltonian selected in
Corollary~\ref{cor:infinite-hj} admits a real $C^1$ extension to the
whole action space that is convex in $p$, its characteristics calibrate
$\mathfrak u$ in both time directions. For real $v\in V$, set
\[
    L(q,z,v)=\sup_{p\in P_w\cap\R^{\mathcal I}}
       \bigl\{\langle p,v\rangle-H(q,p,z)\bigr\}.
\]
The Hamilton--Jacobi equation gives
$L(q,\mathfrak u(q),v)\geqslant D\mathfrak u(q)[v]$. For $v=b(q)$,
convexity gives equality, with the supremum attained at
$p=D\mathfrak u(q)$. Along every characteristic $\gamma$ we therefore have
\[
    \mathfrak u(\gamma(t_1))-\mathfrak u(\gamma(t_0))
    =\int_{t_0}^{t_1}L(\gamma(t),\mathfrak u(\gamma(t)),\dot\gamma(t))\,dt,
    \quad -\infty<t_0<t_1<+\infty.
\]

\subsection{Analytic classical solutions: weak KAM and Aubry--Mather consequences}

We now derive consequences of the analytic classical solution
for weak KAM theory and Aubry--Mather sets. Under global Tonelli
and strict monotonicity assumptions, this solution is the unique
forward and backward weak KAM solution, and its analytic
Legendrian graph coincides with both the Aubry and Mather sets.
The dynamics on this graph is analytically conjugate, in the
original time, to the prescribed translation.

We use the contact weak KAM theory developed in
\cite{WWY2019}. Fix one of the selected finite-dimensional
Hamiltonians $H$ and assume that it admits a $C^3$ extension
to $T^*\T^n\times\R$, also denoted by $H$, satisfying
\begin{equation}\label{eq:hjtonelli}
    H_{pp}>0,\quad
    \lim_{|p|\to+\infty}\frac{H(q,p,z)}{|p|}=+\infty,
    \quad 0<H_z\leqslant C_0.
\end{equation}
Superlinearity is uniform for $(q,z)$ in compact sets.
These hypotheses place the selected Hamiltonian in the strictly
increasing contact Tonelli class.

The contact Lagrangian is
$L(q,z,v)=\sup_{p\in\R^n}\{p\cdot v-H(q,p,z)\}$.
The supremum is attained at the unique momentum
$p(q,z,v)$ satisfying $H_p(q,p(q,z,v),z)=v$.
In particular,
$L_z(q,z,v)=-H_z(q,p(q,z,v),z)<0$.
The constant term $-e$ is included in $H$, and hence the
corresponding term $+e$ is included in $L$.

A continuous function $f$ is dominated by $L$ if, for every
piecewise $C^1$ curve $\gamma:[a,b]\to\T^n$,
\begin{equation}\label{eq:hjdomination}
    f(\gamma(b))-f(\gamma(a))
    \leqslant
    \int_a^b L(\gamma(t),f(\gamma(t)),\dot\gamma(t))\,dt.
\end{equation}
A $C^1$ curve is $f$-calibrated if equality holds on every
compact subinterval of its domain.
A dominated function is a backward, respectively forward,
weak KAM solution if every point admits an $f$-calibrated
curve ending there on $(-\infty,0]$, respectively starting
there on $[0,+\infty)$.
We denote these solution sets by $\mathcal S_-$ and
$\mathcal S_+$, and use the contact Aubry set
$\widetilde{\mathcal A}$, Mather measures, and Mather set
$\widetilde{\mathcal M}$ as defined in \cite{WWY2019}.

\begin{proposition}\label{prop:hjweak}
Let $\mathfrak u,Q$ be given by
Corollary~\ref{cor:hamiltonjacobi}, and assume
\eqref{eq:hjtonelli} for the selected Hamiltonian $H$.
Then $\mathfrak u$ is the unique continuous viscosity
solution of $H(q,Df,f)=0$, and
\[
    \mathcal S_-=\mathcal S_+=\{\mathfrak u\},
    \quad
    \widetilde{\mathcal A}=\widetilde{\mathcal M}
       =j^1\mathfrak u=K(\T^n).
\]
The unique Mather measure is
$\mu=K_*\bigl(d\theta/(2\pi)^n\bigr)$, and
$j^1\mathfrak u$ is locally attracting.
All uniqueness and attraction assertions concern this
fixed Hamiltonian, including its selected parameters
$\xi$ and $e$.
\end{proposition}

\begin{proof}
Fenchel's inequality and
$H(q,D\mathfrak u(q),\mathfrak u(q))=0$ imply
$L(q,\mathfrak u(q),v)\geqslant D\mathfrak u(q)\cdot v$,
with equality precisely when $v=b(q)$.
Integration proves \eqref{eq:hjdomination} for $\mathfrak u$.
The smooth field $b$ is complete on the compact torus, and
its trajectories calibrate $\mathfrak u$ in both time
directions. Hence $\mathfrak u\in\mathcal S_-\cap\mathcal S_+$.

The classical solution $\mathfrak u$ is a viscosity solution,
so \cite[Theorem~B.1, p.~1020]{WWY2019} implies the
admissibility required in that theory.
By \cite[Proposition~A.1, p.~1019]{WWY2019},
$\mathcal S_-=\{\mathfrak u\}$.
The equivalence between backward weak KAM solutions and
viscosity solutions in
\cite[Proposition~2.7, p.~997]{WWY2019} then proves the
claimed viscosity uniqueness.

Let $f\in\mathcal S_+$. By
\cite[Theorem~1.2, p.~987]{WWY2019}, $f\leqslant\mathfrak u$.
The coincidence set
$I_f=\{q\in\T^n:f(q)=\mathfrak u(q)\}$ is nonempty.
Indeed, a forward $f$-calibrated curve lifts to a positive
orbit in the compact pseudograph $G_f$, and its nonempty
$\omega$-limit set is contained in the lifted coincidence
set $\widetilde I_f$ by
\cite[Theorems~1.1 and~1.5(2), pp.~987, 990]{WWY2019};
see also \cite[Remark~4.2, p.~1010]{WWY2019}.

Choose $q_0\in I_f$ and a forward $f$-calibrated curve
$\gamma$ with $\gamma(0)=q_0$.
Set $\delta=\mathfrak u-f\geqslant0$.
For $0\leqslant a<b$, subtracting the calibration identity
for $f$ from the domination inequality for $\mathfrak u$
gives
\[
    \delta(\gamma(b))-\delta(\gamma(a))\leqslant
    \int_a^b\bigl[L(\gamma(t),\mathfrak u(\gamma(t)),\dot\gamma(t))
    -L(\gamma(t),f(\gamma(t)),\dot\gamma(t))\bigr]\,dt\leqslant0.
\]
The last inequality follows from $L_z<0$.
Thus $\delta\circ\gamma$ is nonnegative and nonincreasing.
Since $\delta(q_0)=0$, it vanishes for every $t\geqslant0$,
so the entire positive curve remains in $I_f$.
Substituting $f=\mathfrak u$ along $\gamma$ into the
calibration identity on any interval $[a,b]$ shows that
$\gamma$ also calibrates $\mathfrak u$.
This is the coincidence-set propagation described by
\cite[Theorem~1.4, p.~990; Lemma~4.7 and
Remark~4.2, pp.~1008--1010]{WWY2019}.

By the chain rule and this calibration identity,
\[
    \int_a^b
    \bigl[
    L(\gamma(t),\mathfrak u(\gamma(t)),\dot\gamma(t))
    -D\mathfrak u(\gamma(t))\cdot\dot\gamma(t)
    \bigr]\,dt=0 .
\]
The integrand is continuous and nonnegative by Fenchel's
inequality. It therefore vanishes identically.
Strict convexity gives $\dot\gamma=b(\gamma)$, and
\eqref{eq:hjfrequency} yields
$\gamma(t)=Q(\theta_0+\omega t)$.
Nonresonance of $\omega$ makes this positive semiorbit
dense in $\T^n$. Since $f$ and $\mathfrak u$ agree on it,
continuity implies $f=\mathfrak u$ everywhere.
Hence $\mathcal S_+=\{\mathfrak u\}$.

By \cite[Theorem~1.3, p.~988]{WWY2019}, the Aubry set is
the intersection of the pseudographs of the backward
solution and the maximal forward solution.
Both solutions equal the analytic function $\mathfrak u$,
so $\widetilde{\mathcal A}=j^1\mathfrak u$.
The graph flow is conjugate to the nonresonant translation
flow, whose unique invariant probability is normalized
Haar measure. Thus its unique invariant probability is
$\mu=K_*\bigl(d\theta/(2\pi)^n\bigr)$, and
$\operatorname{supp}\mu=j^1\mathfrak u$.
Mather measures are the invariant Borel probabilities
supported in the Aubry set, and the Mather set is the
closure of the union of their supports; see
\cite[Section~3.2, p.~1001, after Definition~3.3]{WWY2019}.
This proves the assertions about the Mather measure and set.

Finally, the damping identity gives
$\beta=\av{H_z\circ K}>0$, with the average taken in the
parameter angle $\theta$.
The local contact linearization and asymptotic phase in
Section~\ref{sec:linearization} imply local attraction.
\end{proof}

Minimality of the characteristic flow propagates equality with
$\mathfrak u$ to the entire base torus. For each selected Hamiltonian
the attraction rate is determined by the mean damping
$\beta=\av{H_z\circ K}$. In the mechanical model
\eqref{eq:discountedhj}, this gives $\beta=\lambda$.

The dependence of the attraction rate on damping is explicit for
$H_\lambda(q,p,z)=\omega\cdot p+\tfrac12|p|^2+\lambda z$.
The invariant graph is $j^1 0$, and solutions with
$p(0)=0$ satisfy $q(t)=q(0)+\omega t$ and
$z(t)=z(0)e^{-\lambda t}$.
Consequently no positive exponential decay rate can be
chosen uniformly for $\lambda\downarrow0$.
At $\lambda=0$, every nonzero initial value $z(0)$ remains
constant, so this graph is not attracting; for $\lambda<0$,
these transverse displacements grow.

For $\lambda>0$, the Hamiltonian associated with
\eqref{eq:discountedhj} satisfies \eqref{eq:hjtonelli}, so
Proposition~\ref{prop:hjweak} applies.
The variational and KAM constructions then determine the
same solution: its characteristics are calibrated, and the
KAM conjugacy specifies their analytic parameterization
and frequency.
The example in Section~\ref{sec:comparison} illustrates
the role of the frequency requirement.
There $\mathfrak u=0$ is a smooth solution calibrated in
both directions, but its characteristic frequency
$\omega+\varepsilon e_1$ becomes resonant for arbitrarily
small $\varepsilon$.
The adjustable parameter $\xi$ restores the prescribed
rotation.
\appendix
\section{Supplementary estimates for the KAM iteration}
\label{app:verification}

\subsection{Evolution of the first jet}\label{app:jets}

Fix one KAM step and write $G^\sharp=G+\Delta G$ for the
parameter-updated Hamiltonian in the proof of
Lemma~\ref{lem:step}. Let $F=u+v\cdot p+wz$ be the corresponding
generator. For its weighted pullback, set
$G_\tau=\mathscr P_\tau G^\sharp$, where $0\leqslant\tau\leqslant1$.
The subscript $\tau$ in this subsection denotes time within this
single generating flow. Write
\[
    \mathcal J_1G_\tau=a_\tau+h_\tau\cdot p+d_\tau z,\quad
    S_\tau=(G_\tau)_{pp}|_0,\quad
    T_\tau=(G_\tau)_{pz}|_0,\quad
    \mathfrak r_\tau=(G_\tau)_{zz}|_0.
\]
The identity $\partial_\tau G_\tau=\mathscr L_FG_\tau$ gives
\begin{align}
    \partial_\tau a_\tau
        &=v\cdot Da_\tau-Du\cdot h_\tau-u d_\tau+w a_\tau,
        \label{eq:jetconstant}\\
    \partial_\tau h_\tau
        &=Dh_\tau\,v-Dv\,h_\tau-S_\tau Du-T_\tau u,
        \label{eq:jetp}\\
    \partial_\tau d_\tau
        &=v\cdot Dd_\tau-h_\tau\cdot Dw
           -T_\tau\cdot Du-\mathfrak r_\tau u.
        \label{eq:jetz}
\end{align}
For example,
$(G_\tau)_p=h_\tau+S_\tau p+T_\tau z+O_2$ and
$(G_\tau)_z=d_\tau+T_\tau\cdot p+\mathfrak r_\tau z+O_2$.
Substitution into $\mathscr L_FG_\tau=X_FG_\tau+wG_\tau$
yields the three identities. The terms $-wp\cdot h_\tau$ and
$-wz d_\tau$ cancel their counterparts in $wG_\tau$.
The Hessian blocks are those of the full evolving Hamiltonian;
they are not held fixed along the generating flow.

At $\tau=0$, the first-jet coefficients are
\[
    a_0=a+\ell\cdot\eta-\chi\zeta,\quad
    h_0=\omega+b+B\eta,\quad
    d_0=\beta+c+m\cdot\eta.
\]
Since $\Delta G$ is normal-affine, the initial Hessian blocks are
unchanged by the parameter update. The homological equations are
therefore equivalent to
\[
        a_0=(D_\omega+\beta)u,\quad h_0-\omega=D_\omega v+S_0Du+T_0u,\quad d_0-\beta-\delta\beta =D_\omega w+T_0\cdot Du+\mathfrak r_0u.
\]
Integrating \eqref{eq:jetconstant}, \eqref{eq:jetp}, and \eqref{eq:jetz} and using
these initial identities gives
\[
    \begin{aligned}
        a_1
        &=\int_0^1\bigl[
            v\cdot Da_\tau+w a_\tau
            -Du\cdot(h_\tau-\omega)-u(d_\tau-\beta)
            \bigr]\,d\tau,\\
        h_1-\omega
        &=\int_0^1\bigl[
            Dh_\tau\,v-Dv(h_\tau-\omega)
            -(S_\tau-S_0)Du-(T_\tau-T_0)u
            \bigr]\,d\tau,\\
        d_1-\beta-\delta\beta
        &=\int_0^1\bigl[
            v\cdot Dd_\tau-(h_\tau-\omega)\cdot Dw
            -(T_\tau-T_0)\cdot Du
            -(\mathfrak r_\tau-\mathfrak r_0)u
            \bigr]\,d\tau.
    \end{aligned}
\]
These formulas express the first-jet cancellation coefficientwise.
The operator estimate \eqref{eq:twoerrors} bounds the same
remainders without separate estimates for the evolving Hessian
blocks.

\subsection{The composition estimate}

The product rule gives, on any common angular strip,
\[
    \|fg\|_{C^2}\leqslant C\|f\|_{C^2}\|g\|_{C^2}.
\]
This also holds for matrix products and tensor contractions, with a
constant depending on dimension. We detail the derivative loss in
the composition estimate used in Lemma~\ref{lem:flow}.
Let $\mathcal Q(q)=q+h(q)$, with
$\|h\|_{s-12\sigma}\leqslant\sigma/2$ and
$\|h\|_{C^2(s-12\sigma)}\leqslant1$.
For a periodic coefficient $f$ of the old accumulated map,
\begin{align}
    D(f\circ\mathcal Q)
        &=(Df\circ\mathcal Q)(I+Dh),\notag
        \\
    D^2(f\circ\mathcal Q)[\nu_1,\nu_2]
        &=D^2f(\mathcal Q)[(I+Dh)\nu_1,(I+Dh)\nu_2]
          +Df(\mathcal Q)D^2h[\nu_1,\nu_2].
        \label{eq:chaintwo}
\end{align}
All line segments $q+th(q)$, $0\leqslant t\leqslant1$, lie in
$\T^n_{s-11\sigma}$. For $j=0,1,2$,
\[
    D^jf(q+h(q))-D^jf(q)
    =\int_0^1D^{j+1}f(q+th(q))[h(q)]\,dt,
\]
where $[h(q)]$ contracts one derivative slot.
Subtracting $D^2f(q)$ in \eqref{eq:chaintwo} and expanding the
two factors $I+Dh$ bounds the second-order difference by
\[
    C\bigl(
       \|D^3f\|\|h\|+\|D^2f\|\|Dh\|
       +\|D^2f\|\|Dh\|^2+\|Df\|\|D^2h\|
      \bigr).
\]
Here the norms of $f$ are taken on $s-11\sigma$ and those of
$h$ on $s-12\sigma$. Including the lower derivatives yields
\begin{equation}\label{eq:C2compositiondetail}
    \|f\circ\mathcal Q-f\|_{C^2(s-12\sigma)}
    \leqslant C\|f\|_{C^3(s-11\sigma)}
             \|h\|_{C^2(s-12\sigma)}.
\end{equation}
For $0<\sigma<1$, Cauchy's estimate gives
\begin{equation}\label{eq:C3old}
    \|f\|_{C^3(s-11\sigma)}
    \leqslant C\sigma^{-1}\|f\|_{C^2(s)}.
\end{equation}
Thus only one extra angular derivative of the old map is required.

Write the old base map as $Q=\mathrm{id}+k$ and the increment
triple as $(\mathcal Q,1+a,b)$. Formula \eqref{eq:composition}
gives
\begin{equation}\label{eq:tripledifferences}
    \begin{aligned}
 &Q_+-Q=h+k\circ\mathcal Q-k,\quad \rho_+-\rho=(\rho\circ\mathcal Q-\rho) +(\rho\circ\mathcal Q)a,\\
 &U_+-U=(U\circ\mathcal Q-U)+(\rho\circ\mathcal Q)b.
 \end{aligned}
\end{equation}
Apply \eqref{eq:C2compositiondetail}--\eqref{eq:C3old} to the
differences and the product estimate to the remaining terms.
When $Y_s\leqslant\delta$ and
$\|h\|_{C^2}+\|a\|_{C^2}+\|b\|_{C^2}\leqslant CA_F$,
this proves \eqref{eq:compositionbound} on $s-12\sigma$.
The identity part of $Q$ is treated by the explicit term $h$ in
\eqref{eq:tripledifferences}, rather than by a periodic-function
norm.

\subsection{Domains and fixed constants}\label{app:budgets}

The domains in one KAM step are summarized below. The normal
radius $r$ is fixed before the iteration.
\begin{center}
\small
\begin{tabular}{@{}p{0.62\linewidth}p{0.30\linewidth}@{}}
\hline
Object or estimate & Domain \\[2pt]
\hline
Current and parameter-updated Hamiltonians & $D_{s,8r}$ \\[3pt]
Scalar correction $u$ & $\T^n_{s-\sigma}$ \\[3pt]
Corrections $v,w$ & $\T^n_{s-3\sigma}$ \\[3pt]
Generating coefficients through two derivatives & $\T^n_{s-5\sigma}$ \\[3pt]
$\mathscr L_FG^\sharp$ & $D_{s-6\sigma,6r}$ \\[3pt]
$\mathscr L_F^2G^\sharp$ & $D_{s-7\sigma,4r}$ \\[3pt]
Pulled-back integral remainder & $D_{s-8\sigma,2r}$ \\[3pt]
New first-jet defect and accumulated triple & $\T^n_{s-12\sigma}$ \\[3pt]
Restored Hamiltonian for the next step & $D_{s-12\sigma,8r}$ \\[2pt]
\hline
\end{tabular}
\end{center}

The geometric bound $\bar\delta$ and $r$ are chosen first from
\eqref{eq:backgrounddomain}. They determine $B_0$ and a bound
for the normal Hessian. The fixed real-torus margin in
\eqref{eq:fixedmargin} then permits a choice of
$\delta\leqslant\bar\delta$ for which the average matrix is
invertible. Reducing this last bound does not increase $B_0$ or
the Hessian bound.

Let $C_{\rm par}$ and $C_{\rm lin}$ be the constants in the
parameter and generating-function estimates of
Lemma~\ref{lem:linear}, let $c_{\rm fl}$ and $C_{\rm comp}$
be the constants of Lemma~\ref{lem:flow}, and fix $C_Q$ so that
\[
    E_+\leqslant C_Q\sigma^{-(4\kappa+8)}E_s^2.
\]
The derivation of \eqref{eq:sharpstep} fixes $C_Q$ in terms of
these constants, $B_0$, and $r$, independently of the current
$\sigma$ and the iteration index. A sufficient choice is
\[
    0<c_*\leqslant\min\left\{
       1,\frac{c_{\rm fl}\min\{1,r\}}{C_{\rm lin}},
       \frac{\delta}{C_{\rm comp}C_{\rm lin}},
       \frac{1}{C_{\rm par}}
       \right\}.
\]
Indeed, $E_s\leqslant c_*\sigma^{\ell_*}$, with
$\ell_*=2\kappa+6$, implies
\[
    A_F\leqslant C_{\rm lin}c_*\sigma^3 \leqslant c_{\rm fl}\min\{\sigma,r,1\},\quad Y_{s-12\sigma}(\Phi_+) \leqslant\delta+C_{\rm comp}C_{\rm lin}c_*\sigma^2 \leqslant2\delta.
\]
The parameter increment is at most $C_{\rm par}c_*\leqslant1$,
so the updated parameters stay in the larger range used to bound
the original family. These bounds are obtained before extending
$G_+$ to $D_{s-12\sigma,8r}$ by its accumulated affine formula.
On the smaller flow domain this extension equals
$\mathscr P_1G^\sharp$, so its first jet satisfies the same
quadratic estimate.

One may then take
\[
    C_*\geqslant\max\{1,C_{\rm par},C_{\rm comp}C_{\rm lin},C_Q\}
\]
and choose the initial threshold by \eqref{eq:epsilonbudget}.
The distinction between the one-step bound $2\delta$ and the
iteration bound $\delta$ is used here: the former defines the
extended Hamiltonian, while summability of the increments keeps
all subsequent steps within the hypotheses of
Lemma~\ref{lem:step}.

\subsection{Tail estimates and limiting derivatives}

Use the sequences $\sigma_j$, $s^{(j)}$, $E_j$, and $d_j$
from Section~\ref{sec:iteration}, with
$\mathcal A=C_*\sigma_0^{-m_*}$,
$\mathcal B=2^{m_*}$, and
$d_j=\mathcal A\mathcal B^{j+1}E_j$.
For $d_0\leqslant1/2$, \eqref{eq:recursion} gives
$E_j\leqslant E_0\mathcal B^{-j}2^{-(2^j-1)}$.
For each fixed $L\geqslant0$, define
$\mathfrak d_j(L)=2^{(L-m_*)j}2^{-(2^j-1)}$. Then
\[
    \frac{\mathfrak d_{j+1}(L)}{\mathfrak d_j(L)}
    =2^{L-m_*}2^{-2^j}.
\]
If $\mathfrak q_J(L):=2^{L-m_*}2^{-2^J}<1$, all subsequent
ratios are bounded by $\mathfrak q_J(L)$, and hence
\[
    \sum_{j=J}^{\infty}\sigma_j^{-L}E_j
    \leqslant E_0\sigma_0^{-L}
          \frac{\mathfrak d_J(L)}{1-\mathfrak q_J(L)}.
\]
For $L=0$ and $L=\ell_*$, the inequality $L<m_*$ makes this
bound available from $J=0$. It is uniform in the real damping
parameter.

A tail bound can also be expressed in terms of the defect at a
constructed stage $J$. If $E_J=0$, every later correction vanishes.
Otherwise $0<d_J<1$, and repeated squaring gives
\[
    E_{J+k}\leqslant E_J\mathcal B^{-k}d_J^{\,2^k-1},
    \quad k\geqslant0.
\]
Since $\sigma_{J+k}=\sigma_J2^{-k}$ and $2^k-1\geqslant k$,
for $L=0$ or $L=\ell_*$ one obtains
\begin{equation}\label{eq:aposteriorsimple}
    \sum_{k=0}^{\infty}\sigma_{J+k}^{-L}E_{J+k}
    \leqslant\frac{\sigma_J^{-L}E_J}
                 {1-2^{L-m_*}d_J}.
\end{equation}
After multiplication by $C_*$, the cases $L=0$ and $L=\ell_*$
bound the remaining parameter and accumulated-triple corrections,
respectively. These are tail estimates for the iteration already
constructed; the inequality $d_J<1$ alone does not replace its
domain and smallness hypotheses.

Let $f_j$ denote any of $Q_j-\mathrm{id}$, $\rho_j-1$, or $U_j$.
For a fixed $s'<s_\infty$, put $\Delta=s_\infty-s'>0$.
Cauchy's estimate applied to the second derivatives gives, for
any fixed integer $k\geqslant3$,
\[
    \|f_{j+1}-f_j\|_{C^k(s')}
    \leqslant C(k,n)\bigl(1+\Delta^{-(k-2)}\bigr)
                  \|f_{j+1}-f_j\|_{C^2(s_\infty)}.
\]
Together with \eqref{eq:tripleincrement} and \eqref{eq:aposteriorsimple},
this bounds every fixed-order derivative of the remaining tail on
$s'$. In particular, it supplies the third angular derivatives
needed for the full-space $C^2$ estimate of \eqref{eq:affine}.
The same uniform tail control, combined with finite-stage
continuity, gives the parameter-continuity argument in
Section~\ref{sec:iteration}.

\end{document}